\documentclass{svmult}
\usepackage{bbm}
\usepackage{amsfonts}
\usepackage{mathrsfs}
\usepackage{amssymb}
\usepackage{cases}
\newenvironment{enumerate-roman}{\begin{enumerate}}{\end{enumerate}}
\usepackage{makeidx}         
\usepackage{graphicx}        
\usepackage{multicol}        
\usepackage[bottom]{footmisc}

\makeindex             

\newtheorem{thm}{Theorem}[section]
\newtheorem{lem}[thm]{Lemma}
\newtheorem{prop}[thm]{Proposition}
\newtheorem{cor}[thm]{Corollary}
\newtheorem{rmk}[thm]{Remark}

\newtheorem{defi}[thm]{Definition}

\newcommand{\R}{\mathbb{R}}

\begin{document}

\baselineskip14pt

\title*{SPDEs with Polynomial Growth Coefficients and Malliavin Calculus Method}
\author{Qi Zhang\inst{\ {\rm a}}, Huaizhong Zhao\inst{\ {\rm
b}}
}
\authorrunning{Q. Zhang and H.Z. Zhao}
\institute{$^{\rm a}$ School of Mathematical Sciences, Fudan
University, Shanghai, 200433, China\\
$^{\rm b}$ Department of Mathematical Sciences,
Loughborough University, Loughborough, LE11 3TU, UK\\
\texttt{Emails: qzh@fudan.edu.cn}; \texttt{H.Zhao@lboro.ac.uk}}
\maketitle
\newcounter{bean}
\begin{abstract}
In this paper we study the existence and uniqueness of the
$L_{\rho}^{2p}({\mathbb{R}^{d}};{\mathbb{R}^{1}})\times
L_{\rho}^2({\mathbb{R}^{d}};{\mathbb{R}^{d}})$ valued solution of
backward doubly stochastic differential equations with polynomial growth
coefficients using week convergence, equivalence of norm principle
and Wiener-Sobolev compactness arguments. Then we establish a new
probabilistic representation of the weak solutions of SPDEs with
polynomial growth coefficients through the solutions of the corresponding backward doubly stochastic differential equations
(BDSDEs). This probabilistic representation is then used to prove the existence of stationary solution of SPDEs on $\mathbb{R}^d$ via infinite
horizon BDSDEs.
The convergence of the solution of BDSDE to the solution of infinite horizon BDSDEs is shown to be equivalent to the convergence of the pull-back of the solutions of SPDEs. With this we obtain the stability of the stationary solutions as well.
\end{abstract}

\vskip5pt

\vskip5pt

\renewcommand{\theequation}{\arabic{section}.\arabic{equation}}

\section{Introduction}

In this paper, we study SPDEs on $\R^d$ with a polynomial growth coefficients of
the following type
\begin{eqnarray}\label{ch}
dv(t,x)=[\mathscr{L}v(t,x)+f\big(x,v(t,x)\big)]dt+g\big(x,v(t,x)\big)dB_t.
\end{eqnarray}
Here $\mathscr{L}$ is a second order differential operator given by
\begin{eqnarray}\label{h}
\mathscr{L}={1\over2}\sum_{i,j=1}^da_{ij}(x){{\partial^2}\over{\partial
x_i\partial x_j}}+\sum_{i=1}^db_i(x){\partial\over{\partial x_i}}.
\end{eqnarray}
Brownian motion $B$ is a $Q$-Wiener process with values in a
separable Hilbert space $U$. Denote the countable base of $U$ by
$\{e_i\}_{i=1}^{\infty}$, then $Q\in L(U)$ is a symmetric
nonnegative trace class operator such that $Qe_i=\lambda _ie_i$ and
$\sum \limits _{i=1}^{\infty}\lambda _i<\infty$. The coefficients
$f:\R^d\times\R^1\ni(x,v)\mapsto f(x,v)\in\R^1$ is a real-valued
function of polynomial growth $(p\geq2)$; $g:\R^d\times\R^1\ni(x,v)\mapsto
g(x,v)\in\mathcal{L}^2_{U_0}(\mathbb{R}^{1})$ is a Lipschitz
continuous function, where $U_0=Q^{1\over 2}(U)\subset U$ is a
separable Hilbert space with the norm $<u,v>_{U_0}=<Q^{-{1\over
2}}u, Q^{-{1\over 2}}v>_U$ and the complete orthonormal base
$\{\sqrt{\lambda_i}e_i\}_{i=1}^{\infty}$ and
$\mathcal{L}^2_{U_0}(\mathbb{R}^{1})$ is the space of all
Hilbert-Schmidt operators from $U_0$ to $\mathbb{R}^{1}$ with the
Hilbert-Schmidt norm.

One of the goals of this article is to study the probabilistic
representation to the solution of this equation via the
corresponding backward doubly stochastic differential equation
(BDSDE).
%
Apart from using the
Feynman-Kac formula, the solutions of the backward stochastic
differential equations (BSDEs), when they have some regularities
(continuous and differentiable in classical sense, or the solutions
and their weak derivatives exist in certain weighted
$L^q(\R^1,dtdP\rho^{-1}(x)dx)\times
L^2(\R^d,dtdP\rho^{-1}(x)dx)$ space) can give a probabilistic
representation of the corresponding PDEs. This has been achieved for classical solutions when the
coefficients are smooth enough in Pardoux and Peng \cite{pa-pe2} and
for viscosity solution when the coefficients are Lipschitz
continuous in \cite{pa-pe2} and for weak solutions in \cite{ba-ma},
\cite{ba-le}, \cite{zh-zh1}. When the coefficients are
non-Lipschitz, researchers have made some significant progress. In
\cite{le-sa}, Lepeltier and San Martin assumed that the
$\mathbb{R}^1$-valued function $f(r,x,y,z)$ satisfies the measurable
condition, the $y,z$ linear growth condition and the $y,z$
continuous condition, then they proved the existence of the solution
of the corresponding BSDEs. But the uniqueness of the solution
failed to be proved since the comparison theorem cannot be used
under the non-Lipschitz condition. In Zhang and Zhao \cite{zh-zh2}, we
proved the existence and uniqueness of the solution in the space
$S^2(\mathbb{R}^1,dtdP\rho^{-1}(x)dx)\times
L^2(\R^d,dtdP\rho^{-1}(x)dx)$ to BDSDEs under the monotonicity and
linear growth conditions, without assuming Lipschitz condition. We
also gave the probabilistic representation of the weak solutions of the
corresponding SPDEs. Along the line of the viscosity solution, in
\cite{ko}, Kobylanski was able to solve the BSDE when the
coefficients $f(y,z)$ is of quadratic growth in $z$, with the help
of Hopf-Cole transformation. In \cite{pa1}, Pardoux used weak
convergence in a finite dimensional space to study the viscosity
solution of the PDEs and the corresponding BSDEs when $f(y,z)$ is of
polynomial growth in $y$. As to the weak solutions of PDEs with
polynomial growth coefficients and the corresponding BSDEs, the
existing methods in BSDE were not adequate to solving the equation.
In \cite{zh-zh3}, we developed a new method to BSDEs using Alaoglu
weak convergence theorem and Rellich-Kondrachov compactness
embedding theorem to get a strongly convergent subsequence in the
space $L^2(\mathbb{R}^1,dt\rho^{-1}(x)dx)$, and then using
equivalence of norm principle to get the convergent subsequence for
BSDEs. We therefore established the correspondence of the solutions
of BSDEs and the weak solutions of such PDEs.

To solve the BDSDEs corresponding to SPDE (\ref{ch}) when $f$ is of
polynomial growth in $y$, we can still use the Alaoglu weak
convergence argument. But the key to make it work is to find a
strong convergent subsequence. In the deterministic case, we use the
estimate for Sobolev norm of the solutions of the sequence of BSDEs
to get a strong convergent subsequence in $L^2(\rho^{-1}(x)dx)$. But
this method does not work for the BDSDEs as the subsequence choice
may depend on $\omega\in\Omega$. In this paper, we will develop a
method using Wiener-Sobolev compactness argument to tackle the
compactness problem for BDSDEs. First we estimate the Sobolev norm
of the solution as well as the Malliavin derivative to get
convergent sequence in $L^2(dtdP_B\rho^{-1}(x)dx)$, then from
equivalence of norm principle to pass the compactness to the
solutions of BDSDEs in $L^2(dtdP_BdP_W\rho^{-1}(x)dx)$. The
Wiener-Sobolev compact embedding theorem is a powerful tool in
proving the relatively compactness of a random field. The random
version (independent of time and spatial variables) was obtained in
Da Prato, Malliavin and Nualart \cite{da-ma}, Peszat \cite{pes}. This was
extended later by Bally and Saussereau \cite{ba-sa}. This has been extended to the space $C\big([0,T];L^2(dPdx)\big)$ and
applied to study the existence of an infinite horizon stochastic
integral equation arising in the study of random period solution in
Feng, Zhao and Zhou \cite{fe-zh-zh}, Feng and Zhao \cite{fz1}.

Our motivation to study the probabilistic representation is to use
it to study the dynamics of the random dynamical system generated by the
SPDEs, although the probabilistic representation and weak solution of BDSDEs with polynomial growth coefficients were interesting problems themselves and have their own
interest. Stationary solution is one of the central concepts in the
study of the long term behaviour of the stochastic dynamical systems
generated by SPDEs. It is a pathwise equilibrium which is invariant,
over time, along its measurable and $P$-preserving metric dynamical
system $\theta_t:\Omega\longrightarrow\Omega$. In deterministic
case, it gives the solution of elliptic equation. Due to the nature
of the noise that is pumped to the system constantly, the stationary
solution is random and changes along times. Therefore the study is
rather difficult and no universal method, but important in order to understand the equilibrium and long time behaviour of stochastic systems. There are many works in
the literature on the local behaviour of the solutions near a
stationary solution, if exists (e.g. Arnold \cite{ar}, Duan, Lu and Schmalfuss \cite{du-lu-sc1}, Mohammed, Zhang and Zhao \cite{mo-zh-zh}, Lian and Lu \cite{li-lu} to name but a few). So the existence is a key to
understand complexity of many random dynamical systems. Although there is no universal method applicable to generic problems, researchers have obtained many results to a variety of SPDEs e.g. Sinai \cite{si1}, \cite{si2}, Mattingly \cite{ma},  E, Khanin, Mazel and Sinai \cite{kh-ma-si}, Caraballo, Kloeden and Schmalfuss \cite{kloeden}, Zhang and Zhao \cite{zh-zh1}, \cite{zh-zh2}. The cases of non-dissipative stochastic differential equations
and SPDEs with additive noise has been obtained in Feng, Zhao and Zhou \cite{fe-zh-zh}, Feng and Zhao \cite{fz1}. In
applications, stationary solutions also appear in many other real world
problems, e.g. in the interpolation of data and image processing,
the stationary solution of the stochastic parabolic infinity
Laplacian equation gives the final restored image of the image
processing in a random model (Wei and Zhao \cite{we-zh}). Note
that in the corresponding deterministic model, the elliptic infinity Laplacian equation gives
the final restored image as the limit of the solution of the
infinity Laplacian equation (Caselles, Morel and Sbert
\cite{ca-mo-sb}). In this paper, we will solve the infinite horizon
BDSDEs with the polynomial growth coefficients and therefore obtain
the stationary solutions of SPDEs (\ref{ch}). We will also prove that the convergence of the solution of BDSDE to the solution of infinite horizon BDSDEs is equivalent to the convergence of the pull-back of the solutions of SPDEs. Therefore, we obtain the convergence of the solution with a class of initial condition $h$ converges to the stationary solution.


\section{Preliminaries and definitions}
\setcounter{equation}{0}

We will study the weak solutions of the SPDE (\ref{ch}) and the
corresponding BDSDE in a Hilbert space ($\rho$-weighted $L^2(dx)$
space). Utilizing this correspondence, we will give the
probabilistic representation of the weak solution of SPDE (\ref{ch})
on finite horizon with a given initial value and find the stationary
solution of SPDE (\ref{ch}).

For this purpose, we study first backward SPDE. Let $(\Omega,\mathscr{F},P)$ be a probability space and $\hat{B}$, $W$ be
mutually independent Brownian motions in $U$ and $\mathbb{R}^d$,
respectively. In Section 6, we choose $\hat{B}$ to be time reversal Brownian motion of $B$ so establish connection with forward SPDEs, especially its stationary solution. Here we consider general Brownian motion $\hat{B}$.



We first consider the following backward SPDE:
\begin{eqnarray}\label{bz}
u(t,x)&=&h(x)+\int_{t}^{T}[\mathscr{L}u(s,x)+f\big(s,x,u(s,x)\big)]ds\nonumber\\
&&-\int_{t}^{T}g\big(s,x,u(s,x)\big)d^\dagger \hat{B}_s,\ \ \ \
0\leq t\leq T.
\end{eqnarray}
Here $\mathscr{L}$ is given by (\ref{h}) with $b:\mathbb{R}^d\longrightarrow\mathbb{R}^d$,
$a=\sigma\sigma^*:\mathbb{R}^d\longrightarrow\mathbb{R}^{d\times d}$. Assume
$h:\mathbb{R}^d\longrightarrow\mathbb{R}^1$,
$f:[0,T]\times\mathbb{R}^d\times\mathbb{R}^1\longrightarrow\mathbb{R}^1$
and
$g:[0,T]\times\mathbb{R}^d\times\mathbb{R}^1\longrightarrow\mathcal{L}^2_{U_0}(\mathbb{R}^{1})$
are measurable. The stochastic integral $\int_{t}^{T}g(s,x,u(s,x))d^\dagger\hat{B}_s$ is a backward stochastic integral which will be made clear later.

Denote by $L_{\rho}^2({\mathbb{R}^{d}};{\mathbb{R}^{1}})$ the space of measurable functions $l:\mathbb{R}^d\longrightarrow\mathbb{R}^1$ such that $\int_{\mathbb{R}^d}l^2(x)\rho^{-1}(x)dx<\infty$. Define the inner product
\begin{eqnarray*}
\langle
l_1,l_2\rangle=\int_{\mathbb{R}^d}l_1(x)l_2(x)\rho^{-1}(x)dx,
\end{eqnarray*}
then $L_{\rho}^2({\mathbb{R}^{d}};{\mathbb{R}^{1}})$ is a Hilbert space.
Here $\rho(x)=(1+|x|)^q$, $q>d+{8p}$, is a weight function
and $p$ is given in Condition (H.1). Similarly, denote by
$L_{\rho}^{k}({\mathbb{R}^{d}};{\mathbb{R}^{1}})$, $k\geq2$, the
weighted $L^k$ space with the norm
$||l||_{L_{\rho}^{k}(R^d)}=(\int_{\mathbb{R}^{d}}l^{k}(x)\rho^{-1}(x)dx)^{1\over
{k}}$. It is easy to see that
$\rho(x):\mathbb{R}^d\longrightarrow\mathbb{R}^1$ is a continuous
positive function satisfying
$\int_{\mathbb{R}^{d}}|x|^{{8p}}\rho^{-1}(x)dx<\infty$. We can
consider more general $\rho(x)$ as in \cite{ba-ma} and all the
results of this paper still hold. But this is not the purpose of this paper. Note that due to the polynomial growth of
$f$, we need $u(t,\cdot)\in
L_{\rho}^{2p}({\mathbb{R}^{d}};{\mathbb{R}^{1}})$ for $t\in[0,T]$.

Now define $X_s^{t,x}$ to be the solution of the following stochastic differential equations for any given $t\geq0$ and $x\in\mathbb{R}^d$:
\begin{numcases}{}\label{a}
X_{s}^{t,x}=x+\int_{t}^{s}b(X_{r}^{t,x})dr+\int_{t}^{s}\sigma(X_{r}^{t,x})dW_r,\ \ \ s\geq t,\nonumber\\
X_s^{t,x}=x,\ \ \ 0\leq s<t.
\end{numcases}
The BDSDE associated with SPDE (\ref{bz}) is
\begin{eqnarray}\label{spdespgrowth2}
Y_{s}^{t,x}&=&h(X_{T}^{t,x})+\int_{s}^{T}f(r,X_{r}^{t,x},Y_{r}^{t,x})dr\nonumber\\
&&-\int_{s}^{T}g(r,X_{r}^{t,x},Y_{r}^{t,x})d^\dagger{\hat{B}}_r-\int_{s}^{T}\langle
Z_{r}^{t,x},dW_r\rangle,\ \ \ 0\leq t\leq s\leq T.
\end{eqnarray}
It is well known that $\hat B$ has the following expansion
(\cite{pr-za1}): for each $r$,
\begin{eqnarray}\label{bmexpansion}
\hat B_r=\sum\limits _{j=1}^{\infty} \sqrt {\lambda _j}\hat \beta
_j(r)e_j,
\end{eqnarray}
where
\begin{eqnarray*}
\hat \beta _j(r)={1\over \sqrt{\lambda _j}}<\hat B_r,e_j>_U, \ \
j=1,2,\cdots
\end{eqnarray*}
are mutually independent real-valued Brownian motions on $(\Omega,
{\mathscr{F}}, P)$ and the series (\ref{bmexpansion}) is convergent
in $L^2(\Omega, {\mathscr{F}}, P)$. Set $g_j\triangleq
g\sqrt{\lambda_j}e_j:\mathbb{R}^{d}\times\mathbb{R}^1\times\mathbb{R}^{d}{\longrightarrow{\mathbb{R}^1}}$,
then BDSDE (\ref{spdespgrowth2}) is equivalent to
\begin{eqnarray}\label{cm}
Y_{s}^{t,x}&=&h(X_{T}^{t,x})+\int_{s}^{T}f(r,X_{r}^{t,x},Y_{r}^{t,x})dr\nonumber\\
&&-\sum_{j=1}^{\infty}\int_{s}^{T}g_j(r,X_{r}^{t,x},Y_{r}^{t,x})d^\dagger{\hat{\beta}}_j(r)-\int_{s}^{T}\langle
Z_{r}^{t,x},dW_r\rangle,\ \ \ 0\leq t\leq s\leq T.
\end{eqnarray}

For the convenience of readers, we need to recall the definitions of
weak solutions of SPDEs and the
$L_{\rho}^2({\mathbb{R}^{d}};{\mathbb{R}^{1}})\times
L_{\rho}^2({\mathbb{R}^{d}};{\mathbb{R}^{d}})$ valued solutions of
BDSDEs. Denote by ${\cal N}$ the class of $P$-null sets of
${\mathscr{F}}$ and let
\begin{eqnarray*}
\mathscr{F}_{s,T}\triangleq{\mathscr{F}_{s,T}^{\hat{B}}}\bigvee
\mathscr{F}_{t,s}^W,\ \ \ {\rm for}\ 0\leq t\leq s\leq T,\ \ \ \ \ \
\mathscr{F}_s\triangleq{\mathscr{F}_{s,\infty}^{\hat{B}}}\bigvee
\mathscr{F}_{t,s}^W,\ \ \ {\rm for}\ 0\leq t\leq s,
\end{eqnarray*}
where for any process $(\eta_s)_{s\geq0}$,
$\mathscr{F}_{t,s}^\eta=\sigma\{\eta_r-\eta_t$; ${0\leq t\leq r\leq
s}\}\bigvee{\cal N}$, $\mathscr{F}_{s,\infty}^{\eta}=\bigvee_{T\geq
s}{\mathscr{F}_{s,T}^\eta}$.

First recall
\begin{defi}\label{qi00} (Definitions 2.1, \cite{zh-zh1})
Let $\mathbb{S}$ be a separable Banach space with norm $\|\cdot\|_\mathbb{S}$
and Borel $\sigma$-field $\mathscr{S}$ and $q\geq2$, $K>0$. We
denote by $M^{q,-K}([t,\infty);\mathbb{S})$ the set of
$\mathscr{B}([t,\infty))\otimes\mathscr{F}/\mathscr{S}$ measurable
random processes $\{\phi(s)\}_{s\geq t}$ with values in $\mathbb{S}$
satisfying
\begin{enumerate-roman}
\item $\phi(s):\Omega\longrightarrow\mathbb{S}$ is $\mathscr{F}_{s}$ measurable for $s\geq t$;
\item $E[\int_{t}^{\infty}{\rm e}^{-Ks}\|\phi(s)\|_\mathbb{S}^qds]<\infty$.
\end{enumerate-roman}
Also we denote by $S^{q,-K}([t,\infty);\mathbb{S})$ the set of
$\mathscr{B}([t,\infty))\otimes\mathscr{F}/\mathscr{S}$ measurable
random processes $\{\psi(s)\}_{s\geq t}$ with values in $\mathbb{S}$
satisfying
\begin{enumerate-roman}
\item $\psi(s):\Omega\longrightarrow\mathbb{S}$ is
$\mathscr{F}_{s}$ measurable for $s\geq t$ and $\psi(\cdot,\omega)$
is a.s. continuous;
\item $E[\sup_{s\geq t}{\rm
e}^{-Ks}\|\psi(s)\|_\mathbb{S}^q]<\infty$.
\end{enumerate-roman}
\end{defi}
If we replace time interval $[t,\infty)$ by $[t,T]$ in the above
definition, we denote the spaces by $M^{q,0}([t,T];\mathbb{S})$ and
$S^{q,0}([t,T];\mathbb{S})$, respectively. Note that here ${\rm
e}^{-Ks}$ does not play role as $T$ is finite, so we can always take
$K=0$.

For the backward stochastic integral, let
$\{g(s)\}_{s\geq0}$ be a stochastic process with values
in ${\mathcal{L}^2_{U_0}(H)}$ such that $g(s)$ is
$\mathscr{F}_s$ measurable for any $s\geq0$ and locally square
integrable, i.e. for any $0\leq a\leq b<\infty$,
$\int_{a}^{b}\|g(s)\|_{\mathcal{L}^2_{U_0}(H)}^2ds<\infty$ a.s.
Since ${\mathscr{F}}_s$ is a backward filtration with respect to
$\hat B$, so from the one-dimensional backward It$\hat {\rm o}$'s
integral and relation with forward integral, for $0\leq T\leq
T^{\prime}$, we have
\begin{eqnarray*}
\int _t^T\sqrt{\lambda_j}<g(s)e_j,f_k>d^\dagger\hat \beta_j(s)=-\int
_{T^{\prime}-T}^{T^{\prime}-t}\sqrt{\lambda_j}<g(T^{\prime}-s)e_j,f_k>
d \beta_j(s), \ \ j,k=1,2,\cdots
\end{eqnarray*}
where $\beta _j(s)=\hat \beta _j(T^{\prime}-s)-\hat \beta
_j(T^{\prime})$, $j=1,2,\cdots$, and so $B_s=\hat
B_{T^{\prime}-s}-\hat B_{T^{\prime}}$. Here $\{f_k\}$ is the
complete orthonormal basis in $H$. From approximation theorem of the
stochastic integral in a Hilbert space (cf. \cite{pr-za1}),
we have
\begin{eqnarray*}
\int _{T^{\prime}-T}^{T^{\prime}-t}g(T^{\prime}-s)d
B_s=\sum_{j,k=1}^{\infty}\int
_{T^{\prime}-T}^{T^{\prime}-t}\sqrt{\lambda_j}<g(T^{\prime}-s)e_j,f_k>
d\beta_j(s)f_k.
\end{eqnarray*}
Similarly we also have
\begin{eqnarray*}
\int _t^Tg(s)d^\dagger\hat B_s=\sum_{j,k=1}^{\infty} \int
_t^T\sqrt{\lambda_j}<g(s)e_j,f_k>d^\dagger\hat \beta_j(s)f_k.
\end{eqnarray*}
It turns out that
\begin{eqnarray*}\label{bmrelation}
\int_{t}^{T}g(s)d^\dagger{\hat{B}}_s=-\int_{T'-T}^{T'-t}
g({T'-s})d{B}_s\ \ \ \ {\rm a.s.}
\end{eqnarray*}

\begin{defi}\label{12}
A function $u$ is called a weak solution of SPDE (\ref{bz}) if
$(u,\sigma^*\nabla u)\in
L^{2p}([0,T];L_{\rho}^{2p}\\({\mathbb{R}^{d}};{\mathbb{R}^{1}}))\times
L^{2}([0,T];L_{\rho}^2({\mathbb{R}^{d}};{\mathbb{R}^{d}}))$ and for
an arbitrary $\varphi\in C_c^{\infty}(\mathbb{R}^d;\mathbb{R}^1)$,
\begin{eqnarray}\label{ae}
&&\int_{\mathbb{R}^{d}}u(t,x)\varphi(x)dx-\int_{\mathbb{R}^{d}}h(x)\varphi(x)dx-{1\over2}\int_{t}^{T}\int_{\mathbb{R}^{d}}(\sigma^*\nabla u)(s,x)(\sigma^*\nabla\varphi)(x)dxds\nonumber\\
&&-\int_{t}^{T}\int_{\mathbb{R}^{d}}u(s,x)div\big((b-\tilde{A})\varphi\big)(x)dxds\\
&=&\int_{t}^{T}\int_{\mathbb{R}^{d}}f\big(s,x,u(s,x)\big)\varphi(x)dxds-\int_{t}^{T}\int_{\mathbb{R}^{d}}g\big(s,x,u(s,x)\big)\varphi(x)dxd^\dagger\hat{B}_s,\ \ t\in[0,T].\nonumber
\end{eqnarray}
Here $\tilde{A}_j\triangleq{1\over2}\sum_{i=1}^d{\partial
a_{ij}(x)\over\partial x_i}$, and
$\tilde{A}=(\tilde{A}_1,\tilde{A}_2,\cdot\cdot\cdot,\tilde{A}_d)^*$.
\end{defi}
\begin{rmk}
The weak solution of forward SPDE (\ref{ch}) with initial value
$v(0,\cdot)$ can be defined similarly. We also represent it in a
form like $v(t,\cdot,v(0,\cdot))$ to emphasize its dependence on the
initial value $v(0,\cdot)$, when it is necessary.
\end{rmk}

We then give the definition for the
$L_{\rho}^{2p}({\mathbb{R}^{d}};{\mathbb{R}^{1}})\times
L_{\rho}^2({\mathbb{R}^{d}};{\mathbb{R}^{d}})$ valued solution of
BDSDE (\ref{spdespgrowth2}).
\begin{defi}\label{qi051}
A pair of processes $(Y_{s}^{t,x},Z_{s}^{t,x})$ is called a solution
of BDSDE (\ref{spdespgrowth2}) if
$(Y_{\cdot}^{t,\cdot},Z_{\cdot}^{t,\cdot})\in
S^{2p,0}([t,T];L_{\rho}^{2p}({\mathbb{R}^{d}};{\mathbb{R}^{1}}))\times
M^{2,0}([t,T];L_{\rho}^2({\mathbb{R}^{d}};{\mathbb{R}^{d}}))$ and
$(Y_s^{t,x},Z_s^{t,x})$ satisfies
(\ref{spdespgrowth2}) for a.e. $x\in{\mathbb{R}^{d}}$ a.s.
\end{defi}
\begin{rmk}\label{16}
Due to the density of $C_c^{0}(\mathbb{R}^d;\mathbb{R}^1)$ in
$L_{\rho}^2({\mathbb{R}^{d}};{\mathbb{R}^{1}})$, we have, for a.e. $x$,
$(Y_s^{t,x,n},Z_s^{t,x,n})$ satisfies (\ref{spdespgrowth2}) is
equivalent to for an arbitrary $\varphi\in
C_c^{0}(\mathbb{R}^d;\mathbb{R}^1)$, $({Y}_s^{t,x},{Z}_s^{t,x})$
satisfies
\begin{eqnarray*}\label{zhang682}
\int_{\mathbb{R}^{d}}Y_s^{t,x}\varphi(x)dx&=&\int_{s}^{T}\int_{\mathbb{R}^{d}}f(r,X_{r}^{t,x},Y_{r}^{t,x})\varphi(x)dxdr-\int_{s}^{T}\int_{\mathbb{R}^{d}}g(r,X_{r}^{t,x},Y_{r}^{t,x})\varphi(x)dxd^\dagger{\hat{B}}_r\nonumber\\
&&-\int_{s}^{T}\langle
\int_{\mathbb{R}^{d}}Z_{r}^{t,x}\varphi(x)dx,dW_r\rangle\ \ \ P-{\rm
a.s.}
\end{eqnarray*}
\end{rmk}

To find the stationary solution of SPDE (\ref{ch}), we need to
consider its corresponding infinite horizon BDSDE:
\begin{eqnarray}\label{bp}
{\rm e}^{-Ks}Y_{s}^{t,x}&=&\int_{s}^{\infty}{\rm e}^{-Kr}f(X_{r}^{t,x},Y_{r}^{t,x})dr+\int_{s}^{\infty}K{\rm e}^{-Kr}Y_{r}^{t,x}dr\nonumber\\
&&-\int_{s}^{\infty}{\rm
e}^{-Kr}g(X_{r}^{t,x},Y_{r}^{t,x})d^\dagger{\hat{B}}_r-\int_{s}^{\infty}{\rm
e}^{-Kr}\langle Z_{r}^{t,x},dW_r\rangle.
\end{eqnarray}
For the existence and uniqueness of the solution, we can study a more
general form of the above infinite horizon BDSDE with time variable
dependent coefficients and $X_s^{t,x}$ is still the flow generated
by (\ref{a}):
\begin{eqnarray}\label{qi30}
{\rm e}^{-Ks}Y_{s}^{t,x}&=&\int_{s}^{\infty}{\rm e}^{-Kr}f(r,X_{r}^{t,x},Y_{r}^{t,x})dr+\int_{s}^{\infty}K{\rm e}^{-Kr}Y_{r}^{t,x}dr\nonumber\\
&&-\int_{s}^{\infty}{\rm
e}^{-Kr}g(r,X_{r}^{t,x},Y_{r}^{t,x})d^\dagger{\hat{B}}_r-\int_{s}^{\infty}{\rm
e}^{-Kr}\langle Z_{r}^{t,x},dW_r\rangle.
\end{eqnarray}
Here
$f:[0,\infty)\times\mathbb{R}^{d}\times\mathbb{R}^1{\longrightarrow{\mathbb{R}^1}}$,
$g:[0,\infty)\times\mathbb{R}^{d}\times\mathbb{R}^1\longrightarrow\mathcal{L}^2_{U_0}(\mathbb{R}^{1})$.
\begin{defi}\label{4}
A pair of processes $(Y_{s}^{t,x},Z_{s}^{t,x})$ is called a solution
of BDSDE (\ref{qi30}) if
$(Y_{\cdot}^{t,\cdot},Z_{\cdot}^{t,\cdot})\in S^{2p,-K}\bigcap
M^{2p,-K}([t,\infty);L_{\rho}^{2p}({\mathbb{R}^{d}};{\mathbb{R}^{1}}))\times
M^{2,-K}([t,\infty);L_{\rho}^2({\mathbb{R}^{d}};{\mathbb{R}^{d}}))$
and $(Y_s^{t,x},Z_s^{t,x})$
satisfies (\ref{qi30}) for a.e. $x\in{\mathbb{R}^{d}}$ a.s.
\end{defi}
For arbitrary $T\geq t$, the general connection between the solution
of BDSDE (\ref{bp}) and stationary solution of SPDE (\ref{ch}) was established in Zhang and Zhao \cite{zh-zh2}. Such a connection is proved independent of $T$ at which $B$ is reversed to $\hat{B}$ as shown in
\cite{zh-zh1}, \cite{zh-zh2}. Therefore to find the solution of the infinite horizon BDSDE (\ref{bp}) is the key to construct the stationary solution
of SPDE (\ref{ch}).

For $k\geq0$, denote by $C^k_{l,b}$ the set
of $C^k$-functions whose partial derivatives up to $k$th order are bounded,
but the functions themselves need not be bounded, otherwise if the
functions themselves are also bounded, we denote this subspace of
$C^k_{l,b}$ by $C_{b}^k$. The following generalized equivalence of norm principle is an
extension of equivalence of norm principle given in \cite{ku1},
\cite{ba-le}, \cite{ba-ma} to the case when $\varphi$ and $\Psi$ are random.
\begin{lem}\label{qi045}(generalized equivalence
of norm principle \cite{zh-zh1}) Let $X$ be the diffusion process
defined in (\ref{a}) with $b\in
C_{b,l}^2(\mathbb{R}^{d};\mathbb{R}^{d})$, $\sigma\in
C_{b}^3(\mathbb{R}^{d};\mathbb{R}^{d}\times\mathbb{R}^{d})$. If
$s\in[t,T]$,
$\varphi:\Omega\times\mathbb{R}^d\longrightarrow\mathbb{R}^1$ is
independent of the $\sigma$-field $\mathscr{F}^W_{t,s}$ and $\varphi\rho^{-1}\in L^1(\Omega\times\mathbb{R}^{d})$, then
there exist two constants $c>0$ and $C>0$ s.t.
\begin{eqnarray*}
cE[\int_{\mathbb{R}^{d}}|\varphi(x)|\rho^{-1}(x)dx]\leq
E[\int_{\mathbb{R}^{d}}|\varphi(X_{s}^{t,x})|\rho^{-1}(x)dx]\leq
CE[\int_{\mathbb{R}^{d}}|\varphi(x)|\rho^{-1}(x)dx].
\end{eqnarray*}
Moreover if
$\Psi:\Omega\times[t,T]\times\mathbb{R}^d\longrightarrow\mathbb{R}^1$,
$\Psi(s,\cdot)$ is independent of $\mathscr{F}^W_{t,s}$ and
$\Psi\rho^{-1}\in L^1(\Omega\times[t,T]\times\mathbb{R}^{d})$, then
\begin{eqnarray*}
&&cE[\int_{t}^{T}\int_{\mathbb{R}^{d}}|\Psi(s,x)|\rho^{-1}(x)dxds]\leq E[\int_{t}^{T}\int_{\mathbb{R}^{d}}|\Psi(s,X_{s}^{t,x})|\rho^{-1}(x)dxds]\\
&\leq&CE[\int_{t}^{T}\int_{\mathbb{R}^{d}}|\Psi(s,x)|\rho^{-1}(x)dxds].
\end{eqnarray*}
\end{lem}

In the process of obtaining the stationary solution of SPDE
(\ref{ch}), the proof of existence and uniqueness of solution to
BDSDE (\ref{spdespgrowth2}) is a crucial and challenging
step. For this, we will start from studying BDSDE
(\ref{spdespgrowth2}) with finite dimensional noise in next two
sections:
\begin{eqnarray}\label{cn}
Y_{s}^{t,x,N}&=&h(X_{T}^{t,x})+\int_{s}^{T}f(r,X_{r}^{t,x},Y_{r}^{t,x,N})dr\nonumber\\
&&-\sum_{j=1}^{N}\int_{s}^{T}g_j(r,X_{r}^{t,x},Y_{r}^{t,x,N})d^\dagger{\hat{\beta}}_j(r)-\int_{s}^{T}\langle
Z_{r}^{t,x,n N},dW_r\rangle,\ \ \ 0\leq t\leq s\leq T.\ \ \ \ \ \
\end{eqnarray}
Actually, when $N$ tends to infinity, the solution of BDSDE (\ref{cn}) converges to the solution of BDSDE (\ref{cm}) which is equivalent to BDSDE
(\ref{spdespgrowth2}).

The rest of this paper is organized as follows. In Sections 3, we consider approximating BDSDE with Lipschitz coefficients and then use Alaoglu lemma to get a weakly convergent subsequence. We further utilize the equivalence of norm
principle and Malliavin derivatives to get a strongly convergent
subsequence and prove the existence and uniqueness of the solution to
BDSDE (\ref{cn}) in Section 4. In Section 5 we prove that
BDSDE (\ref{spdespgrowth2}), its corresponding backward SPDE
(\ref{bz}) and hence, by variable changes, SPDE (\ref{ch}), have a unique weak solution. The stationary
properties of solutions of BDSDE (\ref{bp}) and SPDE (\ref{ch}) are
shown in Section 6 after proving the existence and uniqueness for
solution of infinite horizon BDSDE (\ref{qi30}).


\section{The weak convergence}
\setcounter{equation}{0}

In this section, we consider BDSDE (\ref{cn}) which can be written as  
\begin{eqnarray}\label{co}
Y_{s}^{t,x}&=&h(X_{T}^{t,x})+\int_{s}^{T}f(r,X_{r}^{t,x},Y_{r}^{t,x})dr
-\int_{s}^{T}\langle
g(r,X_{r}^{t,x},Y_{r}^{t,x}),d^\dagger{\hat{B}}_r\rangle-\int_{s}^{T}\langle
Z_{r}^{t,x},dW_r\rangle,
\end{eqnarray}
for $0\leq t\leq s\leq T$.
Here $\hat{B}=(\beta_1,\beta_2,\cdot\cdot\cdot,\beta_l)$ is a $l$-dimensional Brownian motion.
We assume 
%
\begin{description}
\item[(H.1).] There exists a constant $p\geq2$ and a function ${f}_0:[0,T]\times\mathbb{R}^d\longrightarrow\mathbb{R}^1$ with\\ $\int_0^T\int_{\mathbb{R}^d}|{f}_0(s,x)|^{{8p}}\rho^{-1}(x)dxds<\infty$ s.t. for any $s\in[0,T]$, $x\in{\mathbb{R}^{d}}$, $y\in
{\mathbb{R}^{1}}$, $|f(s,x,y)|\leq
L(|{f}_0(s,x)|+|y|^p)$ and $|\partial_yf(s,x,y)|\leq L(1+|y|^{p-1})$.
\item[(H.2).] 
For the above $p\geq2$, and for any $s,s_1,s_2\in[0,T]$, $x,x_1,x_2\in{\mathbb{R}^{d}}$, $y,y_1,y_2\in
{\mathbb{R}^{1}}$,
\begin{eqnarray*}
&&|f(s,x_1,y)-f(s,x_2,y)|\leq L(1+|y|^{p})|x_1-x_2|,\\
&&|g(s_1,x_1,y_1)-g(s_2,x_2,y_2)|\leq
L(|s_1-s_2|+|x_1-x_2|+|y_1-y_2|).
\end{eqnarray*}
Moreover, $\partial_yf$, $\partial_yg$ exist and satisfy
\begin{eqnarray*}
&&|\partial_yf(s,x_1,y)-\partial_yf(s,x_2,y)|\leq L(1+|y|^{p-1})|x_1-x_2|,\\
&&|\partial_yf(s,x,y_1)-\partial_yf(s,x,y_2)|\leq L(1+|y_1|^{p-2}+|y_2|^{p-2})|y_1-y_2|,\\
&&|\partial_yg(s,x,y)|\leq L,\\
&&|\partial_yg(s,x_1,y_1)-\partial_yg(s,x_2,y_2)|\leq L(|x_1-x_2|+|y_1-y_2|).
\end{eqnarray*}
\item[(H.3).] There exists a
constant $\mu\in\mathbb{R}^{1}$ s.t. for any $s\in[0,T]$, $x\in{\mathbb{R}^{d}}$, $y_1,
y_2\in {\mathbb{R}^{1}}$,
\begin{eqnarray*}
(y_1-y_2)\big(f(s,x,y_1)-f(s,x,y_2)\big)\leq\mu{|y_1-y_2|}^2.
\end{eqnarray*}
\item[(H.4).] For the above $p\geq 2$, $\int_{\mathbb{R}^d}|h(x)|^{{8p}}\rho^{-1}(x)dx<\infty$ 
     and $E[\int_{\mathbb{R}^d}|h(X_{T}^{t',x})-h(X_{T}^{t,x})|^q\rho^{-1}(x)dx]\leq L|t'-t|^{q\over2}$ for any $2\leq q\leq{8p}$ and $X$ defined by (\ref{a}).
\item[(H.5).] The diffusion coefficients $b\in C_{l,b}^2(\mathbb{R}^{d};\mathbb{R}^{d})$, $\sigma\in
C_{b}^3(\mathbb{R}^{d};\mathbb{R}^{d}\times\mathbb{R}^{d})$. 
\item[(H.6).] The matrix $\sigma(x)$ is uniformly elliptic, i.e. there exists a constant $\varepsilon>0$ s.t. $\sigma\sigma^*(x)\geq \varepsilon I_d$.
\end{description}
\begin{rmk}
(i) In (H.1) and (H.4), the power ${8p}$ is determined by the estimates
in Theorem \ref{2}.\\
(ii) Condition (H.4) is weaker than Lipschitz condition of $h$. This assumption is not for the sake to weaken the Lipschitz condition of $h$. When we consider the stationary solution in Section 6, we cannot prove that the stationary solution is Lipschitz in $x$, but we can prove that it satisfies Condition (H.4) in Lemma \ref{29}.\\
(iii) The smooth condition (H.5) guarantees the existence of the flow of diffeomorphism. This is essential in the equivalence of norm principle (Lemma \ref{qi045}), which together with the uniform ellipticity condition (H.6) implies the equivalence of the norm between the solution of SPDE and the solution of BDSDE. See (\ref{dc}) in Section 4.
\end{rmk}
From (H.2) and the fact that
$\int_{\mathbb{R}^{d}}|x|^{{8p}}\rho^{-1}(x)dx<\infty$, we have
\begin{eqnarray}\label{au}
\sup_{s\in[0,T]}\int_{\mathbb{R}^d}|g(s,x,0)|^{{8p}}\rho^{-1}(x)dx<\infty.
\end{eqnarray}

It is easy to see that
$(Y_s^{t,x},Z_s^{t,x})$ solves BDSDE (\ref{co}) for a.e. $x\in\mathbb{R}^d$ if and only if
$(\tilde{Y}_s^{t,x},\tilde{Z}_s^{t,x})=({\rm e}^{\mu
s}Y_s^{t,x},{\rm e}^{\mu s}Z_s^{t,x})$ solves the following BDSDE:
\begin{eqnarray}\label{e}
\tilde{Y}_s^{t,x}=\tilde{h}(X_T^{t,x})+\int_s^T\tilde{f}(r,X_r^{t,x},\tilde{Y}_r^{t,x})dr-\int_{s}^{T}\langle\tilde{g}(r,X_{r}^{t,x},\tilde{Y}_{r}^{t,x}),d^\dagger{\hat{B}}_r\rangle-\int_s^T\langle
\tilde{Z}_r^{t,x},dW_r\rangle,
\end{eqnarray}
where $\tilde{h}(x)={\rm e}^{\mu T}h(x)$, $\tilde{f}(r,x,y)={\rm
e}^{\mu r}f(r,x,{\rm e}^{-\mu r}y)-\mu y$ and $\tilde{g}(r,x,y)={\rm
e}^{\mu r}f(r,x,{\rm e}^{-\mu r}y)$. We can verify that $\tilde{h}$,
$\tilde{f}$ and $\tilde{g}$ satisfy Conditions (H.1)--(H.4) with a
possibly different constant $L$. But, by Condition (H.3), for any
$s\in[0,T]$, $y_1, y_2\in {\mathbb{R}^{1}}$, $x\in{\mathbb{R}^{d}}$,
\begin{eqnarray*}
&&(y_1-y_2)\big(\tilde{f}(s,x,y_1)-\tilde{f}(s,x,y_2)\big)\\
&=&{\rm e}^{2\mu s}({\rm e}^{-\mu s}y_1-{\rm e}^{-\mu s}y_2)\big(f(s,x,{\rm e}^{-\mu s}y_1)-f(s,x,{\rm e}^{-\mu s}y_2)\big)-\mu(y_1-y_2)(y_1-y_2)\leq0.
\end{eqnarray*}

Since $(Y_\cdot^{t,\cdot},Z_\cdot^{t,\cdot})\in
S^{2p,0}([t,T];L_{\rho}^{2p}({\mathbb{R}^{d}};{\mathbb{R}^{1}}))\times
M^{2,0}([t,T];L_{\rho}^2({\mathbb{R}^{d}};{\mathbb{R}^{d}}))$ if and
only if $(\tilde{Y}_\cdot^{t,\cdot},\tilde{Z}_\cdot^{t,\cdot})\\\in
S^{2p,0}([t,T];L_{\rho}^{2p}({\mathbb{R}^{d}};{\mathbb{R}^{1}}))\times
M^{2,0}([t,T];L_{\rho}^2({\mathbb{R}^{d}};{\mathbb{R}^{d}}))$, so we
claim $(Y_s^{t,x},Z_s^{t,x})$ is the solution of BDSDE (\ref{co}) if
and only if $(\tilde{Y}_s^{t,x},\tilde{Z}_s^{t,x})$ is the solution
of BSDE (\ref{e}). Therefore we can replace, without losing any
generality, Condition (H.3) by
\begin{description}
\item[${\rm (H.3)}^*$.] For any $s\in[0,T]$, $x\in{\mathbb{R}^{d}}$, $y_1,
y_2\in {\mathbb{R}^{1}}$,
\begin{eqnarray*}
(y_1-y_2)\big(f(s,x,y_1)-f(s,x,y_2)\big)\leq0.
\end{eqnarray*}
\end{description}

The main task of Sections 3 and 4 is to prove the following theorem about the
existence and uniqueness of the solution of BDSDE (\ref{co}).
\begin{thm}\label{21}
Under Conditions (H.1)--(H.2), (H.3)$^*$, (H.4)--(H.6), BDSDE
(\ref{co}) has a unique solution $(Y_\cdot^{t,\cdot},
Z_\cdot^{t,\cdot})\in
S^{{8p},0}([t,T];L_{\rho}^{{8p}}({\mathbb{R}^{d}};{\mathbb{R}^{1}}))\times
M^{2,0}([t,T];L_{\rho}^2({\mathbb{R}^{d}};{\mathbb{R}^{d}}))$.
\end{thm}
For this, a sequence of BDSDEs with linear growth coefficients
are constructed. Assume that $f$ in BDSDE (\ref{co}) satisfy
Conditions (H.1)-(H.2) and ${\rm (H.3)}^*$. Firstly, for each $n\in N$, define
\begin{eqnarray*}
f_n(s,x,y)=f\big(s,x,\Pi_n(y)\big)+\partial_yf(s,x,{n\over{|y|}}y)(y-{n\over{|y|}}y)I_{\{|y|>n\}},
\end{eqnarray*}
where $\Pi_n(y)={\inf(n,|y|)\over|y|}y$. 
Obviously, for
any $s\in[0,T]$, $x\in\mathbb{R}^{d}$, $y\in\mathbb{R}^{1}$,
$$f_n(s,x,y)\longrightarrow f(s,x,y),\ \ \ {\rm as}\ n\to\infty,$$ and it is easy to check that
$f_n$ satisfies the following
conditions:
\begin{description}
\item[(H.1)$'$.] For any $s\in[0,T]$, $x\in{\mathbb{R}^{d}}$, $y\in
{\mathbb{R}^{1}}$, $|f_n(s,x,y)|\leq
L(|{f}_0(s,x)|+(2{(n\wedge|y|)}^{p-1}+1)|y|)$ and $|\partial_yf_n(s,x,y)|\leq
L(1+|y|^{p-1})$. 
\item[(H.2)$'$.] 
For any $s\in[0,T]$,
$x,x_1,x_2\in{\mathbb{R}^{d}}$, $y,y_1,y_2\in {\mathbb{R}^{1}}$,
\begin{eqnarray*}
&&|f_n(s,x_1,y)-f_n(s,x_2,y)|\leq3L(1+|y|^p)|x_1-x_2|,\\
&&|\partial_yf_n(s,x_1,y)-\partial_yf_n(s,x_2,y)|\leq L(1+|y|^{p-1})|x_1-x_2|,\\
&&|\partial_yf_n(s,x,y_1)-\partial_yf_n(s,x,y_2)|\leq L(1+|y_1|^{p-2}+|y_2|^{p-2})|y_1-y_2|.
\end{eqnarray*}
\item[(H.3)$'$.] For any $s\in[0,T]$, $x\in{\mathbb{R}^{d}}$, $y_1,
y_2\in {\mathbb{R}^{1}}$,
\begin{eqnarray*}
(y_1-y_2)\big(f_n(s,x,y_1)-f_n(s,x,y_2)\big)\leq0.
\end{eqnarray*}
\end{description}
We then study the following BDSDE with
coefficient $f_n$:
\begin{eqnarray}\label{c}
Y_s^{t,x,n}&=&h(X_T^{t,x})+\int_s^Tf_n(r,X_r^{t,x},Y_r^{t,x,n})dr-\int_s^T\langle
g(r,X_r^{t,x},Y_r^{t,x,n}),d\hat{B}_r\rangle-\int_s^T\langle
Z_r^{t,x,n},dW_r\rangle.\ \ \ \
\end{eqnarray}
Notice that the coefficients $h$, $f_n$ and $g$ satisfy
Conditions (H.1)$'$--(H.3)$'$, (H.2), (H.4). Hence by Theorems 2.2 and 2.3 in
\cite{zh-zh2}, we have the following proposition:
\begin{prop}\label{7} (\cite{zh-zh2})
Under the conditions of Theorem \ref{21}, BDSDE (\ref{c}) has a
unique solution $(Y_s^{t,x,n}, Z_s^{t,x,n})\in
S^{2,0}([t,T];L_{\rho}^{2}({\mathbb{R}^{d}};{\mathbb{R}^{1}}))\times
M^{2,0}([t,T];L_{\rho}^2({\mathbb{R}^{d}};{\mathbb{R}^{d}}))$. If we
define $Y_t^{t,x,n}=u_n(t,x)$, then $u_n(t,x)$ is the unique
strong solution of the following SPDE
\begin{eqnarray}\label{f}
u_n(t,x)&=&h(x)+\int_{t}^{T}\{\mathscr{L}u_n(s,x)+f_n\big(s,x,u_n(s,x)\big)\}ds-\int_{t}^{T}\langle
g\big(s,x,u_n(s,x)\big),d^\dagger\hat{B}_s\rangle,
\end{eqnarray}
for $0\leq t\leq T$.
Moreover,
\begin{eqnarray*}
u_n(s,X_s^{t,x})=Y_s^{t,x,n},\ (\sigma^*\nabla
u_n)(s,X^{t,x}_s)=Z_s^{t,x,n}\ {\rm for}\ {\rm a.e.}\ s\in[t,T],\
x\in\mathbb{R}^{d}\ {\rm a.s.}
\end{eqnarray*}
\end{prop}

The key is to pass the limits in (\ref{c}) and (\ref{f}) in some
desired sense. For this, we need some estimates. 
\begin{lem}\label{1}
Under Conditions (H.1)--(H.2), (H.3)$^*$, (H.4)--(H.6), if
$(Y_\cdot^{t,\cdot,n},Z_\cdot^{t,\cdot,n})$ is the solution of BDSDE
(\ref{c}), then we have for any $2\leq m\leq {8p}$,
\begin{eqnarray*}
&&\sup_nE[\sup_{s\in[t,T]}\int_{\mathbb{R}^d}|Y_{s}^{t,x,n}|^{m}\rho^{-1}(x)dx]+\sup_nE[\int_t^T\int_{\mathbb{R}^d}|Y_s^{t,x,n}|^{m}\rho^{-1}(x)dxds]\\
&&+\sup_nE[\int_{t}^{T}\int_{\mathbb{R}^{d}}{{|{Y}_s^{t,x,n}|}^{m-2}}|{Z}_{s}^{t,x,n}|^2\rho^{-1}(x)dxds]
+\sup_nE[\big(\int_t^T\int_{\mathbb{R}^d}|Z_s^{t,x,n}|^{2}\rho^{-1}(x)dxds\big)^{{m\over2}}]<\infty.
\end{eqnarray*}
\end{lem}

The proof of the lemma follows some standard It$\hat {\rm o}$'s formula computation. So it is omitted here.

Taking $m=2$ in Lemma \ref{1}, we know
\begin{eqnarray*}\label{an}
\sup_nE[\int_t^T\int_{\mathbb{R}^d}|Y_s^{t,x,n}|^2\rho^{-1}(x)dxds]+\sup_nE[\int_t^T\int_{\mathbb{R}^d}|Z_s^{t,x,n}|^2\rho^{-1}(x)dxds]<\infty.
\end{eqnarray*}
Define $U_s^{t,x,n}=f_n(s,X_s^{t,x},Y_s^{t,x,n})$ and
$V_s^{t,x,n}=g(s,X_s^{t,x},Y_s^{t,x,n})$, $s\geq t$. Using Lemma
\ref{1} again, we also have
\begin{eqnarray}\label{p}
&&\sup_nE[\int_t^T\int_{\mathbb{R}^d}(|Y_s^{t,x,n}|^2+|Z_s^{t,x,n}|^2+|U_s^{t,x,n}|^2+|V_s^{t,x,n}|^2)\rho^{-1}(x)dxds]\nonumber\\
&\leq&\sup_nC_pE[\int_t^T\int_{\mathbb{R}^d}(1+|f_0(s,X_s^{t,x})|^2+|g(s,X_s^{t,x},0)|^2+|Y_s^{t,x,n}|^{2p}+|Z_s^{t,x,n}|^2)\rho^{-1}(x)dxds]\nonumber\\
&<&\infty.
\end{eqnarray}
Here and in the following $C_p$ is a generic constant.
Then, according to the Alaoglu lemma, we know that there exists a
subsequence, still denoted by $(Y_\cdot^{t,\cdot,n},
Z_\cdot^{t,\cdot,n}, U_\cdot^{t,\cdot,n}, V_\cdot^{t,\cdot,n})$,
converging weakly to a
limit $(Y_\cdot^{t,\cdot}, Z_\cdot^{t,\cdot}, U_\cdot^{t,\cdot},
V_\cdot^{t,\cdot})$ in
$L^2_\rho(\Omega\times[t,T]\times\mathbb{R}^d;\mathbb{R}^1\times\mathbb{R}^d\times\mathbb{R}^1\times\mathbb{R}^l)$, or equivalently
$L^2(\Omega\times[t,T];L_\rho^2(\mathbb{R}^d;\mathbb{R}^1)\times
L_\rho^2(\mathbb{R}^d;\mathbb{R}^d)\times
L_\rho^2(\mathbb{R}^d;\mathbb{R}^1)\times
L_\rho^2(\mathbb{R}^d;\mathbb{R}^l))$. Now we take the weak limit in
$L^2_\rho(\Omega\times[t,T]\times\mathbb{R}^d;\mathbb{R}^1)$) to
BDSDEs (\ref{c}), we can verify that $(Y_s^{t,x}, Z_s^{t,x},
U_s^{t,x}, V_s^{t,x})$ satisfies the following BDSDE:
\begin{eqnarray}\label{d}
Y_s^{t,x}=h(X_T^{t,x})+\int_s^TU_r^{t,x}dr-\int_s^T\langle
V_r^{t,x},d\hat{B}_r\rangle-\int_s^T\langle Z_r^{t,x},dW_r\rangle.
\end{eqnarray}
For this, we will check the weak convergence term by term. The weak
convergence of $Y_s^{t,x,n}$ is deduced by the definition of
$Y_s^{t,x}$. 
We check the weak convergence of $\int_s^TU_r^{t,x,n}dr$. Let
$\eta\in
L^2_\rho(\Omega\times[t,T]\times\mathbb{R}^d;\mathbb{R}^1)$.
Noticing
$\int_t^T\sup_nE[\int_s^T\int_{\mathbb{R}^d}|U_r^{t,x,n}|^2\rho^{-1}(x)dxdr]ds<\infty$
due to (\ref{p}), by Lebesgue's dominated convergence theorem, we
have
\begin{eqnarray*}
&&|E[\int_t^T\int_{\mathbb{R}^d}\int_s^T(U_r^{t,x,n}-U_r^{t,x})dr\eta(s,x)\rho^{-1}(x)dxds]|\nonumber\\
&\leq&\int_t^T|E[\int_s^T\int_{\mathbb{R}^d}(U_r^{t,x,n}-U_r^{t,x})\eta(s,x)\rho^{-1}(x)dxdr]|ds\longrightarrow0,\
\ \ {\rm as}\ n\to\infty.
\end{eqnarray*}
To prove the weak convergence of $\int_s^T\langle
V_r^{t,x,n},d\hat{B}_r\rangle$, first note that for fixed $s$ and
$x$, $\eta(s,x)\in L^2(\Omega;\mathbb{R}^1)$, hence there exist
$\psi\in
L^2(\Omega\times[t,T]\times\mathbb{R}^d\times[t,T];\mathbb{R}^l)$
and $\phi\in
L^2(\Omega\times[t,T]\times\mathbb{R}^d\times[t,T];\mathbb{R}^d)$
s.t.
$\eta(s,x)=E[\eta(s,x)]+\int_t^T\langle\psi(s,x,r),d\hat{B}_r\rangle+\int_t^T\langle\phi(s,x,r),dW_r\rangle$.
Noticing that for a.e. $s\in[t,T]$, $\psi(s,\cdot,\cdot)\in
L^2(\Omega\times\mathbb{R}^d\times[t,T];\mathbb{R}^l)$,
$\phi(s,\cdot,\cdot)\in
L^2(\Omega\times\mathbb{R}^d\times[t,T];\mathbb{R}^d)$ and
$\int_t^T\sup_nE[\int_s^T\int_{\mathbb{R}^d}|V_r^{t,x,n}|^2\rho^{-1}(x)dxdr]ds<\infty$,
by Lebesgue's dominated convergence theorem again, we obtain
\begin{eqnarray*}
&&|E[\int_t^T\int_{\mathbb{R}^d}\int_s^T\langle V_r^{t,x,n}-V_r^{t,x},d\hat{B}_r\rangle\eta(s,x)\rho^{-1}(x)dxds]|\nonumber\\
&=&|\int_t^T\int_{\mathbb{R}^d}E[\int_s^T\langle V_r^{t,x,n}-V_r^{t,x},\psi(s,x,r)\rangle dr]\rho^{-1}(x)dxds|\nonumber\\
&\leq&\int_t^T|E[\int_s^T\int_{\mathbb{R}^d}\langle
V_r^{t,x,n}-V_r^{t,x},\psi(s,x,r)\rangle
\rho^{-1}(x)dxdr]|ds\longrightarrow0,\ \ \ {\rm as}\ n\to\infty.
\end{eqnarray*}
For the weak convergence of last term, we can deduce similarly that
\begin{eqnarray*}
&&|E[\int_t^T\int_{\mathbb{R}^d}\int_s^T\langle Z_r^{t,x,n}-Z_r^{t,x},dW_r\rangle\eta(s,x)\rho^{-1}(x)dxds]|\nonumber\\
&\leq&\int_t^T|E[\int_s^T\int_{\mathbb{R}^d}\langle
Z_r^{t,x,n}-Z_r^{t,x},\phi(s,x,r)\rangle
\rho^{-1}(x)dxdr]|ds\longrightarrow0,\ \ \ {\rm as}\ n\to\infty.
\end{eqnarray*}

Needless to say, if we can show that BDSDE (\ref{c}) converges
weakly to BDSDE (\ref{co}) as $n\to\infty$, then we can say
$(Y_s^{t,x},Z_s^{t,x})$ is a solution of BDSDE (\ref{co}). The key
is to prove that $U^{t,x}_s=f(s,X_s^{t,x},Y_s^{t,x})$ and
$V^{t,x}_s=g(s,X_s^{t,x},Y_s^{t,x})$ for a.e. $s\in[t,T]$,
$x\in\mathbb{R}^d$ a.s. However, the weak convergence of $Y^n$,
$U^n$, $V^n$ and $Z^n$ is far from enough for this purpose. The real difficulty in this analysis is to establish the strong convergence of
$Y^n$ and
$Z^n$, at least along a subsequence. 

\section{The strong convergence}
\setcounter{equation}{0}

To obtain the strongly convergent subsequence of $Y^n$ and $Z^n$, we need to
estimate the Malliavin derivatives to prove the relative compactness
of $Y_n$ first. Let $\mathcal{O}$ be a bounded domain in
$\mathbb{R}^d$. Denote $C^k_c(\mathcal{O})$ the class of $k$-times
differentiable functions which have a compact support inside
$\mathcal{O}$. For $\varphi\in C^k_c(\mathcal{O})$, we define
$v^\varphi(s,\omega)=\int_\mathcal{O}v(s,x,\omega)\varphi(x)dx$. The
following theorem proved in Bally and Saussereau \cite{ba-sa} can be regarded as an
extension of Rellich-Kondrachov compactness theorem to stochastic
case. This kind of Wiener-Sobolev compactness theorem for time and
space independent case was considered by Da Prato and Malliavin
\cite{da-ma}, Peszat \cite{pes}. One extension was given in Feng,
Zhao and Zhou \cite{fe-zh-zh}, Feng and Zhao \cite{fz1} to replace $L^2$ norm in two variable by $\sup$ norm in the two variables, in order to apply it to infinite
horizon stochastic integral equations. For the purpose of this paper, the $L^2$ norm used by Bally and Saussereau is enough.

Denote by $C^\infty_p(\mathbb{R}^n)$ the set of infinitely differentiable functions $f: \mathbb{R}^n\longrightarrow\mathbb{R}^1$ such
that $f$ and all its partial derivatives have polynomial growth. Let $\mathbb{K}$ be the class of smooth random
variables $F$ that is $F = f (W(h_1),\cdot\cdot\cdot,W(h_n))$ with $n\in\mathbb{N}$, $h_1,\cdot\cdot\cdot,h_n\in L^2([0,T])$ and $f\in C^\infty_p(\mathbb{R}^n)$. The derivative operator of a smooth random variable $F$ is the stochastic process $\{D_t F,\ t\in[0,T]\}$ defined
by (cf. \cite{nu})
$$D_t F =\sum_{i=1}^n{{\partial f}\over{\partial x^i}}\big(W(h_1),\cdot\cdot\cdot,W(h_n)\big)h_i(t).$$
We will denote $\mathbb{D}^{1,2}$ the domain of $D$ in $L^2(\Omega)$, i.e. $\mathbb{D}^{1,2}$ is the closure of $\mathbb{K}$ with respect to the norm
$$\|F\|_{1,2}^2= E[|F|^2] + E[\|D_t F\|^2_{L^2([0,T])}].$$
Recall
\begin{thm}\label{3}(Theorem 2, \cite{ba-sa})
Let $(u_n)_{n\in\mathbb{N}}$ be a sequence of
$L^2([0,T]\times\Omega;H^1(\mathcal{O}))$. Suppose that\\
(1)
$\sup_{n}E[\int_0^T\|u_n(s,\cdot)\|^2_{H^1(\mathcal{O})}ds]<\infty$.\\
(2) For all $\varphi\in C^k_c(\mathcal{O})$ and $t\in[0,T]$,
$u_n^\varphi(s)\in\mathbb{D}^{1,2}$ and
$\sup_{n}\int_0^T\|u_n^\varphi(s)\|^2_{\mathbb{D}^{1,2}}ds<\infty$.\\
(3) For all $\varphi\in C^k_c(\mathcal{O})$, the sequence
$(E[u_n^\varphi])_{n\in\mathbb{N}}$ of $L^2([0,T])$ satisfies\\
(3i) For any $\varepsilon>0$, there exists $0<\alpha<\beta<T$ s.t.
$$\sup_{n}\int_{[0,T]\setminus(\alpha,\beta)}|E[u_n^\varphi(s)]|^2ds<\varepsilon.$$
(3ii) For any $0<\alpha<\beta<T$ and $h\in\mathbb{R}^1$ s.t.
$|h|<min(\alpha,T-\beta)$, it holds
$$\sup_{n}\int_\alpha^\beta|E[u_n^\varphi(s+h)]-E[u_n^\varphi(s)]|^2ds<C_p|h|.$$
(4) For all $\varphi\in C^k_c(\mathcal{O})$, the following
conditions are satisfied:\\
(4i) For any $\varepsilon>0$, there exists $0<\alpha<\beta<T$ and
$0<\alpha'<\beta'<T$ s.t.
$$\sup_{n}E[\int_{[0,T]^2\setminus(\alpha,\beta)\times(\alpha',\beta')}|D_\theta u_n^\varphi(s)|^2d\theta ds]<\varepsilon.$$
(4ii) For any $0<\alpha<\beta<T$, $0<\alpha'<\beta'<T$ and
$h,h'\in\mathbb{R}^1$ s.t.
$max(|h|,|h'|)<min(\alpha,\alpha',T-\beta,T-\beta')$, it holds that
$$\sup_{n}E[\int_\alpha^\beta\int_{\alpha'}^{\beta'}|D_{\theta+h}u_n^\varphi(s+h')-D_\theta
u_n^\varphi(s)|^2d\theta ds]<C_p(|h|+|h'|).$$ Then
$(u_n)_{n\in\mathbb{N}}$ is relatively compact in
$L^2(\Omega\times[0,T]\times\mathcal{O};\mathbb{R}^1)$.
\end{thm}

Using Theorem \ref{3}, we can verify that the sequence $u_n(s,x)$ in
SPDE (\ref{f}) is relatively compact. In this process, some
estimates on the Malliavin derivative of the random variable $(Y,Z)$
w.r.t. Brownian motion $B$ are needed. In what follows, we will need the following results
whose proofs are deferred to Section 7. Throughout this paper, Malliavin derivative always refers to Malliavin derivative w.r.t. $B$ unless we say otherwise.
\begin{lem}\label{26}
Under Conditions (H.1)--(H.2), (H.3)$^*$, (H.4)--(H.6), the
Malliavin derivative of the solution $(Y_s^{t,x,n},Z_s^{t,x,n})$ of
BDSDE (\ref{c}) exists and satisfies the following linear equation:
\begin{numcases}{}\label{aq}
D_\theta
Y_s^{t,x,n}=g(\theta,X_\theta^{t,x},Y_\theta^{t,x,n})+\int_s^\theta\partial_yf_n(r,X_r^{t,x},Y_r^{t,x,n})D_\theta
Y_r^{t,x,n}dr
\nonumber\\
\ \ \ \ \ \ \ \ \ \ \ \ \ \
-\int_s^\theta\partial_yg(r,X_r^{t,x},Y_r^{t,x,n})D_\theta
Y_r^{t,x,n}d^\dagger\hat{B}_r-\int_s^\theta D_\theta
Z_r^{t,x,n}dW_r,\ \ \ 0\leq\theta\leq T,\nonumber\\
D_\theta Y_s^{t,x,n}=0,\ \ \ t\leq\theta<s.
\end{numcases}
Moreover, for any $2\leq m\leq {8p}$,
\begin{eqnarray*}\label{at}
&&\sup_n\sup_{\theta\in[t,T]}E[\sup_{s\in[t,T]}\int_{\mathbb{R}^d}|D_\theta Y_{s}^{t,x,n}|^{m}\rho^{-1}(x)dx]+\sup_n\sup_{\theta\in[t,T]}E[\int_t^T\int_{\mathbb{R}^d}|D_\theta Y_s^{t,x,n}|^{m}\rho^{-1}(x)dxds]\nonumber\\
&&+\sup_n\sup_{\theta\in[t,T]}E[\int_{t}^{T}\int_{\mathbb{R}^{d}}{{|D_\theta{Y}_s^{t,x,n}|}^{m-2}}|D_\theta{Z}_{s}^{t,x,n}|^2\rho^{-1}(x)dxds]\nonumber\\
&&+\sup_n\sup_{\theta\in[t,T]}E[\big(\int_t^T\int_{\mathbb{R}^d}|D_\theta
Z_s^{t,x,n}|^{2}\rho^{-1}(x)dxds\big)^{{m\over2}}]<\infty.
\end{eqnarray*}
\end{lem}
{\em Proof}. From Condition (H.2)$'$ and the
results of \cite{ba-sa} or \cite{pa-pe3}, it is easy to know that
the Malliavin derivative of $(Y_s^{t,x,n},Z_s^{t,x,n})$ exists and
satisfies (\ref{aq}). The rest of the proof follows some standard computations using It$\hat {\rm o}$'s formula. So it is omitted here.
$\hfill \diamond$\\

Now we are ready to prove the relative compactness of the solutions
of SPDE (\ref{f}) in the following theorem. In the proof
of Theorem \ref{2} and the similar arguments throughout this paper, we will
leave out the similar localization procedures as in the proof of
Lemma \ref{1} when applying generalized It$\hat {\rm
o}$'s formula, due to the limitation of the length of the paper. 
\begin{thm}\label{2}
Under Conditions (H.1)--(H.2), (H.3)$^*$, (H.4)--(H.6), if
$(Y_\cdot^{t,\cdot,n},Z_\cdot^{t,\cdot,n})$ is the solution of BDSDE
(\ref{c}) and $\mathcal{O}$ is a bounded domain in $\mathbb{R}^d$,
then the sequence $u_n(s,x)\triangleq Y_s^{s,x,n}$ is relatively
compact in $L^2(\Omega\times[0,T]\times\mathcal{O};\mathbb{R}^1)$.
\end{thm}
{\em Proof}. We verify that $u_n$ satisfies Conditions (1)--(4) in
Theorem \ref{3}.

Step 1. We first verify Condition (1). By Conditions (H.5)--(H.6), Lemma \ref{qi045},
Proposition \ref{7} and Lemma \ref{1}, we have
\begin{eqnarray}\label{dc}
\sup_{n}E[\int_0^T\|u_n(s,\cdot)\|^2_{H^1(\mathcal{O})}ds]
&\leq&C_p\sup_{n}E[\int_0^T\int_\mathcal{O}(|u_n(s,x)|^2+|\nabla
u_n(s,x)|^2)\rho^{-1}(x)dxds]\nonumber\\
&\leq&
C_p\sup_{n}E[\int_0^T\int_{\mathbb{R}^d}(|Y_s^{0,x,n}|^2+|Z_s^{0,x,n}|^2)\rho^{-1}(x)dxds]<\infty.
\end{eqnarray}

Step 2. We then verify Condition (2). It is easy to see that
$D_\theta u_n^\varphi(s)=\int_\mathcal{O}D_\theta
u_n(s,x)\varphi(x)dx$. By Lemma \ref{26}, $D_\theta
u_n(s,x)=D_\theta Y_s^{s,x,n}$ exists. Indeed, we can further prove
$u_n^\varphi(s)\in\mathbb{D}^{1,2}$. 
By stochastic calculus, we have 
\begin{eqnarray}\label{db}
\|u_n^\varphi(s)\|^2_{\mathbb{D}^{1,2}}
&\leq&C_pE[\int_{\mathbb{R}^d}|u_n(s,x)|^2\rho^{-1}(x)dx]+C_pE[\int_s^T\int_{\mathbb{R}^d}|D_\theta u_n(s,x)|^2\rho^{-1}(x)dxd\theta]\nonumber\\
&\leq&C_pE[\int_{\mathbb{R}^d}|Y_s^{0,x,n}|^2\rho^{-1}(x)dx]+C_pE[\int_s^T\int_{\mathbb{R}^d}|D_\theta Y_s^{s,x,n}|^2\rho^{-1}(x)dxd\theta]\nonumber\\
&\leq&C_p\sup_n\int_{\mathbb{R}^d}|h(x)|^2\rho^{-1}(x)dx+C_p\int_{0}^{T}\int_{\mathbb{R}^d}|f_0(r,x)|^2\rho^{-1}(x)dxdr\nonumber\\
&&+C_p\int_{0}^{T}\int_{\mathbb{R}^d}|g(r,x,0)|^2\rho^{-1}(x)dxdr<\infty.
\end{eqnarray}
Also the right hand side of the above inequality is independent of
$s$ and $n$, so
$$\sup_{n}\int_0^T\|u_n^\varphi(s)\|^2_{\mathbb{D}^{1,2}}ds<\infty.$$

Step 3. Let us verify Condition (3). 
First (3i) follows immediately from (\ref{db}). To see (3ii), assume $h>0$
without losing any generality. From
(\ref{c}) and Cauchy-Schwarz inequality,  we have
\begin{eqnarray}\label{df}
\sup_{n}\int_\alpha^\beta|E[u_n^\varphi(s+h)]-E[u_n^\varphi(s)]|^2ds
&\leq&C_p\sup_{n}\int_\alpha^\beta E[\int_{\mathbb{R}^d}|u_n(s+h,x)-u_n(s,x)|^2\rho^{-1}(x)dx]ds\nonumber\\
&\leq&C_p\sup_{n}\int_\alpha^\beta E[\int_{\mathbb{R}^d}|Y_{s+h}^{0,x,n}-Y_s^{0,x,n}|^2\rho^{-1}(x)dx]ds\nonumber\\
&\leq&C_p\int_\alpha^\beta\int_s^{s+h}\int_{\mathbb{R}^d}|f_0(r,x)|^2\rho^{-1}(x)dxdrds\nonumber\\
&&+C_p\sup_{n}\int_\alpha^\beta\int_s^{s+h}\sup_{r\in[0,T]}(1+E[\int_{\mathbb{R}^d}|Y_r^{0,x,n}|^{2p}\rho^{-1}(x)dx])drds\nonumber\\
&&+C_p\sup_{n}\int_\alpha^\beta\int_s^{s+h}E[\int_{\mathbb{R}^d}|Z_r^{0,x,n}|^{2}\rho^{-1}(x)dx]drds\nonumber\\
&&+C_p\int_\alpha^\beta\int_s^{s+h}\sup_{r\in[0,T]}\int_{\mathbb{R}^d}|g(r,x,0)|^{2}\rho^{-1}(x)dxdrds.
\end{eqnarray}
Note that by changing integration order,
\begin{eqnarray}\label{di}
&&\sup_{n}\int_\alpha^\beta\int_s^{s+h}E[\int_{\mathbb{R}^d}|Z_r^{0,x,n}|^{2}\rho^{-1}(x)dx]drds\nonumber\\
&=&\sup_{n}\big(\int_{\alpha}^{\alpha+h}\int_{\alpha}^{r}E[\int_{\mathbb{R}^d}|Z_r^{0,x,n}|^2\rho^{-1}(x)dx]dsdr+\int_{\alpha+h}^{\beta}\int_{r-h}^{r}E[\int_{\mathbb{R}^d}|Z_r^{0,x,n}|^2\rho^{-1}(x)dx]dsdr\nonumber\\
&&\ \ \ \ \ \ \ \ +\int_{\beta}^{\beta+h}\int_{r-h}^{\beta}E[\int_{\mathbb{R}^d}|Z_r^{0,x,n}|^2\rho^{-1}(x)dx]dsdr\big)\nonumber\\
&=&\sup_{n}\big((r-\alpha)\int_{\alpha}^{\alpha+h}E[\int_{\mathbb{R}^d}|Z_r^{0,x,n}|^2\rho^{-1}(x)dx]dr+h\int_{\alpha+h}^{\beta}E[\int_{\mathbb{R}^d}|Z_r^{0,x,n}|^2\rho^{-1}(x)dx]dr\nonumber\\
&&\ \ \ \ \ \ \ \ +(\beta+h-r)\int_{\beta}^{\beta+h}E[\int_{\mathbb{R}^d}|Z_r^{0,x,n}|^2\rho^{-1}(x)dx]dr\big)\nonumber\\
&=&C_ph\sup_{n}E[\int_{0}^{T}\int_{\mathbb{R}^d}|Z_r^{0,x,n}|^2\rho^{-1}(x)dxdr].
\end{eqnarray}
A similar calculation can be done to $\int_\alpha^\beta\int_s^{s+h}\int_{\mathbb{R}^d}|f_0(r,x)|^2\rho^{-1}(x)dxdrds$. Hence it follows from (\ref{df}) that
\begin{eqnarray*}
&&\sup_{n}\int_\alpha^\beta|E[u_n^\varphi(s+h)]-E[u_n^\varphi(s)]|^2ds\nonumber\\
&\leq&C_ph\int_0^T\int_{\mathbb{R}^d}|f_0(r,x)|^2\rho^{-1}(x)dxdr+C_ph\sup_{n}\sup_{r\in[0,T]}(1+E[\int_{\mathbb{R}^d}|Y_r^{0,x,n}|^{2p}\rho^{-1}(x)dx])\nonumber\\
&&+C_ph\sup_{n}E[\int_{0}^{T}\int_{\mathbb{R}^d}|Z_r^{0,x,n}|^2\rho^{-1}(x)dxdr]+C_ph\sup_{r\in[0,T]}\int_{\mathbb{R}^d}|g(r,x,0)|^{2}\rho^{-1}(x)dx.
\end{eqnarray*}
Also noticing Condition (H.1), (\ref{au}) and Lemma \ref{1}, we can conclude that (3ii) holds.

Step 4. We now verify Condition (4). For (4i), since by the equivalence of norm principle 
it turns out that
\begin{eqnarray*}
\sup_n\sup_{\theta\in[t,T]}\sup_{s\in[t,T]}E[|D_\theta u_n^\varphi(s)|^2]
&\leq&C_p\sup_n\sup_{\theta\in[t,T]}\sup_{s\in[t,T]}E[\int_{\mathbb{R}^d}|D_\theta u_n(s,x)|^2\rho^{-1}(x)dx]\\
&\leq&C_p\sup_n\sup_{\theta\in[t,T]}\sup_{s\in[t,T]}E[\int_{\mathbb{R}^d}|D_\theta
Y_s^{0,x,n}|^2\rho^{-1}(x)dx]<\infty.
\end{eqnarray*}
So (4i) follows. 
To see (4ii), assume without losing any generality
that $h,h'>0$, then 
\begin{eqnarray}\label{ap}
&&\sup_{n}E[\int_\alpha^\beta\int_{\alpha'}^{\beta'}|D_{\theta+h}u_n^\varphi(s+h')-D_\theta u_n^\varphi(s)|^2d\theta ds]\nonumber\\
&\leq&C_p\sup_{n}E[\int_{\alpha+h'}^{\beta+h'}\int_{\alpha'}^{\beta'}\int_\mathcal{O}|D_{\theta+h}u_n(s,x)-D_{\theta}u_n(s,x)|^2\rho^{-1}(x)dxd\theta ds]\nonumber\\
&&+C_p\sup_{n}E[\int_\alpha^\beta\int_{\alpha'}^{\beta'}\int_\mathcal{O}|D_{\theta}u_n(s+h',x)-D_{\theta}u_n(s,x)|^2\rho^{-1}(x)dxd\theta
ds].
\end{eqnarray}
For the first term on the right hand side of (\ref{ap}), by the equivalence of norm principle,
\begin{eqnarray}\label{bc}
&&\sup_{n}E[\int_{\alpha+h'}^{\beta+h'}\int_{\alpha'}^{\beta'}\int_\mathcal{O}|D_{\theta+h}u_n(s,x)-D_{\theta}u_n(s,x)|^2\rho^{-1}(x)dxd\theta ds]\nonumber\\
&=&\sup_{n}E[\int_{\alpha+h'}^{\beta+h'}\int_{\alpha'}^{\beta'}\int_{\mathbb{R}^d}|D_{\theta+h}Y_s^{s,x,n}-D_{\theta}Y_s^{s,x,n}|^2\rho^{-1}(x)dxd\theta
ds]\nonumber\\
&\leq&C_p\sup_{n}E[\int_{\alpha+h'}^{\beta+h'}\int_{\alpha'}^{\beta'}\int_{\mathbb{R}^d}|D_{\theta+h}Y_s^{0,x,n}-D_{\theta}Y_s^{0,x,n}|^2\rho^{-1}(x)dxd\theta
ds].
\end{eqnarray}
By BDSDE (\ref{aq}) we know that
\begin{eqnarray*}
(D_{\theta+h}-D_\theta)Y_s^{0,x,n}
&=&H(\theta,\theta+h)+\int_s^{\theta}\partial_yf_n(r,X_r^{0,x},Y_r^{0,x,n})(D_{\theta+h}-D_\theta)
Y_r^{0,x,n}
dr\nonumber\\
&&-\int_s^{\theta}\partial_yg(r,X_r^{0,x},Y_r^{0,x,n})(D_{\theta+h}-D_\theta)
Y_r^{0,x,n}d^\dagger\hat{B}_r-\int_s^{\theta}(D_{\theta+h}-D_\theta)
Z_r^{0,x,n}dW_r,
\end{eqnarray*}
where
\begin{eqnarray*}
H(\theta,\theta+h)&=&g(\theta+h,X_{\theta+h}^{0,x},Y_{\theta+h}^{0,x,n})-g(\theta,X_{\theta}^{0,x},Y_{\theta}^{0,x,n})+\int_\theta^{\theta+h}\partial_yf_n(r,X_r^{0,x},Y_r^{0,x,n})D_{\theta+h}
Y_r^{0,x,n}dr\nonumber\\
&&
-\int_\theta^{\theta+h}\partial_yg(r,X_r^{0,x},Y_r^{0,x,n})D_{\theta+h}
Y_r^{0,x,n}d^\dagger\hat{B}_r-\int_\theta^{\theta+h}D_{\theta+h}Z_r^{0,x,n}dW_r.
\end{eqnarray*}
Applying It$\hat {\rm o}$'s formula to ${\rm
e}^{Kr}|(D_{\theta+h}-D_\theta)Y_s^{0,x,n}|^2$, 
we have
\begin{eqnarray}\label{bb}
&&E[\int_{\mathbb{R}^d}|(D_{\theta+h}-D_\theta)Y_s^{0,x,n}|^2\rho^{-1}(x)dx]+E[\int_s^\theta\int_{\mathbb{R}^d}|(D_{\theta+h}-D_\theta)
Y_r^{0,x,n}|^2\rho^{-1}(x)dxdr]\nonumber\\
&&+E[\int_s^\theta\int_{\mathbb{R}^d}|(D_{\theta+h}-D_\theta)
Z_r^{0,x,n}|^2\rho^{-1}(x)dxdr]\nonumber\\
&\leq&C_pE[\int_{\mathbb{R}^d}|H(\theta,\theta+h)|^2\rho^{-1}(x)dx].
\end{eqnarray}
Next we prove that
\begin{eqnarray}\label{dg}
\sup_{n}E[\int_{\alpha+h'}^{\beta+h'}\int_{\alpha'}^{\beta'}\int_{\mathbb{R}^d}|H(\theta,\theta+h)|^2\rho^{-1}(x)dxd\theta
ds]\leq C_ph.
\end{eqnarray}
First note that
\begin{eqnarray}\label{ba}
&&E[\int_{\mathbb{R}^d}|H(\theta,\theta+h)|^2\rho^{-1}(x)dx]\\
&\leq&C_ph^2+C_pE[\int_{\mathbb{R}^d}|X_{\theta+h}^{0,x}-X_{\theta}^{0,x}|^2\rho^{-1}(x)dx]+C_pE[\int_{\mathbb{R}^d}|Y_{\theta+h}^{0,x,n}-Y_{\theta}^{0,x,n}|^2\rho^{-1}(x)dx]\nonumber\\
&&+C_pE[\int_\theta^{\theta+h}\int_{\mathbb{R}^d}|\partial_yf_n(r,X_r^{0,x},Y_r^{0,x,n})|^2|D_{\theta+h}
Y_r^{0,x,n}|^2\rho^{-1}(x)dxdr]\nonumber\\
&&+C_p\int_\theta^{\theta+h}\sup_{s\in[0,T]}\sup_{r\in[0,T]}E[\int_{\mathbb{R}^d}|D_sY_r^{0,x,n}|^2\rho^{-1}(x)dx]dr+C_pE[\int_\theta^{\theta+h}\int_{\mathbb{R}^d}|D_sZ_r^{0,x,n}|^2\rho^{-1}(x)dxdr].\nonumber
\end{eqnarray}
We need to estimate each term in the above formula. From (\ref{a}),
we have
\begin{eqnarray*}\label{aw}
&&E[\int_{\mathbb{R}^d}|X_{\theta+h}^{0,x}-X_{\theta}^{0,x}|^2\rho^{-1}(x)dx]\nonumber\\
&\leq&C_pE[\int_{\mathbb{R}^d}(\int_\theta^{\theta+h}|b(X_{u}^{0,x})|du)^2\rho^{-1}(x)dx]+C_p\int_{\mathbb{R}^d}E[\int_\theta^{\theta+h}|\sigma(X_{u}^{0,x})|^2du]\rho^{-1}(x)dx\nonumber\\
&\leq&C_pE[\int_{\mathbb{R}^d}h\int_\theta^{\theta+h}(1+|X_{u}^{0,x}|)^{2}du\rho^{-1}(x)dx]+C_pE[\int_{\mathbb{R}^d}\int_\theta^{\theta+h}L^2du\rho^{-1}(x)dx]\nonumber\\
&\leq&C_ph\int_\theta^{\theta+h}E[\int_{\mathbb{R}^d}(1+|X_{u}^{0,x}|)^{2}\rho^{-1}(x)dx]du+C_phE[\int_{\mathbb{R}^d}L^2\rho^{-1}(x)dx]\nonumber\\
&\leq&C_ph.
\end{eqnarray*}
By (H.1)$'$, Lemma \ref{1} and Lemma \ref{26}, we have
\begin{eqnarray*}\label{ay}
&&E[\int_\theta^{\theta+h}\int_{\mathbb{R}^d}|\partial_yf_n(r,X_r^{0,x},Y_r^{0,x,n})|^2|D_{\theta+h}
Y_r^{0,x,n}|^2\rho^{-1}(x)dxdr]\nonumber\\
&\leq&C_p\int_\theta^{\theta+h}\sup_{n}\sup_{\theta\in[0,T]}\sup_{r\in[0,T]}\big(\sqrt{E[\int_{\mathbb{R}^d}(1+|Y_r^{0,x,n}|^{4p-4})\rho^{-1}(x)dx]}\sqrt{E[\int_{\mathbb{R}^d}|D_{\theta}
Y_r^{0,x,n}|^4\rho^{-1}(x)dx]}\big)dr\nonumber\\
&\leq &
C_ph.
\end{eqnarray*}
By Lemma \ref{26} again, we also have that
\begin{eqnarray*}\label{az}
\int_\theta^{\theta+h}\sup_{n}\sup_{s\in[0,T]}\sup_{r\in[0,T]}E[\int_{\mathbb{R}^d}|D_sY_r^{0,x,n}|^2\rho^{-1}(x)dx]dr\leq
C_ph.
\end{eqnarray*}
Hence, from (\ref{ba}), to prove (\ref{dg}) is reduced to prove
\begin{eqnarray}\label{dh}
&&\sup_{n}\int_{\alpha+h'}^{\beta+h'}\int_{\alpha'}^{\beta'}E[\int_\theta^{\theta+h}\int_{\mathbb{R}^d}|Y_{\theta+h}^{0,x,n}-Y_{\theta}^{0,x,n}|^2\rho^{-1}(x)dxdr]d\theta
ds\nonumber\\
&&+\sup_{n}\int_{\alpha+h'}^{\beta+h'}\int_{\alpha'}^{\beta'}E[\int_\theta^{\theta+h}\int_{\mathbb{R}^d}|D_sZ_r^{0,x,n}|^2\rho^{-1}(x)dxdr]d\theta
ds\leq C_ph.
\end{eqnarray}
From (\ref{c}), we have
\begin{eqnarray}\label{dj}
&&E[\int_{\mathbb{R}^d}|Y_{\theta+h}^{0,x,n}-Y_{\theta}^{0,x,n}|^2\rho^{-1}(x)dx]\nonumber\\
&\leq&C_pE[\int_\theta^{\theta+h}\int_{\mathbb{R}^d}|f_n(r,X_r^{0,x},Y_r^{0,x,n})|^2\rho^{-1}(x)dxdr]\\
&&+C_pE[\int_\theta^{\theta+h}\int_{\mathbb{R}^d}|g(r,X_r^{0,x},Y_r^{0,x,n})|^2\rho^{-1}(x)dxdr]+C_pE[\int_\theta^{\theta+h}\int_{\mathbb{R}^d}|Z_r^{0,x,n}|^2\rho^{-1}(x)dxdr]\nonumber\\
&\leq&C_p\int_\theta^{\theta+h}\int_{\mathbb{R}^d}|f_0(r,x)|^2\rho^{-1}(x)dxdr+C_p\int_\theta^{\theta+h}\sup_n\sup_{r\in[0,T]}E[\int_{\mathbb{R}^d}|Y_r^{0,x,n}|^{2p}\rho^{-1}(x)dx]dr\nonumber\\
&&+C_p\int_\theta^{\theta+h}E[\int_{\mathbb{R}^d}|Z_r^{0,x,n}|^2\rho^{-1}(x)dx]dr+C_p\int_\theta^{\theta+h}\sup_{r\in[0,T]}E[\int_{\mathbb{R}^d}|g(r,x,0)|^2\rho^{-1}(x)dx]dr.\nonumber
\end{eqnarray}
A similar calculation of changing the integrations order as in (\ref{di}) leads to
\begin{eqnarray*}
&&\sup_{n}\int_{\alpha+h'}^{\beta+h'}\int_{\alpha'}^{\beta'}\int_\theta^{\theta+h}E[\int_{\mathbb{R}^d}(|f_0(r,x)|^2+|Z_r^{0,x,n}|^2)\rho^{-1}(x)dx]drd\theta
ds\nonumber\\
&\leq&C_ph\sup_{n}E[\int_{0}^{T}\int_{\mathbb{R}^d}(|f_0(r,x)|^2+|Z_r^{0,x,n}|^2)\rho^{-1}(x)dxdr].
\end{eqnarray*}
Moreover, by Condition (H.1), Lemma \ref{1} and (\ref{au}) we conclude from (\ref{dj}) that
\begin{eqnarray*}\label{ax}
\sup_{n}\int_{\alpha+h'}^{\beta+h'}\int_{\alpha'}^{\beta'}E[\int_\theta^{\theta+h}\int_{\mathbb{R}^d}|Y_{\theta+h}^{0,x,n}-Y_{\theta}^{0,x,n}|^2\rho^{-1}(x)dxdr]d\theta
ds\leq C_ph.
\end{eqnarray*}
Furthermore, by changing the integrations order again and Lemma \ref{26}, we have
\begin{eqnarray*}
\sup_{n}\int_{\alpha+h'}^{\beta+h'}\int_{\alpha'}^{\beta'}E[\int_\theta^{\theta+h}\int_{\mathbb{R}^d}|D_sZ_r^{0,x,n}|^2\rho^{-1}(x)dxdr]d\theta
ds
\leq C_ph.
\end{eqnarray*}
Hence (\ref{dh}) follows. So (\ref{dg}) holds. Now
by (\ref{bc}) and (\ref{bb}) we can deduce that
\begin{eqnarray}\label{bd}
\sup_{n}E[\int_{\alpha+h'}^{\beta+h'}\int_{\alpha'}^{\beta'}\int_{\mathbb{R}^d}|D_{\theta+h}u_n(s,x)-D_{\theta}u_n(s,x)|^2\rho^{-1}(x)dxd\theta ds]
\leq C_ph.
\end{eqnarray}

Now we deal with the second term on the right hand side of
(\ref{ap}). Notice
\begin{eqnarray}\label{bf}
&&\sup_{n}E[\int_{\alpha}^{\beta}\int_{\alpha'}^{\beta'}\int_\mathcal{O}|D_{\theta}u_n(s+h',x)-D_{\theta}u_n(s,x)|^2\rho^{-1}(x)dxd\theta
ds]\nonumber\\
&\leq&\sup_{n}2E[\int_{\alpha}^{\beta}\int_{\alpha'}^{\beta'}\int_{\mathbb{R}^d}|D_{\theta}Y_{s+h'}^{s,x,n}-D_{\theta}Y_s^{s,x,n}|^2\rho^{-1}(x)dxd\theta
ds]\nonumber\\
&&+\sup_{n}2E[\int_{\alpha}^{\beta}\int_{\alpha'}^{\beta'}\int_\mathcal{O}|D_{\theta}Y_{s+h'}^{s+h',x,n}-D_{\theta}Y_{s+h'}^{s,x,n}|^2\rho^{-1}(x)dxd\theta
ds].
\end{eqnarray}
For the first term on the right hand side of (\ref{bf}), by (\ref{aq}), Lemma \ref{qi045} and the exchange of the integrations, it is easy to see that
\begin{eqnarray}\label{bn}
&&\sup_{n}E[\int_{\alpha}^{\beta}\int_{\alpha'}^{\beta'}\int_{\mathbb{R}^d}|D_{\theta}Y_{s+h'}^{s,x,n}-D_{\theta}Y_s^{s,x,n}|^2\rho^{-1}(x)dxd\theta
ds]\nonumber\\
&\leq&C_p\sup_{n}\int_{\alpha}^{\beta}\int_{\alpha'}^{\beta'}E[\int_{\mathbb{R}^d}|D_{\theta}(Y_{s+h'}^{0,x,n}-Y_s^{0,x,n})|^2\rho^{-1}(x)dx]d\theta
ds\nonumber\\
&\leq&C_p\sup_{s\in[0,T-h']}\int_s^{s+h'}(1+\sup_{n}E[\int_{\mathbb{R}^d}|Y_r^{0,x,n}|^{4p-4}\rho^{-1}(x)dx])dr\nonumber\\
&&+C_p\sup_{s\in[0,T-h']}\int_s^{s+h'}(1+\sup_{n}\sup_{\theta\in[0,T]}E[\int_{\mathbb{R}^d}|D_\theta
Y_r^{0,x,n}|^4\rho^{-1}(x)dx])dr\nonumber\\
&&+C_ph'\sup_{n}\sup_{\theta\in[0,T]}E[\int_{0}^{T}\int_{\mathbb{R}^d}|D_\theta
Z_r^{t,x,n}|^2\rho^{-1}(x)dxdr]\nonumber\\
&\leq&C_ph'.
\end{eqnarray}
For the second term on the right hand side of (\ref{bf}), 
firstly from BDSDE (\ref{aq}) we know that
\begin{eqnarray*}
&&D_{\theta}(Y_{s+h'}^{s+h',x,n}-Y_{s+h'}^{s,x,n})\nonumber\\
&=&J(s,s+h')+\int_{s+h'}^{\theta}\partial_yf_n(r,X_r^{s,x},Y_r^{s,x,n})D_{\theta}(Y_{r}^{s+h',x,n}-Y_{r}^{s,x,n})
dr\nonumber\\
&&-\int_{s+h'}^{\theta}\partial_yg(r,X_r^{s,x},Y_r^{s,x,n})D_{\theta}(Y_{r}^{s+h',x,n}-Y_{r}^{s,x,n})d^\dagger\hat{B}_r-\int_{s+h'}^{\theta}D_{\theta}(Z_{r}^{s+h',x,n}-Z_{r}^{s,x,n})dW_r,
\end{eqnarray*}
where
\begin{eqnarray*}
J(s,s+h')
&=&g(\theta,X_{\theta}^{s+h',x},Y_{\theta}^{s+h',x,n})-g(\theta,X_{\theta}^{s,x},Y_{\theta}^{s,x,n})\nonumber\\
&&+\int_{s+h'}^{\theta}\big(\partial_yf_n(r,X_r^{s+h',x},Y_r^{s+h',x,n})-\partial_yf_n(r,X_r^{s,x},Y_r^{s,x,n})\big)D_\theta Y_r^{s+h',x,n}
dr\nonumber\\
&&-\int_{s+h'}^{\theta}\big(\partial_yg(r,X_r^{s+h',x},Y_r^{s+h',x,n})-\partial_yg(r,X_r^{s,x},Y_r^{s,x,n})\big)D_{\theta}Y_r^{s+h',x,n}d^\dagger\hat{B}_r.
\end{eqnarray*}
Applying It$\hat {\rm o}$'s formula to ${\rm
e}^{Kr}|D_{\theta}(Y_{r}^{s+h',x,n}-Y_{r}^{s,x,n})|^2$, 
we have
\begin{eqnarray}\label{bh}
&\sup_nE[\int_\mathcal{O}|D_{\theta}(Y_{s+h'}^{s+h',x,n}-Y_{s+h'}^{s,x,n})|^2\rho^{-1}(x)dx]
&+\sup_nE[\int_{s+h'}^\theta\int_\mathcal{O}|D_{\theta}(Y_{r}^{s+h',x,n}-Y_{r}^{s,x,n})|^2\rho^{-1}(x)dxdr]\nonumber\\
&&+\sup_nE[\int_{s+h'}^\theta\int_\mathcal{O}|D_{\theta}(Z_{r}^{s+h',x,n}-Z_{r}^{s,x,n})|^2\rho^{-1}(x)dxdr]\nonumber\\
&&\leq C_p\sup_nE[\int_\mathcal{O}|J(s,s+h')|^2\rho^{-1}(x)dx].
\end{eqnarray}
So we only need to estimate
$E[\int_\mathcal{O}|J(s,s+h')|^2\rho^{-1}(x)dx]$. Note that by
Condition (H.2)$'$,
\begin{eqnarray}\label{bi}
&&E[\int_\mathcal{O}|J(s,s+h')|^2\rho^{-1}(x)dx]\nonumber\\
&\leq&C_pE[\int_\mathcal{O}|X_{\theta}^{s+h',x}-X_{\theta}^{s,x}|^2\rho^{-1}(x)dx]+C_pE[\int_\mathcal{O}|Y_{\theta}^{s+h',x,n}-Y_{\theta}^{s,x,n}|^2\rho^{-1}(x)dx]\nonumber\\
&&+C_pE[\int_{s+h'}^\theta\int_\mathcal{O}|X_r^{s+h',x}-X_r^{s,x}|^2(1+|Y_r^{s+h',x,n}|^{p-1})^2|D_\theta Y_r^{s+h',x,n}|^2\rho^{-1}(x)dxdr]\nonumber\\
&&+C_pE[\int_{s+h'}^\theta\int_\mathcal{O}|Y_r^{s+h',x,n}-Y_r^{s,x,n}|^2(1+|Y_r^{s+h',x,n}|^{p-2}+|Y_r^{s,x,n}|^{p-2})^2|D_\theta Y_r^{s+h',x,n}|^2\rho^{-1}(x)dxdr]\nonumber\\
&\leq&C_pE[\int_\mathcal{O}|X_{\theta}^{s+h',x}-X_{\theta}^{s,x}|^2\rho^{-1}(x)dx]+C_pE[\int_\mathcal{O}|Y_{\theta}^{s+h',x,n}-Y_{\theta}^{s,x,n}|^2\rho^{-1}(x)dx]\nonumber\\
&&+C_p\int_{s+h'}^\theta\sqrt{E[\int_\mathcal{O}(|Y_r^{s+h',x,n}-Y_r^{s,x,n}|^{4}+|X_r^{s+h',x}-X_r^{s,x}|^{4})\rho^{-1}(x)dx]}\nonumber\\
&&\ \ \ \ \ \ \ \ \ \ \ \ \ \ \times(E[\int_{\mathbb{R}^d}(1+|Y_r^{s+h',x,n}|^{8p-8}+|Y_r^{s,x,n}|^{8p-16})\rho^{-1}(x)dx])^{1\over4}\nonumber\\
&&\ \ \ \ \ \ \ \ \ \ \ \ \ \ \times(E[\int_{\mathbb{R}^d}|D_\theta Y_r^{s+h',x,n}|^8\rho^{-1}(x)dx])^{1\over4}dr.
\end{eqnarray}
From (\ref{a}), we have
\begin{eqnarray*}
X_{r}^{s+h',x}-X_{r}^{s,x}&=&-\int_s^{s+h'}b(X_{u}^{s,x})du-\int_s^{s+h'}\sigma(X_{u}^{s,x})dW_u\nonumber\\
&&+\int_{s+h'}^r\big(b(X_{u}^{s+h',x})-b(X_{u}^{s,x})\big)du+\int_{s+h'}^r\big(\sigma(X_{u}^{s+h',x})-\sigma(X_{u}^{s,x})\big)dW_u.
\end{eqnarray*}
For $q=4$ or $8$, applying It$\hat {\rm o}$'s formula to $|X_{r}^{s+h',x}-X_{r}^{s,x}|^q$, we have
\begin{eqnarray*}
&&E[\int_{\mathbb{R}^d}|X_{r}^{s+h',x}-X_{r}^{s,x}|^q\rho^{-1}(x)dx]\nonumber\\
&\leq&C_pE[\int_{\mathbb{R}^d}{h'}^{q\over2}(\int_s^{s+h'}(1+|X_{u}^{s,x}|)^{2}du)^{q\over2}\rho^{-1}(x)dx]+C_p\int_{\mathbb{R}^d}E[\big(\int_s^{s+h'}|\sigma(X_{u}^{s,x})|^2du\big)^{q\over2}]\rho^{-1}(x)dx\nonumber\\
&&+C_pE[\int_{s+h'}^r\int_{\mathbb{R}^d}|X_{u}^{s+h',x}-X_{u}^{s,x}|^q\rho^{-1}(x)dxdu]\nonumber\\
&\leq&C_p{h'}^{q\over2}+C_pE[\int_{s+h'}^r\int_{\mathbb{R}^d}|X_{u}^{s+h',x}-X_{u}^{s,x}|^q\rho^{-1}(x)dxdu].
\end{eqnarray*}
By Gronwall's inequality, we have for $s+h'\leq r\leq T$,
\begin{eqnarray}\label{bj}
E[\int_{\mathbb{R}^d}|X_{r}^{s+h',x}-X_{r}^{s,x}|^q\rho^{-1}(x)dx]\leq
C_p{h'}^{q\over2}.
\end{eqnarray}
Similarly, noticing (\ref{c}) and applying It$\hat {\rm o}$'s
formula to $|Y_{r}^{s+h',x,n}-Y_{r}^{s,x,n}|^4$, we have
\begin{eqnarray*}
&&E[\int_\mathcal{O}|Y_{r}^{s+h',x,n}-Y_{r}^{s,x,n}|^4\rho^{-1}(x)dx]\nonumber\\
&&+6E[\int_{r}^T\int_\mathcal{O}|Y_{u}^{s+h',x,n}-Y_{u}^{s,x,n}|^2|Z_{u}^{s+h',x,n}-Z_{u}^{s,x,n}|^2\rho^{-1}(x)dxdu]\nonumber\\
&\leq&L{h'}^2+(6L+12L^2)E[\int_{r}^T\int_\mathcal{O}|Y_{u}^{s+h',x,n}-Y_{u}^{s,x,n}|^4\rho^{-1}(x)dxdu]\nonumber\\
&&+{12L}E[\int_{r}^T\int_\mathcal{O}|Y_{u}^{s+h',x,n}-Y_{u}^{s,x,n}|^2(1+|Y_{u}^{s,x,n}|^{2p})|X_{u}^{s+h',x}-X_{u}^{s,x}|^2\rho^{-1}(x)dxdu]\nonumber\\
&&+12L^2E[\int_{r}^T\int_\mathcal{O}|Y_{u}^{s+h',x,n}-Y_{u}^{s,x,n}|^2|X_{u}^{s+h',x}-X_{u}^{s,x}|^2\rho^{-1}(x)dxdu]\nonumber\\
&\leq&L{h'}^2+(2\varepsilon+6L+12L^2)E[\int_{r}^T\int_\mathcal{O}|Y_{r}^{s+h',x,n}-Y_{r}^{s,x,n}|^4\rho^{-1}(x)dxdu]\nonumber\\
&&+C_p\sqrt{\sup_nE[\int_{r}^T\int_{\mathbb{R}^d}(1+|Y_{u}^{s,x,n}|^{8p})\rho^{-1}(x)dxdu]E[\int_{r}^T\int_\mathcal{O}|X_{r}^{s+h',x}-X_{r}^{s,x}|^{8}\rho^{-1}(x)dxdu]}\nonumber\\
&&+C_pE[\int_{r}^T\int_\mathcal{O}|X_{r}^{s+h',x}-X_{r}^{s,x}|^4\rho^{-1}(x)dxdu].
\end{eqnarray*}
Therefore, we can deduce from Lemma \ref{1}, (\ref{bj}) and Gronwall's inequality that, for $s+h'\leq r\leq T$,
\begin{eqnarray}\label{bk}
E[\int_\mathcal{O}|Y_{r}^{s+h',x,n}-Y_{r}^{s,x,n}|^4\rho^{-1}(x)dx]
\leq
C_p{h'}^{2}.
\end{eqnarray}
By (\ref{bj}), (\ref{bk}), Lemmas \ref{1} and
\ref{26}, we know from (\ref{bi}) that
\begin{eqnarray}\label{bl}
\sup_{n}E[\int_\mathcal{O}|J(s,s+h')|^2\rho^{-1}(x)dx]\leq C_ph'.
\end{eqnarray}
Therefore, by (\ref{bh}) and (\ref{bl}) we have
\begin{eqnarray}\label{bm}
&&\sup_{n}E[\int_{\alpha}^{\beta}\int_{\alpha'}^{\beta'}\int_\mathcal{O}|D_{\theta}Y_{s+h'}^{s+h',x,n}-D_{\theta}Y_{s+h'}^{s,x,n}|^2\rho^{-1}(x)dxd\theta
ds]\nonumber\\
&\leq&C_p\int_{\alpha}^{\beta}\int_{\alpha'}^{\beta'}\sup_{n}E[\int_\mathcal{O}|D_{\theta}Y_{s+h'}^{s+h',x,n}-D_{\theta}Y_{s+h'}^{s,x,n}|^2\rho^{-1}(x)dx]d\theta
ds\leq C_ph'.
\end{eqnarray}
Finally, by (\ref{ap}), (\ref{bd})--(\ref{bn}) and (\ref{bm}), (4ii) is satisfied. Theorem \ref{2} is proved. $\hfill\diamond$\\

From the relative compactness of $u_n$ in
$L^2(\Omega\times[0,T];L_\rho^2(\mathcal{O};\mathbb{R}^1))$ for a
bounded domain $\mathcal{O}$ in $\mathbb{R}^d$, we can further prove
that there exists a subsequence of $u_n$, still denoted by $u_n$,
which converges strongly in
$L^2(\Omega\times[0,T];L_\rho^2(\mathbb{R}^d;\mathbb{R}^1))$. We
start from an easy lemma.
\begin{lem}\label{27}
Under Conditions (H.1)--(H.2), (H.3)$^*$, (H.4)--(H.6), for
$u_n(t,x)$ defined in Theorem \ref{2}, we have
$\sup_nE[\int_0^T\int_{\mathbb{R}^d}|u_n(s,x)|^{2p}\rho^{-1}(x)dxds]<\infty$.
Furthermore,
\begin{eqnarray*}
\lim_{N\to\infty}\sup_nE[\int_0^T\int_{\mathbb{R}^d}|u_n(s,x)|I_{{U_N}^c}(x)\rho^{-1}(x)dxds]=0,
\end{eqnarray*}
where ${U_N}=\{x\in\mathbb{R}^d:\ |x|\leq N\}$. 
\end{lem}
{\em Proof}. The claim
$\sup_nE[\int_0^T\int_{\mathbb{R}^d}|u_n(s,x)|^{2p}\rho^{-1}(x)dxds]<\infty$
follows immediately from the equivalence of norm principle
Lemma \ref{qi045} and Lemma \ref{1}. Let's prove the second part of
this lemma. Since $\int_{\mathbb{R}^d}\rho^{-1}(x)dx<\infty$, the claim follows from the following inequality
\begin{eqnarray*}
&&\sup_nE[\int_0^T\int_{\mathbb{R}^d}|u_n(s,x)|^2I_{{U_N}^c}(x)\rho^{-1}(x)dxds]\nonumber\\
&\leq&\big(\sup_nE[\int_0^T\int_{\mathbb{R}^d}|u_n(s,x)|^{2p}\rho^{-1}(x)dxds]\big)^{1\over p}\big(E[\int_0^T\int_{\mathbb{R}^d}|I_{{U_N}^c}(x)|^{p\over{p-1}}\rho^{-1}(x)dxds]\big)^{{p-1}\over p}.
\end{eqnarray*} $\hfill\diamond$

\begin{thm}\label{28}
Under Conditions (H.1)--(H.2), (H.3)$^*$, (H.4)--(H.6), if $(Y_s^{t,x,n},
Z_s^{t,x,n})$ is the solution of BSDEs (\ref{c}) and $Y_s^{t,x}$ is
the weak limit of $Y_s^{t,x,n}$ in
$L^2_\rho(\Omega\times[t,T]\times\mathbb{R}^d;\mathbb{R}^1)$, then
there is a subsequence of $Y_s^{t,x,n}$, still denoted by
$Y_s^{t,x,n}$, converging strongly to $Y_s^{t,x}$ in
$L^2(\Omega\times[t,T];L_\rho^2(\mathbb{R}^d;\mathbb{R}^1))$.
\end{thm}
{\em Proof}. By Theorem \ref{2}, we know that for each bounded
domain $\mathcal{O}\subset\mathbb{R}^d$, there exists a subsequence
of $u_n$ which converges strongly in
$L^2(\Omega\times[0,T];L_\rho^2(\mathcal{O};\mathbb{R}^1))$. So for
$U_1$, we are able to extract a subsequence from $u_n(s,x)$, denoted
by $u_{1n}(s,x)$, which converges strongly in
$L^2(\Omega\times[0,T];L_\rho^2({U_1};\mathbb{R}^1))$. Obviously the
subsequence $u_{1n}(s,x)$ still satisfies the conditions in Theorem
\ref{2}. Applying Theorem \ref{2} again, we are able to extract
another subsequence from $u_{1n}(s,x)$, denoted by $u_{2n}(s,x)$,
that converges strongly in
$L^2(\Omega\times[0,T];L_\rho^2({U_2};\mathbb{R}^1))$. Actually we
can do this procedure for all $U_i$, $i=1,2,\cdot\cdot\cdot$. Now we
pick up the diagonal sequence $u_{ii}(s,x)$,
$i=1,2,\cdot\cdot\cdot$, and still denote this sequence by ${u}_n$
for convenience. It is easy to see that ${u}_n$ converges strongly
in all $L^2(\Omega\times[0,T];L_\rho^2({U_i};\mathbb{R}^1))$,
$i=1,2,\cdot\cdot\cdot$. For arbitrary $\varepsilon>0$, noticing
Lemma \ref{27}, we can find $j(\varepsilon)$ large enough s.t.
\begin{eqnarray*}
\sup_nE[\int_0^T\int_{U^{c}_{j(\varepsilon)}}|{u}_n(s,x)|^2\rho^{-1}(x)dxds]<{\varepsilon\over8}.\
\end{eqnarray*}
For this $j(\varepsilon)$, there exists $n^*(\varepsilon)>0$ s.t.
when $m,n\geq n^*(\varepsilon)$, we know
\begin{eqnarray*}
\|{u}_m-{u}_n\|^2_{L^2(\Omega\times[0,T];L_\rho^2(U_{j(\varepsilon)};\mathbb{R}^1))}=\int_0^T\int_{U_{j(\varepsilon)}}|{u}_m(s,x)-{u}_n(s,x)|^2\rho^{-1}(x)dxds<{\varepsilon\over2}.
\end{eqnarray*}
Therefore as $m,n\geq n^*(\varepsilon)$,
\begin{eqnarray*}
\|{u}_m-{u}_n\|^2_{L^2(\Omega\times[0,T];L_\rho^2(\mathbb{R}^d;\mathbb{R}^1))}
&\leq&E[\int_0^T\int_{U_{j(\varepsilon)}}|{u}_m(s,x)-{u}_n(s,x)|^2\rho^{-1}(x)dxds]\\
&&+E[\int_0^T\int_{U^{c}_{j(\varepsilon)}}(2|{u}_m(s,x)|^2+2|{u}_n(s,x)|^2)\rho^{-1}(x)dxds]<\varepsilon.
\end{eqnarray*}
That is to say ${u}_n$ converges strongly in
$L^2(\Omega\times[0,T];L_\rho^2(\mathbb{R}^d;\mathbb{R}^1))$. Then the strong convergence of $Y_\cdot^{t,\cdot,n}$ follows from the standard equivalence of norm principle argument.
On the other hand,  $Y_s^{t,x,n}$ is also weakly
convergent in
$L^2(\Omega\times[t,T];L_\rho^2(\mathbb{R}^d;\mathbb{R}^1))$ with
the weak limit $Y_s^{t,x}$. Therefore $Y_s^{t,x,n}$ converges
strongly to $Y_s^{t,x}$ in
$L^2(\Omega\times[t,T];L_\rho^2(\mathbb{R}^d;\mathbb{R}^1))$. $\hfill\diamond$\\

Considering the strongly convergent subsequence
$\{Y_\cdot^{t,\cdot,n}\}_{n=1}^\infty$ derived from Theorem \ref{28} and using B-D-G inequality to BDSDE
(\ref{c}), we can prove that for arbitrary $m,n$,
\begin{eqnarray}\label{by}
&&E[\sup_{s\in[t,T]}\int_{\mathbb{R}^d}|Y_s^{t,x,m}-Y_s^{t,x,n}|^2\rho^{-1}(x)dx]+E[\int_t^T\int_{\mathbb{R}^d}|Z_r^{t,x,m}-Z_r^{t,x,n}|^2\rho^{-1}(x)dxdr]\nonumber\\
&\leq&
C_p\big(E[\int_t^T\int_{\mathbb{R}^d}|Y_r^{t,x,m}-Y_r^{t,x,n}|^2\rho^{-1}(x)dxdr]\nonumber\\
&&\ \ \ \ \ \times E[\int_t^T\int_{\mathbb{R}^d}(1+|f_0(r,x)|^2+|Y_r^{t,x,n}|^{2p}+|Y_r^{t,x,m}|^{2p})\rho^{-1}(x)dxdr]\big)^{1\over2}\nonumber\\
&&+C_pE[\int_t^T\int_{\mathbb{R}^d}|g(r,X_{r}^{t,x},{Y}_{r}^{t,x,m})-g(r,X_{r}^{t,x},{Y}_{r}^{t,x,n})|^2\rho^{-1}(x)dxdr].
\end{eqnarray}
Using strong
subsequence of $Y_\cdot^{t,\cdot,n}$ and the Lipschitz continuity of $g$, by the dominated convergence
theorem we can conclude from (\ref{by}) that this subsequence $\{Y_\cdot^{t,\cdot,n}\}_{n=1}^\infty$
converges strongly also in $S^{2,0}([t,T];L_{\rho}^2({\mathbb{R}^{d}};{\mathbb{R}^{1}}))$ and the
corresponding subsequence of $\{Z_\cdot^{t,\cdot,n}\}_{n=1}^\infty$
converges strongly $M^{2,0}([t,T];L_{\rho}^2({\mathbb{R}^{d}};{\mathbb{R}^{d}}))$ as well.
Certainly the strong convergence limit should be identified with the
weak convergence limit $Z_\cdot^{t,\cdot}$. Hence the following
corollary follows without a surprise.
\begin{cor}\label{18}
Let $(Y_\cdot^{t,\cdot},Z_\cdot^{t,\cdot})$ be the solution of BDSDE
(\ref{d}) and $(Y_\cdot^{t,\cdot,n},Z_\cdot^{t,\cdot,n})$ be the
subsequence of the solutions of BDSDE (\ref{c}), of which
$Y_\cdot^{t,\cdot,n}$ converges strongly to $Y_\cdot^{t,\cdot}$ in
$L^2(\Omega\times[t,T];L_\rho^2(\mathbb{R}^d;\mathbb{R}^1))$, then $(Y_\cdot^{t,\cdot,n},Z_\cdot^{t,\cdot,n})$ also converges strongly to $(Y_\cdot^{t,\cdot},Z_\cdot^{t,\cdot})$
in $S^{2,0}([t,T];L_{\rho}^2({\mathbb{R}^{d}};{\mathbb{R}^{1}}))\times M^{2,0}([t,T];L_{\rho}^2({\mathbb{R}^{d}};{\mathbb{R}^{d}}))$.
\end{cor}

As for $Y_s^{t,x}$, we further have
\begin{lem}\label{33}
Under Conditions (H.1)--(H.2), (H.3)$^*$, (H.4)--(H.6),
$E[\int_t^T\int_{\mathbb{R}^d}|Y_s^{t,x}|^{2p}\rho^{-1}(x)dxds]<\infty$ and $Y_s^{t,x}=Y_s^{s,X_s^{t,x}}$ for a.e. $s\in[t,T]$, $x\in\mathbb{R}^d$ a.s.
\end{lem}
{\em Proof}. First by Lemma \ref{qi045} and Corollary \ref{18}, we have
\begin{eqnarray*}
&&E[\int_t^T\int_{\mathbb{R}^d}|Y_s^{t,x}-Y_s^{s,X_s^{t,x}}|^2\rho^{-1}(x)dxds]\nonumber\\
&\leq&\lim_{n\to\infty}2E[\int_t^T\int_{\mathbb{R}^d}|Y_s^{t,x,n}-Y_s^{t,x}|^2\rho^{-1}(x)dxds]
+\lim_{n\to\infty}2E[\int_t^T\int_{\mathbb{R}^d}|Y_s^{s,X_s^{t,x},n}-Y_s^{s,X_s^{t,x}}|^2\rho^{-1}(x)dxds]\nonumber\\
&\leq&\lim_{n\to\infty}2E[\int_t^T\int_{\mathbb{R}^d}|Y_s^{t,x,n}-Y_s^{t,x}|^2\rho^{-1}(x)dxds]
+\lim_{n\to\infty}C_pE[\sup_{s\leq r\leq T}\int_{\mathbb{R}^d}|Y_r^{s,x,n}-Y_r^{s,x}|^2\rho^{-1}(x)dx]=0.
\end{eqnarray*}
Hence,
\begin{eqnarray}\label{de}
Y_s^{t,x}=Y_s^{s,X_s^{t,x}}\ {\rm for}\ {\rm a.e.}\ s\in[t,T],\
x\in\mathbb{R}^d\ {\rm a.s.}
\end{eqnarray}
If we define $Y_s^{s,x}=u(s,x)$, then by (\ref{de}) and Lemma \ref{qi045} again we also have
\begin{eqnarray}\label{s}
\lim_{n\to\infty}E[\int_0^T\int_{\mathbb{R}^d}|u_n(s,x)-u(s,x)|^2\rho^{-1}(x)dxds]=0,
\end{eqnarray}
and
\begin{eqnarray*}\label{q}
Y_s^{t,x}=u(s,X_s^{t,x})\ {\rm for}\ {\rm a.e.}\ s\in[t,T],\ x\in\mathbb{R}^d\
{\rm a.s.}
\end{eqnarray*}
Therefore, we claim that the strong limit of $u_n(s,x)$ in
$L^2(\Omega\times[0,T];L_\rho^2(\mathbb{R}^d;\mathbb{R}^1))$ is $u(s,x)$. 

By the equivalence of norm principle, to get
$E[\int_t^T\int_{\mathbb{R}^d}|Y_s^{t,x}|^{2p}\rho^{-1}(x)dxds]<\infty$, we only need to
prove
$E[\int_0^T\int_{\mathbb{R}^d}|u(s,x)|^{2p}\rho^{-1}(x)dxds]<\infty$.
For this, we first derive from $\lim_{n\to\infty}E[\int_{0}^{T}\int_{\mathbb{R}^{d}}|u_n(s,x)-u(s,x)|^2\rho^{-1}(x)dxds]=0$ a subsequence of
$\{u_n(s,x)\}_{n=1}^\infty$, still denoted by $\{u_n(s,x)\}_{n=1}^\infty$, s.t.
\begin{eqnarray*}\label{ao}
u_n(s,x)\longrightarrow u(s,x)\  {\rm and}\
\sup_n|u_n(s,x)|^{2p}<\infty\ \ {\rm for}\ {\rm a.e.}\ s\in[t,T],\
x\in\mathbb{R}^{d}\ {\rm a.s.}
\end{eqnarray*}

By a similar argument
as in Lemma \ref{27}, for this subsequence $u_n$, we can prove,
using H$\ddot{o}$lder inequality, that for any $\delta>0$,
\begin{eqnarray*}\label{cl}
\lim_{N\to\infty}\sup_nE[\int_0^T\int_{\mathbb{R}^d}|u_n(s,x)|^{2p-\delta}I_{\{|u_n(s,x)|^{2p-\delta}>N\}}(s,x)\rho^{-1}(x)dxds]=0.
\end{eqnarray*}
That is to say that $|u_n(s,x)|^{2p-\delta}$ is uniformly
integrable. Moreover by the fact that $u_n(s,x)\longrightarrow
u(s,x)$ for a.e. $s\in[0,T]$, $x\in\mathbb{R}^{d}$ a.s. and
generalized Lebesgue's dominated convergence theorem \cite{ro},
we have
\begin{eqnarray*}
E[\int_0^T\int_{\mathbb{R}^d}|u(s,x)|^{2p-\delta}\rho^{-1}(x)dxds]&=&\lim_{n\to\infty}E[\int_0^T\int_{\mathbb{R}^d}|u_n(s,x)|^{2p-\delta}\rho^{-1}(x)dxds]\\
&\leq&C_p\big(\sup_{n}E[\int_0^T\int_{\mathbb{R}^d}|u_n(s,x)|^{2p}\rho^{-1}(x)dxds]\big)^{{2p-\delta}\over{2p}}\leq
C_p,
\end{eqnarray*}
where the last $C_p<\infty$ is a constant independent of $n$ and
$\delta$. Then using Fatou lemma to take the limit as $\delta\to0$
in the above inequality, we can get
$E[\int_0^T\int_{\mathbb{R}^d}|u(s,x)|^{2p}\rho^{-1}(x)dxds]<\infty$.
$\hfill\diamond$\\

Indeed, with Corollary
\ref{18} and Lemma \ref{33}, using It$\hat {\rm o}$'s formula to
${\rm
e}^{Kr}|Y_{r}^{t,x}|^{2p}$, we can further prove that
$Y_\cdot^{t,\cdot}\in
S^{2p,0}([t,T];L_{\rho}^{2p}({\mathbb{R}^{d}};{\mathbb{R}^{1}}))$
(To see similar calculations, one can refer to the
proof of Lemma \ref{1} in Section 7). 
\begin{prop}\label{19}
For $(Y_\cdot^{t,\cdot},Z_\cdot^{t,\cdot})$ and
$(Y_\cdot^{t,\cdot,n},Z_\cdot^{t,\cdot,n})$ given in Corollary
\ref{18}, $Y_\cdot^{t,\cdot}\in
S^{2p,0}([t,T];L_{\rho}^{2p}({\mathbb{R}^{d}};{\mathbb{R}^{1}}))$.
\end{prop}

Now we are ready to prove the identification of the limiting BDSDEs.
\begin{lem}\label{17}
The random field $U$, $V$, $Y$ and $Z$ have the following relation:
\begin{eqnarray*}\label{ak}
U_s^{t,x}=f(s,X_s^{t,x},Y_s^{t,x}),\
V_s^{t,x}=g(s,X_s^{t,x},Y_s^{t,x})\ {\rm for}\ {\rm a.e.}\
s\in[t,T],\ x\in\mathbb{R}^d\ {\rm a.s.}\ \ \ \
\end{eqnarray*}
\end{lem}
{\em Proof}. Let $\mathcal{K}$ be a set in
$\Omega\times[t,T]\times\mathbb{R}^d$ s.t.
$\sup_n|Y_s^{t,x,n}|+\sup_n|Z_s^{t,x,n}|+|f_0(s,X_s^{t,x})|+|g(s,X_s^{t,x},0)|<K$.
First we can find a subsequence of
$(Y_s^{t,x,n},Z_s^{t,x,n})$, still denoted by
$(Y_s^{t,x,n},Z_s^{t,x,n})$, satisfying
$(Y_s^{t,x,n},Z_s^{t,x,n})\longrightarrow (Y_s^{t,x},Z_s^{t,x})$
a.s. and $\sup_n|Y_s^{t,x,n}|+\sup_n|Z_s^{t,x,n}|<\infty$ for a.e.
$s\in[t,T]$, $x\in\mathbb{R}^{d}$ a.s.
Then it turns out that as $K\to\infty$,
$\mathcal{K}\uparrow\Omega\times[t,T]\times\mathbb{R}^d$. 
Moreover it is easy to see that for this subsequence,
\begin{eqnarray*}
&&E[\int_t^T\int_{\mathbb{R}^d}2(\sup_n|f_n(s,X_s^{t,x},Y_s^{t,x,n})|^2+|f(s,X_s^{t,x},Y_s^{t,x})|^2)I_\mathcal{K}(s,x)\rho^{-1}(x)dxds]\nonumber\\
&\leq&C_pE[\int_t^T\int_{\mathbb{R}^d}(|f_0(s,X_s^{t,x})|^2+\sup_n|Y_s^{t,x,n}|^{2p})I_\mathcal{K}(s,x)\rho^{-1}(x)dxds]\nonumber\\
&&+C_pE[\int_t^T\int_{\mathbb{R}^d}(|f_0(s,X_s^{t,x})|^2+|Y_s^{t,x}|^{2p})I_\mathcal{K}(s,x)\rho^{-1}(x)dxds]<\infty.
\end{eqnarray*}
Thus, we can apply Lebesgue's dominated convergence theorem to the
following estimate:
\begin{eqnarray}\label{zz17}
&&\lim_{n\to\infty}E[\int_t^T\int_{\mathbb{R}^d}|f_n(s,X_s^{t,x},Y_s^{t,x,n})I_\mathcal{K}(s,x)-f(s,X_s^{t,x},Y_s^{t,x})I_\mathcal{K}(s,x)|^2\rho^{-1}(x)dxds]\nonumber\\
&\leq&2E[\int_t^T\int_{\mathbb{R}^d}\lim_{n\to\infty}|f_n(s,X_s^{t,x},Y_s^{t,x,n})-f(s,X_s^{t,x},Y_s^{t,x,n})|^2I_\mathcal{K}(s,x)\rho^{-1}(x)dxds]\nonumber\\
&&+2E[\int_t^T\int_{\mathbb{R}^d}\lim_{n\to\infty}|f(s,X_s^{t,x},Y_s^{t,x,n})-f(s,X_s^{t,x},Y_s^{t,x})|^2I_\mathcal{K}(s,x)\rho^{-1}(x)dxds].
\end{eqnarray}
Since $Y_s^{t,x,n}\longrightarrow Y_s^{t,x}$ for a.e. $s\in[t,T]$,
$x\in\mathbb{R}^{d}$ a.s., there exists a $N(s,x,\omega)$ s.t. when
$n\geq N(s,x,\omega)$, $|Y_s^{t,x,n}|\leq|Y_s^{t,x}|+1$.
We take $f_n(s,x,y)$ to be the
smootherized truncations to $f(s,x,y)$ on
variable $y$, so taking
$n\geq\max\{N(s,x,\omega),\ |Y_s^{t,x}|+1\}$, we have $f_n(s,X_s^{t,x},Y_s^{t,x,n})$ $=f(s,X_s^{t,x},{\inf(n,|Y_s^{t,x,n}|)\over|Y_s^{t,x,n}|}Y_s^{t,x,n})+\partial_yf(s,x,{n\over{|Y_s^{t,x,n}|}}Y_s^{t,x,n})(Y_s^{t,x,n}-{n\over{|Y_s^{t,x,n}|}}Y_s^{t,x,n})I_{\{|Y_s^{t,x,n}|>n\}}=f(s,X_s^{t,x},Y_s^{t,x,n})$.
That is to say
$$\lim_{n\to\infty}|f_n(s,X_s^{t,x},Y_s^{t,x,n})-f(s,X_s^{t,x},Y_s^{t,x,n})|^2=0\
{\rm for}\ {\rm a.e.}\ s\in[t,T],\ x\in\mathbb{R}^{d}\ {\rm a.s.}$$
On the other hand,
$\lim_{n\to\infty}|f(s,X_s^{t,x},Y_s^{t,x,n})-f(s,X_s^{t,x},Y_s^{t,x})|^2=0$
for a.e. $s\in[t,T]$, $x\in\mathbb{R}^{d}$ a.s. is obvious due to
the continuity of $y\longrightarrow f(s,x,y)$.

Therefore by (\ref{zz17}),
$f_n(s,X_s^{t,x},Y_s^{t,x,n})I_\mathcal{K}(s,x)=U_s^{t,x,n}I_\mathcal{K}(s,x)$
converges strongly to \\$f(s,X_s^{t,x},Y_s^{t,x})I_\mathcal{K}(s,x)$
in $L^2_\rho(\Omega\times[t,T]\times\mathbb{R}^d;\mathbb{R}^1)$, but
$U_s^{t,x,n}I_\mathcal{K}(s,x)$ converges weakly to
$U_s^{t,x}I_\mathcal{K}(s,x)$ in
$L^2_\rho(\Omega\times[t,T]\times\mathbb{R}^d;\mathbb{R}^1)$, so
$f(s,X_s^{t,x},Y_s^{t,x})I_\mathcal{K}(s,x)=U_s^{t,x}I_\mathcal{K}(s,x)$
for a.e. $r\in[t,T]$, $x\in\mathbb{R}^{d}$ a.s. Taking $K\to\infty$,
we obtain the first part of Lemma \ref{17}. The other part of Lemma \ref{17} can be proved similarly. $\hfill\diamond$\\

{\em Proof of Theorem \ref{21}}. With Proposition
\ref{19} and Lemma \ref{17}, the existence of solution of BDSDE (\ref{co}) is easy to
see. The uniqueness can be proved using a standard substraction, It$\hat {\rm
o}$'s formula and Gronwall inequality argument. Here the monotonicity plays an important role.
$\hfill\diamond$\\

By the stochastic flow property $X_r^{s,X^{t,x}_s}=X^{t,x}_r$ for $t\leq s\leq r\leq T$ and the uniqueness of solution of BDSDE (\ref{co}), following a similar argument as Proposition 3.4 in \cite{zh-zh1} we have
\begin{cor}\label{4}
Under the conditions of Theorem \ref{21}, let $(Y_s^{t,x},Z_s^{t,x})$ be the solution of BDSDE
(\ref{co}), then
\begin{eqnarray*}\label{v}
Y_s^{t,x}=Y_s^{s,X_s^{t,x}},\ Z_s^{t,x}=Z_s^{s,X_s^{t,x}}\ \ {\rm
for}\ {\rm any}\ s\in[t,T],\ {\rm a.e.}\ x\in\mathbb{R}^d\ {\rm
a.s.}
\end{eqnarray*}
\end{cor}

Naturally, we can relate the $S^{2p,0}([t,T],L^{2p}_\rho(\mathbb{R}^d;\mathbb{R}^1))\times M^{2,0}([t,T],L^2_\rho(\mathbb{R}^d;\mathbb{R}^d))$ solution of BDSDE (\ref{co}) to the weak solution of SPDE (\ref{bz})
with finite dimensional noise. We have the following theorem: 
\begin{thm}\label{31}
Define $u(t,x)=Y_t^{t,x}$, where $(Y_s^{t,x},Z_s^{t,x})$ is the
solution of BDSDE (\ref{co}) under Conditions (H.1)--(H.2),
(H.3)$^*$, (H.4)--(H.6), then $u(t,x)$ is the unique weak solution
of SPDE (\ref{bz}) with finite dimensional noise. Moreover, let $u$ be a representative in the equivalence class of the solution of the SPDE (\ref{bz}) in $S^{2p,0}([t,T];L_{\rho}^{2p}({\mathbb{R}^{d}};{\mathbb{R}^{1}}))$ with $\sigma^*\nabla u\in
M^{2,0}([t,T];L_{\rho}^2({\mathbb{R}^{d}};{\mathbb{R}^{d}}))$, then $u(t,x)=Y_t^{t,x}$ for a.e. $t\in[0,T]$, $x\in\mathbb{R}^d$  a.s. and
\begin{eqnarray}\label{dd}
u(s,X_s^{t,x})=Y_s^{t,x},\ (\sigma^*\nabla u)(s,X^{t,x}_s)=Z_s^{t,x}
\ {\rm for}\ {\rm a.e.}\ s\in[t,T],\ x\in\mathbb{R}^{d}\ {\rm a.s.}
\end{eqnarray}
\end{thm}
{\em Proof}. Using Corollary \ref{18}, we first prove the
relationship between $(Y,Z)$ and $u$, when we take $u(t,x)=Y_t^{t,x}$. Having proved Lemma \ref{33}, we only need
to prove that $(\sigma^*\nabla u)(s,X^{t,x}_s)=Z_s^{t,x}$ for a.e.
$s\in[t,T]$, $x\in\mathbb{R}^{d}$ a.s. This can be deduced from
(\ref{v}) and the strong convergence of $Z_\cdot^{t,\cdot,n}$ to
$Z_\cdot^{t,\cdot}$ in
$L^2(\Omega\times[t,T];L_\rho^2(\mathbb{R}^d;\mathbb{R}^1))$ by the
similar argument as in Proposition 4.2 in \cite{zh-zh1}.

We then prove that $u(t,x)$ defined above is the unique weak solution of SPDE
(\ref{bz}). We start our proof from smoothed SPDE (\ref{f}). Let
$u_n(s,x)$ be the weak solution of SPDE (\ref{f}), then
$(u_n,\sigma^*\nabla u_n)\in
L^{2}(\Omega\times[0,T];L_{\rho}^2({\mathbb{R}^{d}};{\mathbb{R}^{1}}))\times
L^{2}(\Omega\times[0,T];L_{\rho}^2({\mathbb{R}^{d}};{\mathbb{R}^{d}}))$
and for an arbitrary $\varphi\in
C_c^{\infty}(\mathbb{R}^d;\mathbb{R}^1)$,
\begin{eqnarray}\label{qi16}
&&\int_{\mathbb{R}^{d}}u_n(t,x)\varphi(x)dx-\int_{\mathbb{R}^{d}}h(x)\varphi(x)dx-{1\over2}\int_{t}^{T}\int_{\mathbb{R}^{d}}(\sigma^*\nabla u_n)(s,x)(\sigma^*\nabla\varphi)(x)dxds\nonumber\\
&&-\int_{t}^{T}\int_{\mathbb{R}^{d}}u_n(s,x)div\big((b-\tilde{A})\varphi\big)(x)dxds\\
&=&\int_{t}^{T}\int_{\mathbb{R}^{d}}f_n\big(s,x,u_n(s,x)\big)\varphi(x)dxds-\int_{t}^{T}\int_{\mathbb{R}^{d}}\langle
g\big(s,x,u_n(s,x)\big)\varphi(x)dx,d^\dagger\hat{B}_s\rangle.\nonumber
\end{eqnarray}
We can prove that along a subsequence each term of (\ref{qi16})
converges weakly to the corresponding term of (\ref{ae}) in
$L^2(\Omega;\mathbb{R}^1)$. By (\ref{s}), we know that $u_n$
converges strongly to $u$ in
$L^2_\rho(\Omega\times[0,T]\times\mathbb{R}^d;\mathbb{R}^1)$, thus 
$u_n$ also converges weakly. Moreover, note
$\sup_{x\in\mathbb{R}^d}(|div\big((b-\tilde{A})\varphi\big)(x)|)<\infty$
and $\rho$ is a continuous function in $\mathbb{R}^d$. So it is
obvious that in the sense of the weak convergence in
$L^2(\Omega;\mathbb{R}^d)$,
\begin{eqnarray*}
\lim_{n\to\infty}\int_{t}^{T}\int_{\mathbb{R}^{d}}u_n(s,x)div\big((b-\tilde{A})\varphi\big)(x)dxds=\int_{t}^{T}\int_{\mathbb{R}^{d}}u(s,x)div\big((b-\tilde{A})\varphi\big)(x)dxds.
\end{eqnarray*}
Also it is easy to see that
\begin{eqnarray*}
\lim_{n\to\infty}{1\over2}\int_{t}^{T}\int_{\mathbb{R}^{d}}(\sigma^*\nabla
u_n)(s,x)(\sigma^*\nabla\varphi)(x)dxds
&=&-{1\over2}\int_{t}^{T}\int_{\mathbb{R}^{d}}u(s,x)\nabla(\sigma\sigma^*\nabla\varphi\sigma)(x)\rho(x)\rho^{-1}(x)dxds\nonumber\\&=&{1\over2}\int_{t}^{T}\int_{\mathbb{R}^{d}}(\sigma^*\nabla
u)(s,x)(\sigma^*\nabla\varphi)(x)dxds.
\end{eqnarray*}
Note that $f_n(s,X_s^{t,x},Y_s^{t,x,n})$
converges weakly to $f(s,X_s^{t,x},Y_s^{t,x})$ in
$L_\rho^2(\Omega\times[t,T]\times\mathbb{R}^d;\mathbb{R}^1)$. In
fact we can use the same procedures to prove that
$f_n\big(s,x,u_n(s,x)\big)$ converges weakly to
$f\big(s,x,u(s,x)\big)$ and $g\big(s,x,u_n(s,x)\big)$ converges
weakly to $g\big(s,x,u(s,x)\big)$ in
$L_\rho^2(\Omega\times[t,T]\times\mathbb{R}^d;\mathbb{R}^1)$. Then
following the proof of BDSDE (\ref{c}) converging weakly to BDSDE
(\ref{d}) and taking weak limit here in $L^2(\Omega;\mathbb{R}^1)$,
we obtain the weak convergence of three terms:
\begin{eqnarray*}
&&\lim_{n\to\infty}\int_t^T\int_{\mathbb{R}^d}f_n\big(s,x,u_n(s,x)\big)\varphi(x)dxds-\lim_{n\to\infty}\int_{t}^{T}\int_{\mathbb{R}^{d}}\langle
g\big(s,x,u_n(s,x)\big)\varphi(x)dx,d^\dagger\hat{B}_s\rangle\\
&=&\int_t^T\int_{\mathbb{R}^d}f\big(s,x,u(s,x)\big)\varphi(x)dxds-\int_{t}^{T}\int_{\mathbb{R}^{d}}\langle
g\big(s,x,u(s,x)\big)\varphi(x)dx,d^\dagger\hat{B}_s\rangle.
\end{eqnarray*}
Finally, that for any $t\in[0,T]$,
$\lim_{n\to\infty}\int_{\mathbb{R}^{d}}u_n(t,x)\varphi(x)dx=\int_{\mathbb{R}^{d}}u(t,x)\varphi(x)dx$
in the sense of weak convergence in $L^2(\Omega;\mathbb{R}^1)$ can
be deduced from Corollary \ref{18}:
\begin{eqnarray*}
\lim_{n\to\infty}|E[\int_{\mathbb{R}^{d}}(u_n(t,x)-u(t,x))\varphi(x)dx]|^2
&\leq&\lim_{n\to\infty}C_pE[\int_{\mathbb{R}^{d}}|u_n(t,X_t^{0,x})-u(t,X_t^{0,x})|^2\rho^{-1}(x)dx]\nonumber\\
&\leq&\lim_{n\to\infty}C_pE[\sup_{0\leq t\leq
T}\int_{\mathbb{R}^{d}}|Y_t^{0,x,n}-Y_t^{0,x}|^2\rho^{-1}(x)dx]=0.\nonumber
\end{eqnarray*}
Here the convergence in the $S^{2p,0}$ space gives us a strong result about the convergence $\int_{\mathbb{R}^d}u_n(t,x)\varphi(x)dx \linebreak
\longrightarrow \int_{\mathbb{R}^d}u(t,x)\varphi(x)dx$ in $L^2(\Omega;\mathbb{R}^d)$ uniformly in $t$ as $n\to\infty$. Therefore we prove that (\ref{ae}) is satisfied for all $t\in[0,T]$,
hence $u(t,x)$ is a weak solution of SPDE (\ref{bz}).

The uniqueness of weak solution of SPDE (\ref{bz}) can be derived
from the uniqueness of solution of BDSDE (\ref{spdespgrowth2}). For this, let
$u$ be a weak solution of SPDE (\ref{bz}). Define
$F(s,x)=f\big(s,x,u(s,x)\big)$ and $G(s,x)=g\big(s,x,u(s,x)\big)$.
Since $u$ is the solution, so
$\int_{0}^{T}\int_{\mathbb{R}^d}\big(|u(s,x)|^{2p}+|(\sigma^*\nabla
u)(s,x)|^2\big)\rho^{-1}(x)dxds<\infty$ and
\begin{eqnarray*}\label{zz14}
&&E[\int_{0}^{T}\int_{\mathbb{R}^d}\big(|F(s,x)|^2+|G(s,x)|^2\big)\rho^{-1}(x)dxds]\nonumber\\
&\leq&C_pE[\int_{0}^{T}\int_{\mathbb{R}^d}\big(1+|f_0(s,x)|^2+|g(s,x,0)|^2+|u(s,x)|^{2p}\big)\rho^{-1}(x)dxds]<\infty.\nonumber
\end{eqnarray*}
Then we get a SPDE with the generator $(F,G)\in L^{2}(\Omega\times[0,T];L_{\rho}^2({\mathbb{R}^{d}};{\mathbb{R}^{1}}))\times L^{2}(\Omega\times[0,T];L_{\rho}^2({\mathbb{R}^{d}};{\mathbb{R}^{l}}))$. For this generator $(F,G)$, we claim that $(Y_s^{t,x},Z_s^{t,x})\triangleq (u(s,X_s^{t,x}),(\sigma^*\nabla u)(s,X^{t,x}_s))$ solves the following linear BDSDE for a.e. $x\in\mathbb{R}^d$ with probability one:
\begin{eqnarray}\label{zz15}
Y_s^{t,x}=h(X_{T}^{t,x})+\int_{s}^{T}F(r,X_{r}^{t,x})dr-\int_{s}^{T}\langle
G(r,X_{r}^{t,x}),d^\dagger\hat{B}_r\rangle-\int_{s}^{T}\langle
Z^{t,x}_r,dW_r\rangle.
\end{eqnarray}
First we use the mollifier to smootherize $(h,F,G)$, then we get a smootherized sequence $(h^m,F^m,G^m)$ such that $(h^m(\cdot),F^m(s,\cdot),G^m(s,\cdot))\longrightarrow (h(\cdot),F(s,\cdot),G(s,\cdot))$ in
$L_{\rho}^2({\mathbb{R}^{d}};{\mathbb{R}^{1}})\times L_{\rho}^2({\mathbb{R}^{d}};{\mathbb{R}^{1}})\times L_{\rho}^2({\mathbb{R}^{d}};{\mathbb{R}^{1}})$. Denote
by ${u}_m(t,x)$ the solution
of SPDE on $[0,T]$ with terminal value $h^m(x)$ and generator \linebreak $(F^m(s,x),G^m(s,x))$ and by $(Y_{s,m}^{t,x},Z_{s,m}^{t,x})$ the solution of BDSDE with terminal value $h^m(X_T^{t,x})$ and generator $(F^m(s,X_s^{t,x}),G^m(s,X_s^{t,x}))$, then following classical results of Pardoux and Peng \cite{pa-pe3}, we
have ${Z}_{t,m}^{t,x}=\sigma^*\nabla{u}_m(t,x)$, and ${Y}_{s,m}^{t,x}={u}_m(s,X_s^{t,x})={Y}_{s,m}^{s,X_s^{t,x}}$,
${Z}_{s,m}^{t,x}=\sigma^*\nabla{u}_m(s,X_s^{t,x})={Z}_{s,m}^{s,X_s^{t,x}}$.
But by standard estimates $({Y}_{s,m}^{t,x},{Z}_{s,m}^{t,x})$ is a Cauchy sequence in $M^{2,0}([0,T];L_{\rho}^2({\mathbb{R}^{d}};{\mathbb{R}^{1}}))\times
M^{2,0}([0,T];L_{\rho}^2({\mathbb{R}^{d}};{\mathbb{R}^{d}}))$. By equivalence of norm principle, $u_m(s,x)$ is also a Cauchy sequence in $\mathcal{H}$, where $\mathcal{H}$ is the set of random fields $\{w(s,x);\
s\in[0,T],\ x\in\mathbb{R}^{d}\}$ such that $(w,\sigma^*\nabla w)\in
M^{2,0}([0,T];L_{\rho}^2({\mathbb{R}^{d}};{\mathbb{R}^{1}}))\times
M^{2,0}([0,T];L_{\rho}^2({\mathbb{R}^{d}};{\mathbb{R}^{d}}))$ with
the norm
$
\left \{E[\int_{0}^{T}\int_{\mathbb{R}^{d}}(|w(s,x)|^2
\right . \linebreak
\left .+|(\sigma^*\nabla)w(s,x)|^2)\rho^{-1}(x)dxds)]\right \}^{1\over 2}<\infty.
$
 So there exists ${u}\in\mathcal{H}$ such that
$({u}_m,\sigma^*\nabla
{u}_m)\rightarrow({u},\sigma^*\nabla{u})$
in
$M^{2,0}([0,T];L_{\rho}^2({\mathbb{R}^{d}};{\mathbb{R}^{1}}))\times
M^{2,0}([0,T];L_{\rho}^2({\mathbb{R}^{d}};{\mathbb{R}^{d}}))$ due to the completeness of $\mathcal{H}$. By the equivalence of norm principle again, we know that $(Y_s^{t,x},Z_s^{t,x})\triangleq (u(s,X_s^{t,x}),(\sigma^*\nabla u)(s,X^{t,x}_s))$ is the limit of Cauchy sequence of $({Y}_{s,m}^{t,x},{Z}_{s,m}^{t,x})$. Now it is
easy to pass the limit as $m\rightarrow\infty$ on the BDSDE which $({Y}_{s,m}^{t,x},{Z}_{s,m}^{t,x})$ satisfies
and conclude that $(Y_s^{t,x},Z_s^{t,x})$ is a solution
of BDSDE (\ref{zz15}).

Noting the
definition of $F(s,x)$, $G(s,x)$, $Y_s^{t,x}$ and $Z_s^{t,x}$, we have that
$(Y_s^{t,x},Z_s^{t,x})$ solves
BDSDE (\ref{spdespgrowth2}) for a.e. $x\in\mathbb{R}^d$ with probability one.
Moreover, by Lemma \ref{qi045},
\begin{eqnarray*}
&&E[\int_{t}^{T}\int_{\mathbb{R}^d}(|Y_s^{t,x}|^{2p}+|Z_s^{t,x}|^2)\rho^{-1}(x)dxds]
\leq C_p\int_{t}^{T}\int_{\mathbb{R}^d}|u(s,x)|^{2p}+|(\sigma^*\nabla
u)(s,x)|^2\rho^{-1}(x)dxds<\infty.
\end{eqnarray*}
As Proposition \ref{19}, we can further deduce
that $Y_\cdot^{t,\cdot}\in
S^{2p,0}([t,T];L_{\rho}^{2p}({\mathbb{R}^{d}};{\mathbb{R}^{1}}))$ and
therefore $(Y_s^{t,x},Z_s^{t,x})$ is a solution of BDSDE (\ref{spdespgrowth2}). If there is another solution $\hat{u}$ to
SPDE (\ref{bz}), then similarly we can find another solution
$(\hat{Y}_s^{t,x},\hat{Z}_s^{t,x})$ to BDSDE (\ref{spdespgrowth2}),
where
\begin{eqnarray*}
\hat{Y}_s^{t,x}=\hat{u}(s,X_s^{t,x})\ {\rm and}\
\hat{Z}_s^{t,x}=(\sigma^*\nabla \hat{u})(s,X^{t,x}_s).
\end{eqnarray*}
By Theorem \ref{21}, the solution of BDSDE (\ref{spdespgrowth2}) is
unique, therefore
\begin{eqnarray*}
Y_s^{t,x}=\hat{Y}_s^{t,x}\ {\rm for}\ {\rm a.e.}\ s\in[t,T],\
x\in\mathbb{R}^{d}\ {\rm a.s.}
\end{eqnarray*}
Especially for $t=0$,
\begin{eqnarray*}
Y_s^{0,x}=\hat{Y}_s^{0,x}\ {\rm for}\ {\rm a.e.}\ s\in[0,T],\
x\in\mathbb{R}^{d}\ {\rm a.s.}
\end{eqnarray*}
By Lemma \ref{qi045} again,
\begin{eqnarray*}
E[\int_{0}^{T}\int_{\mathbb{R}^d}|u(s,x)-\hat{u}(s,x)|^2\rho^{-1}(x)dxds]\leq
C_pE[\int_{0}^{T}\int_{\mathbb{R}^d}|Y_s^{0,x}-\hat{Y}_s^{0,x}|^2)\rho^{-1}(x)dxds]=0.
\end{eqnarray*}
So $u(s,x)=\hat{u}(s,x)$ for a.e. $s\in[0,T]$, $x\in\mathbb{R}^{d}$
a.s. The uniqueness is proved. The uniqueness implies that for any selection $u$ in the equivalence class of solution of the SPDE (\ref{bz}), $u(s,x)=Y_s^{s,x}$ for a.e. $s\in[0,T]$, $x\in\mathbb{R}^d$ a.s. Moreover, noting that $(u(s,X_s^{t,x}),\sigma^*\nabla u(s,X_s^{t,x}))$ solves the BDSDE (\ref{spdespgrowth2}) and using the uniqueness of solution of BDSDE (\ref{spdespgrowth2}) in the equivalence class, we have (\ref{dd}) for any representative $Y$ in the equivalence class of the solution of BDSDE (\ref{spdespgrowth2}).
$\hfill\diamond$

\section{BDSDEs and SPDEs with infinite dimensional noise}
\setcounter{equation}{0}

In this section, the main tasks are to prove the existence and
uniqueness of solution to BDSDE with infinite dimensional noise
(\ref{spdespgrowth2}) and give the probabilistic representation of
SPDE (\ref{ch}) with an initial value. Considering BDSDE (\ref{cm}),
the equivalent form of BDSDE (\ref{spdespgrowth2}), we assume
Conditions (H.1)--(H.6) except for (H.2) which will be replaced by the following refined condition.
\begin{description}
\item[(H.2)$^*$.] 
Assume $\|g(0,0,0)\|_{\mathcal{L}^{2}_{U_0}(\mathbb{R}^{1})}^{2}<\infty$,
and there exist constants $L,L_j\geq0$ with $\sum_{j=1}^\infty
L^{2}_j<\infty$ s.t. 
for any $s,s_1,s_2\in[0,T]$, $x,x_1,x_2\in{\mathbb{R}^{d}}$,
$y,y_1,y_2\in {\mathbb{R}^{1}}$,
\begin{eqnarray*}
&&|f(s,x_1,y)-f(s,x_2,y)|\leq L(1+|y|^{p})|x_1-x_2|,\\
&&|g_j(s_1,x_1,y_1)-g_j(s_2,x_2,y_2)|\leq
L_j(|s_1-s_2|+|x_1-x_2|+|y_1-y_2|).
\end{eqnarray*}
Moreover, we assume $\partial_yf, \partial_yg_j$ exist and satisfy
\begin{eqnarray*}
&&|\partial_yf(s,x_1,y)-\partial_yf(s,x_2,y)|\leq L(1+|y|^{p-1})|x_1-x_2|,\\
&&|\partial_yf(s,x,y_1)-\partial_yf(s,x,y_2)|\leq L(1+|y_1|^{p-2}+|y_2|^{p-2})|y_1-y_2|,\\
&&|\partial_yg_j(s,x,y)|\leq L_j,\\
&&|\partial_yg_j(s,x,y_1)-\partial_yg_j(s,x,y_2)|\leq L_j(|x_1-x_2|+|y_1-y_2|).
\end{eqnarray*}
\end{description}
\begin{rmk}
(i) An equivalent transformation of BDSDE (\ref{spdespgrowth2})
similarly as in (\ref{e}) allows us to take the monotone constant
$\mu=0$ in Condition (H.3) without losing any generality.\\
(ii) 
Similar to (\ref{au}), by Condition (H.2)$^*$ we can get
\begin{eqnarray}\label{cp}
\sup_{s\in[0,T]}\int_{\mathbb{R}^d}(\sum_{j=1}^\infty|g_j(s,x,0)|^{2})^p\rho^{-1}(x)dx<\infty.
\end{eqnarray}
\end{rmk}
\begin{thm}\label{30}
Under Conditions (H.1), (H.2)$^*$, (H.3)--(H.6), BDSDE (\ref{spdespgrowth2}) has a
unique solution $(Y_\cdot^{t,\cdot}, Z_\cdot^{t,\cdot})\in
S^{2p,0}([t,T];L_{\rho}^{2p}({\mathbb{R}^{d}};{\mathbb{R}^{1}}))\times
M^{2,0}([t,T];L_{\rho}^2({\mathbb{R}^{d}};{\mathbb{R}^{d}}))$.
\end{thm}
{\em Proof}. For every $N\in\mathbb{N}$, it follows from Theorem
\ref{21} that BDSDE (\ref{cn}) has a unique solution
$(Y_\cdot^{t,\cdot,N}, Z_\cdot^{t,\cdot,N})\in
S^{2p,0}([t,T];L_{\rho}^{2p}({\mathbb{R}^{d}};{\mathbb{R}^{1}}))\times
M^{2,0}([t,T];L_{\rho}^2({\mathbb{R}^{d}};{\mathbb{R}^{d}}))$. We
then prove that $(Y_\cdot^{t,\cdot,N}, Z_\cdot^{t,\cdot,N})$ is a
Cauchy sequence.

First note that if we do a similar estimates on
$(Y_\cdot^{t,\cdot,N}, Z_\cdot^{t,\cdot,N})$ as Lemma \ref{1}, by
Condition (H.1), (H.4) and (\ref{cp}) we can get
\begin{eqnarray}\label{cr}
&&\sup_NE[\sup_{s\in[t,T]}\int_{\mathbb{R}^d}|Y_{s}^{t,x,N}|^{2p}\rho^{-1}(x)dx]+\sup_NE[\int_t^T\int_{\mathbb{R}^d}|Z_s^{t,x,N}|^{2}\rho^{-1}(x)dxds]\nonumber\\
&\leq&C_p\int_{\mathbb{R}^d}|h(x)|^{2p}\rho^{-1}(x)dx+C_p\int_{t}^{T}\int_{\mathbb{R}^d}|f_0(s,x)|^{2p}\rho^{-1}(x)dxds\nonumber\\
&&+C_p\int_{t}^{T}\int_{\mathbb{R}^d}\sum_{j=1}^\infty|g_j(s,x,0)|^{2p}\rho^{-1}(x)dxds<\infty.
\end{eqnarray}
Then for $M,N\in\mathbb{N}$, $j\in\mathbb{N}$, define
\begin{eqnarray*}
(\bar{Y}_r^{t,x,M,N},\bar{Z}_r^{t,x,M,N})&=&(Y_r^{t,x,M}-Y_r^{t,x,N},Z_r^{t,x,M}-Z_r^{t,x,N}),\nonumber\\
\bar{g}^{M,N}_j(r,x)&=&g_j(r,X_r^{t,x},Y_r^{t,x,M})-g_j(r,X_r^{t,x},Y_r^{t,x,N}).
\end{eqnarray*}
Without losing any generality, assume $M\geq N$. Applying
It$\hat {\rm o}$'s formula to ${\rm
e}^{Kr}{|\bar{Y}_r^{t,x,M,N}|}^{2p}$ and noting the monotonicity condition of $f$ and the Lipschitz condition of $g_j$, we obtain
\begin{eqnarray}\label{cq}
&&\int_{\mathbb{R}^{d}}{\rm
e}^{Ks}{|\bar{Y}_s^{t,x,M,N}|}^{2p}\rho^{-1}(x)dx+p(2p-1)\int_s^T\int_{\mathbb{R}^{d}}{\rm
e}^{Kr}{|\bar{Y}_r^{t,x,M,N}|}^{2p-2}{|\bar{Z}_r^{t,x,M,N}|}^{2}\rho^{-1}(x)dxdr\nonumber\\
&&+(K-p(2p-1)\sum_{j=1}^{\infty}L^2_j-2\varepsilon)\int_s^T\int_{\mathbb{R}^{d}}{\rm
e}^{Kr}{|\bar{Y}_r^{t,x,M,N}|}^{2p}\rho^{-1}(x)dxdr\nonumber\\
&\leq&C_p\sum_{j=N+1}^{M}L_j^{2}\int_s^T\int_{\mathbb{R}^{d}}{\rm
e}^{Kr}|Y_r^{t,x,M}|^{2p}\rho^{-1}(x)dxdr+C_p\int_s^T\int_{\mathbb{R}^{d}}{\rm
e}^{Kr}(\sum_{j=N+1}^{M}|g_j(r,X_r^{t,x},0)|^{2})^p\rho^{-1}(x)dxdr\nonumber\\
&&-2p\sum_{j=1}^{N}\int_{s}^{T}\int_{\mathbb{R}^{d}}{\rm
e}^{Kr}{|\bar{Y}_r^{t,x,M,N}|}^{2p-2}\bar{Y}_r^{t,x,M,N}\bar{g}^{M,N}_j(r,x)\rho^{-1}(x)dxd^\dagger{\hat{\beta}}_j(r)\nonumber\\
&&-2p\sum_{j=N+1}^{M}\int_{s}^{T}\int_{\mathbb{R}^{d}}{\rm
e}^{Kr}|\bar{Y}_r^{t,x,M,N}|^{2p-2}\bar{Y}_r^{t,x,M,N}g_j(r,X_r^{t,x},Y_r^{t,x,M})\rho^{-1}(x)dxd^\dagger{\hat{\beta}}_j(r)\nonumber\\
&&-2p\int_{s}^{T}\langle\int_{\mathbb{R}^{d}}{\rm
e}^{Kr}|\bar{Y}_r^{t,x,M,N}|^{2p-2}\bar{Y}_r^{t,x,M,N}\bar{Z}_{r}^{t,x,M,N}\rho^{-1}(x)dx,dW_r\rangle.
\end{eqnarray}
Choosing sufficiently large $K$ and taking expectation on both sides
of (\ref{cq}), by Lemma \ref{qi045}, (\ref{cp}) and (\ref{cr}) we have
\begin{eqnarray}\label{cs}
&&\lim_{M,N\to\infty}E[\int_s^T\int_{\mathbb{R}^{d}}{|\bar{Y}_r^{t,x,M,N}|}^{2p}\rho^{-1}(x)dxdr]+\lim_{M,N\to\infty}E[\int_s^T\int_{\mathbb{R}^{d}}{|\bar{Y}_r^{t,x,M,N}|}^{2p-2}{|\bar{Z}_r^{t,x,M,N}|}^{2}\rho^{-1}(x)dxdr]\nonumber\\
&\leq&\lim_{M,N\to\infty}C_p\sum_{j=N+1}^{M}L_j^{2}\sup_NE[\int_s^T\int_{\mathbb{R}^{d}}|Y_r^{t,x,N}|^{2p}\rho^{-1}(x)dxdr]\nonumber\\
&&+\lim_{M,N\to\infty}C_p\int_s^T\int_{\mathbb{R}^{d}}(\sum_{j=N+1}^{M}|g_j(r,x,0)|^{2})^p\rho^{-1}(x)dxdr=0.
\end{eqnarray}
Considering (\ref{cq}) again and applying the B-D-G inequality, 
by (\ref{cs})
we have
\begin{eqnarray}\label{ct}
&&\lim_{M,N\to\infty}E[\sup_{s\in[t,T]}\int_{\mathbb{R}^{d}}{|\bar{Y}_s^{t,x,M,N}|}^{2p}\rho^{-1}(x)dx]\nonumber\\
&\leq&\lim_{M,N\to\infty}C_pE[\int_t^T\int_{\mathbb{R}^{d}}{|\bar{Y}_r^{t,x,M,N}|}^{2p}\rho^{-1}(x)dxdr]\nonumber\\
&&+\lim_{M,N\to\infty}C_p\sum_{j=N+1}^{M}L_j^{2p}\sup_NE[\int_t^T\int_{\mathbb{R}^{d}}|Y_r^{t,x,N}|^{2p}\rho^{-1}(x)dxdr]\nonumber\\
&&+\lim_{M,N\to\infty}C_p\sum_{j=N+1}^{M}\int_t^T\int_{\mathbb{R}^{d}}|g_j(r,x,0)|^{2p}\rho^{-1}(x)dxdr\nonumber\\
&&+\lim_{M,N\to\infty}C_pE[\int_t^T\int_{\mathbb{R}^{d}}{|\bar{Y}_r^{t,x,M,N}|}^{2p-2}{|\bar{Z}_r^{t,x,M,N}|}^{2}\rho^{-1}(x)dxdr]=0.
\end{eqnarray}
The estimates (\ref{cs}) and (\ref{ct}) imply that
$(Y_\cdot^{t,\cdot,N}, Z_\cdot^{t,\cdot,N})$ is a Cauchy sequence in
$S^{2p,0}([t,T];L_{\rho}^{2p}({\mathbb{R}^{d}};{\mathbb{R}^{1}}))\\ \times
M^{2,0}([t,T];L_{\rho}^2({\mathbb{R}^{d}};{\mathbb{R}^{d}}))$. So
there exists $(Y_\cdot^{t,\cdot}, Z_\cdot^{t,\cdot})\in
S^{2p,0}([t,T];L_{\rho}^{2p}({\mathbb{R}^{d}};{\mathbb{R}^{1}}))\times
M^{2,0}([t,T];L_{\rho}^2({\mathbb{R}^{d}};{\mathbb{R}^{d}}))$ as the limit of $(Y_\cdot^{t,\cdot,N}, Z_\cdot^{t,\cdot,N})$. We
then show that $(Y_\cdot^{t,\cdot}, Z_\cdot^{t,\cdot})$ is the
solution of BDSDE (\ref{spdespgrowth2}). According to Definition
\ref{qi051}, we only need to verify that for a.e.
$x\in{\mathbb{R}^{d}}$, $(Y_s^{t,x},Z_s^{t,x})$ satisfies
(\ref{spdespgrowth2}). For this, we prove that for arbitrary
$\varphi\in C_c^{0}(\mathbb{R}^d;\mathbb{R}^1)$, the integration
form of (\ref{cn}) with $\varphi$ converges to the integration form
of (\ref{spdespgrowth2}) with $\varphi$ in
$L^1(\Omega;\mathbb{R}^{1})$ along a subsequence. Due to the strong
convergence of $(Y_\cdot^{t,\cdot,N}, Z_\cdot^{t,\cdot,N})$ to
$(Y_\cdot^{t,\cdot}, Z_\cdot^{t,\cdot})$ in
$S^{2p,0}([t,T];L_{\rho}^{2p}({\mathbb{R}^{d}};{\mathbb{R}^{1}}))\times
M^{2,0}([t,T];L_{\rho}^2({\mathbb{R}^{d}};{\mathbb{R}^{d}}))$, only
the convergence of the drift term and the diffusion term w.r.t.
$\hat{B}$ are not obvious. In fact, the convergence of the diffusion term
w.r.t. $\hat{B}$ in $L^1(\Omega;\mathbb{R}^{1})$ can be referred to
the proof of Theorem 3.2 in \cite{zh-zh1}. In what follows, we show
the convergence of the drift term of the integration form of
(\ref{cn}) to the corresponding term of (\ref{spdespgrowth2}) in
$L^1(\Omega;\mathbb{R}^{1})$ along a subsequence. Since for
arbitrary $\varphi\in C_c^{0}(\mathbb{R}^d;\mathbb{R}^1)$ and
$0<\delta<1$,
\begin{eqnarray*}
&&E[\ |\int_{s}^{T}\int_{\mathbb{R}^{d}}\big(f(r,X_{r}^{t,x},Y_r^{t,x,N})-f(r,X_{r}^{t,x},Y_r^{t,x})\big)\varphi(x)dxdr|]\nonumber\\
&\leq&C_p\big(E[\int_{s}^{T}\int_{\mathbb{R}^{d}}|f(r,X_{r}^{t,x},Y_r^{t,x,N})-f(r,X_{r}^{t,x},Y_r^{t,x})|^{1+\delta}\rho^{-1}(x)dxdr]\big)^{{1\over{1+\delta}}},
\end{eqnarray*}
we only need to prove that along a subsequence
\begin{eqnarray}\label{cv}
\lim_{N\to\infty}E[\int_{s}^{T}\int_{\mathbb{R}^{d}}|f(r,X_{r}^{t,x},Y_r^{t,x,N})-f(r,X_{r}^{t,x},Y_r^{t,x})|^{1+\delta}\rho^{-1}(x)dxdr]=0.
\end{eqnarray}
First we will find a subsequence of $Y_s^{t,x,N}$, still denoted by $Y_s^{t,x,N}$, s.t.
\begin{eqnarray*}
Y_s^{t,x,N}\longrightarrow Y_s^{t,x}\ {\rm for}\ {\rm a.e.}\ s\in[t,T],\
x\in\mathbb{R}^{d}\ {\rm a.s.}\ {\rm and}\
E[\int_{t}^{T}\int_{\mathbb{R}^{d}}\sup_{N}|Y_s^{t,x,N}|^{2p}\rho^{-1}(x)dxds]<\infty.
\end{eqnarray*}
For this, by the strong convergence of $Y_\cdot^{t,\cdot,N}$ to
$Y_\cdot^{t,\cdot}$ in
$S^{2p,0}([t,T];L_{\rho}^{2p}({\mathbb{R}^{d}};{\mathbb{R}^{1}}))$, we may assume without losing any generality that
$Y_s^{t,x,N}\longrightarrow Y_s^{t,x}$ for a.e. $s\in[t,T]$,
$x\in\mathbb{R}^{d}$ a.s. and extract a subsequence of
$Y_s^{t,x,N}$, still denoted by
$Y_s^{t,x,N}$, s.t.
\begin{eqnarray*}
(E[\int_{0}^{T}\int_{\mathbb{R}^{d}}|Y_s^{t,x,N+1}-Y_s^{t,x,N}|^{2p}\rho^{-1}(x)dxds])^{1\over{2p}}\leq{1\over{2^N}}.
\end{eqnarray*}
For any $N$,
\begin{eqnarray*}
|Y_s^{t,x,N}|\leq|Y_s^{t,x,1}|+\sum_{i=1}^{N-1}|Y_s^{t,x,i+1}-Y_s^{t,x,i}|\leq|Y_s^{t,x,1}|+\sum_{i=1}^\infty|Y_s^{t,x,i+1}-Y_s^{t,x,i}|.
\end{eqnarray*}
Then by the triangle inequality of a norm, we have
\begin{eqnarray*}
&&(E[\int_{t}^{T}\int_{\mathbb{R}^{d}}\sup_N|Y_s^{t,x,N}|^{2p}\rho^{-1}(x)dxds])^{1\over{2p}}\\
&\leq&(E[\int_{t}^{T}\int_{\mathbb{R}^{d}}(|Y_s^{t,x,1}|+\sum_{i=1}^\infty|Y_s^{t,x,i+1}-Y_s^{t,x,i}|)^{2p}\rho^{-1}(x)dxds])^{1\over{2p}}\\
&\leq&(E[\int_{t}^{T}\int_{\mathbb{R}^{d}}|Y_s^{t,x,1}|^{2p}\rho^{-1}(x)dxds])^{1\over{2p}}+\sum_{i=1}^\infty(E[\int_{t}^{T}\int_{\mathbb{R}^{d}}|Y_s^{t,x,i+1}-Y_s^{t,x,i}|^{2p}\rho^{-1}(x)dxds])^{1\over{2p}}\\
&\leq&(E[\int_{t}^{T}\int_{\mathbb{R}^{d}}|Y_s^{t,x,1}|^{2p}\rho^{-1}(x)dxds])^{1\over{2p}}+\sum_{i=1}^\infty{1\over{2^i}}<\infty.
\end{eqnarray*}
Noticing
\begin{eqnarray*}
&&E[\int_{s}^{T}\int_{\mathbb{R}^{d}}\sup_{N}|f(r,X_{r}^{t,x},Y_r^{t,x,N})-f(r,X_{r}^{t,x},Y_r^{t,x})|^{1+\delta}\rho^{-1}(x)dxdr]\\
&\leq&C_pE[\int_{s}^{T}\int_{\mathbb{R}^{d}}(|f_0(r,x)|^{1+\delta}+\sup_{N}|Y_r^{t,x,N}|^{(1+\delta)p}++|Y_r^{t,x}|^{(1+\delta)p})\rho^{-1}(x)dxdr]<\infty,
\end{eqnarray*}
we have
\begin{eqnarray}\label{cu}
\sup_{N}|f(r,X_{r}^{t,x},Y_r^{t,x,N})-f(r,X_{r}^{t,x},Y_r^{t,x})|^{1+\delta}<\infty
\ \ {\rm for}\ {\rm a.e.}\ s\in[t,T],\ x\in\mathbb{R}^{d}\ {\rm
a.s.}
\end{eqnarray}
Therefore, by (\ref{cr}) and (\ref{cu})
\begin{eqnarray*}
&&\lim_{n\to\infty}\sup_NE[\int_s^T\int_{\mathbb{R}^d}|f(r,X_{r}^{t,x},Y_r^{t,x,N})-f(r,X_{r}^{t,x},Y_r^{t,x})|^{1+\delta}\nonumber\\
&&\ \ \ \ \ \ \ \ \ \ \ \ \ \ \ \ \ \ \ \ \ \ \ \ \ \ \times
I_{\{|f(r,X_{r}^{t,x},Y_r^{t,x,N})-f(r,X_{r}^{t,x},Y_r^{t,x})|^{1+\delta}>n\}}(r,x)\rho^{-1}(x)dxdr]\nonumber\\
&\leq&\lim_{n\to\infty}\sup_N\big(E[\int_s^T\int_{\mathbb{R}^d}|f(r,X_{r}^{t,x},Y_r^{t,x,N})-f(r,X_{r}^{t,x},Y_r^{t,x})|^2\rho^{-1}(x)dxdr]\big)^{{1+\delta}\over2}\nonumber\\
&&\ \ \ \ \ \ \ \ \ \ \ \ \ \ \times
\big(E[\int_s^T\int_{\mathbb{R}^d}I_{\{|f(r,X_{r}^{t,x},Y_r^{t,x,N})-f(r,X_{r}^{t,x},Y_r^{t,x})|^{1+\delta}>n\}}(r,x)\rho^{-1}(x)dxdr]\big)^{{1-\delta}\over2}\nonumber\\
&\leq&\big(\sup_NE[\int_s^T\int_{\mathbb{R}^d}(|f_0(r,x)|^2+|Y_r^{t,x,N}|^{2p}+|Y_r^{t,x}|^{2p})\rho^{-1}(x)dxdr]\big)^{{1+\delta}\over2}\nonumber\\
&&\ \ \times
\big(E[\int_s^T\int_{\mathbb{R}^d}\lim_{n\to\infty}I_{\{\sup_{N}|f(r,X_{r}^{t,x},Y_r^{t,x,N})-f(r,X_{r}^{t,x},Y_r^{t,x})|^{1+\delta}>n\}}(r,x)\rho^{-1}(x)dxdr]\big)^{{1-\delta}\over2}=0.
\end{eqnarray*}
The above convergence demonstrates that
$|f(r,X_{r}^{t,x},Y_r^{t,x,N})-f(r,X_{r}^{t,x},Y_r^{t,x})|^{1+\delta}$
is uniformly integrable, which together with the fact that
$\lim_{N\to\infty}|f(r,X_{r}^{t,x},Y_r^{t,x,N})-f(r,X_{r}^{t,x},Y_r^{t,x})|^{1+\delta}=0$
for a.e. $s\in[t,T]$, $x\in\mathbb{R}^{d}$ a.s. leads to (\ref{cv}).
So the existence of solution of BDSDE (\ref{spdespgrowth2}) follows.
As for the uniqueness proof, it is similar to the uniqueness proof
of Theorem \ref{21}. The proof of Theorem \ref{30} is completed. $\hfill\diamond$\\

With the results on BDSDE (\ref{spdespgrowth2}), the results on its
corresponding SPDE (\ref{bz}) follow.
\begin{thm}\label{20}
Define $u(t,x)=Y_t^{t,x}$, where $(Y_s^{t,x},Z_s^{t,x})$ is the
solution of BDSDE (\ref{spdespgrowth2}) under Conditions
(H.1), (H.2)$^*$, (H.3)--(H.6), then $u(t,x)$ is the unique weak solution of SPDE
(\ref{bz}). Moreover,
\begin{eqnarray*}
u(s,X_s^{t,x})=Y_s^{t,x},\ (\sigma^*\nabla u)(s,X^{t,x}_s)=Z_s^{t,x}
\ {\rm for}\ {\rm a.e.}\ s\in[t,T],\ x\in\mathbb{R}^{d}\ {\rm a.s.}
\end{eqnarray*}
\end{thm}
{\em Proof}. Consider BDSDE (\ref{cn}) and the following SPDE with
finite dimensional noise:
\begin{eqnarray}\label{cx}
u_N(t,x)=h(x)+\int_{t}^{T}[\mathscr{L}u_N(s,x)+f\big(s,x,u_N(s,x)\big)]ds-\sum_{j=1}^{N}\int_{t}^{T}g_j\big(s,x,u_N(s,x)\big)d^\dagger{\hat{\beta}}_j(s).
\end{eqnarray}
By Theorem \ref{31} we know that $Y_t^{t,x,N}$ is the weak solution
of SPDE (\ref{cx}) and
\begin{eqnarray*}
u_N(s,X_s^{t,x})=Y_s^{t,x,N},\ (\sigma^*\nabla
u_N)(s,X^{t,x}_s)=Z_s^{t,x,N} \ {\rm for}\ {\rm a.e.}\ s\in[t,T],\
x\in\mathbb{R}^{d}\ {\rm a.s.}
\end{eqnarray*}
The
remaining part of the proof is to verify that $u_N(s,x)$ is a Cauchy sequence
in $\mathcal{H}$ and its limit $u(s,x)$ is the weak solution of
SPDE (\ref{bz}). The procedure of these proofs are actually similar
to Proposition 4.2 and Theorem 4.3 in \cite{zh-zh1} where a
Lipschitz condition to $f(s,x,y)$ on $y$ rather than polynomial growth
condition is assumed. However, the polynomial growth condition in the
arguments brings the trouble only when verifying that for arbitrary
$\varphi\in C_c^{\infty}(\mathbb{R}^d;\mathbb{R}^1)$ the integration
form of the drift term of (\ref{cx}) converges
to the corresponding term of (\ref{bz}) in $L^1(\Omega)$, i.e.
\begin{eqnarray}\label{cy}
\lim_{N\to\infty}E[\
|\int_{0}^{T}\int_{\mathbb{R}^{d}}\big(f(s,x,u_N(s,x))-f(s,x,u(s,x))\big)\varphi(x)dxdr|]=0.
\end{eqnarray}
If we know that for arbitrary $\varphi\in
C_c^{\infty}(\mathbb{R}^d;\mathbb{R}^1)$ and $\delta<1$,
\begin{eqnarray}\label{cz}
\lim_{N\to\infty}E[\int_{0}^{T}\int_{\mathbb{R}^{d}}|f(s,x,u_N(s,x))-f(s,x,u(s,x))|^{1+\delta}\rho^{-1}(x)dxdr]=0,
\end{eqnarray}
then (\ref{cy}) follows from H${\rm \ddot{o}}$lder inequality
immediately. In fact, noting (\ref{cv}) we can prove (\ref{cz}) by the strong
convergence of $u_N(s,x)$ to $u(s,x)$ in $\mathcal{H}$,
$u(s,X_s^{t,x})=Y_s^{t,x}$ for a.e. $s\in[t,T]$,
$x\in\mathbb{R}^{d}$ a.s. and the equivalence of norm principle. $\hfill\diamond$

\begin{rmk}\label{6}
Consider 
a simpler form of SPDE (\ref{bz})
with coefficients $f,g$ being independent of time variables. If we choose Brownian
motion $\hat{B}$ in backward SPDE 
as the time reversal version of
Brownian motion $B$ in SPDE (\ref{ch}), i.e.
$\hat{B}_s=B_{T-s}-B_T$, $0\leq
s\leq T$, and let $u$ be the weak solution of corresponding backward SPDE, 
then we can see easily that $v(t)\triangleq u(T-t)$ is the unique
weak solution of SPDE (\ref{ch}) s.t. $(v,\sigma^*\nabla
v)\in
L^{2p}([0,T];L_{\rho}^{2p}({\mathbb{R}^{d}};{\mathbb{R}^{1}}))\times
L^{2}([0,T];L_{\rho}^2({\mathbb{R}^{d}};{\mathbb{R}^{d}}))$.
\end{rmk}

\section{Infinite horizon BDSDEs and stationary solutions of SPDEs}
\setcounter{equation}{0}

In this section, we first consider the infinite horizon BDSDE with
polynomial growth coefficients. For this, we assume
the previous conditions (H.1), (H.2)$^*$, (H.3) with the following
changes:
\begin{description}
\item[(H.7).] Change ``$s\in[0,T]$" to
``$s\in[0,\infty)$" in (H.1).
\item[(H.8).] Change ``$s,s_1,s_2\in[0,T]$" to ``$s,s_1,s_2\in[0,\infty)$"
in (H.2)$^*$.
\item[(H.9).] Change ``$\mu\in\mathbb{R}^{1}$" to ``$\mu>0$ with $2\mu-{K}-{p(2p-1)}\sum_{j=1}^{\infty}L_j>0$", ``$s\in[0,T]$" to ``$s\in[0,\infty)$"
and ``$\leq\mu|y_1-y_2|^2$" to ``$\leq-\mu|y_1-y_2|^2$" in (H.3).
\end{description}
\begin{rmk}
From Conditions (H.7) and (H.8), for any given $K>0$, we can deduce
\begin{eqnarray*}\label{da}
&&\int_0^\infty\int_{\mathbb{R}^d}{\rm e}^{-Ks}(|f(s,x,0)|^{2p}+(\sum_{j=1}^\infty|g_j(s,x,0)|^{2})^p)\rho^{-1}(x)dxds\nonumber\\
&\leq&C_p\int_0^\infty\int_{\mathbb{R}^d}{\rm
e}^{-Ks}(|{f}_0(s,x)|^{2p}+(\sum_{j=1}^\infty
L^{2}_j)^ps^{2p}+(\sum_{j=1}^\infty
L^{2}_j)^p|x|^{2p}+(\sum_{j=1}^\infty|g_j(0,0,0)|^{2})^p)\rho^{-1}(x)dxds<\infty.
\end{eqnarray*}
\end{rmk}

Then we have the existence and uniqueness theorem for BDSDE
(\ref{qi30}):
\begin{thm}\label{zz47} Under Conditions {\rm(H.5)}--{\rm(H.9)},
BDSDE (\ref{qi30}) has a unique solution
$(Y_{\cdot}^{t,\cdot},Z_{\cdot}^{t,\cdot})\in S^{2p,-K}\bigcap
M^{2p,-K}([t,\infty];L_{\rho}^{2p}({\mathbb{R}^{d}};{\mathbb{R}^{1}}))\times
M^{2,-K}([t,\infty];L_{\rho}^2({\mathbb{R}^{d}};{\mathbb{R}^{d}}))$.
\end{thm}
{\em Proof}.
Here we only prove the existence of solution as the uniqueness is
similar to the proof of uniqueness in Theorem \ref{21}. For the same
reason of Remark \ref{16}, for a.e. $x$, $(Y_s^{t,x},Z_s^{t,x})$
satisfies (\ref{qi30}) is equivalent to that for arbitrary
$\varphi\in C_c^{0}(\mathbb{R}^d;\mathbb{R}^1)$,
$({Y}_s^{t,x},{Z}_s^{t,x})$ satisfies 
\begin{eqnarray}\label{cg}
\int_{\mathbb{R}^{d}}{\rm e}^{-Ks}Y_s^{t,x}\varphi(x)dx&=&\int_{s}^{\infty}\int_{\mathbb{R}^{d}}{\rm e}^{-Kr}f(r,X_{r}^{t,x},Y_{r}^{t,x})\varphi(x)dxdr+\int_{s}^{\infty}\int_{\mathbb{R}^{d}}K{\rm e}^{-Kr}Y_{r}^{t,x}\varphi(x)dxdr\nonumber\\
&&-\int_{s}^{\infty}\int_{\mathbb{R}^{d}}{\rm e}^{-Kr}g(r,X_{r}^{t,x},Y_{r}^{t,x})\varphi(x)dxd^\dagger{\hat{B}}_r\nonumber\\
&&-\int_{s}^{\infty}\langle \int_{\mathbb{R}^{d}}{\rm
e}^{-Kr}Z_{r}^{t,x}\varphi(x)dx,dW_r\rangle\ \ \ P-{\rm a.s.}
\end{eqnarray}
For each $n\in\mathbb{N}$, we define a
sequence of BDSDEs by setting $h=0$ and $T=n$ in BDSDE (\ref{spdespgrowth2}): 
\begin{eqnarray}\label{zz62}
Y_{s}^{t,x,n}&=&\int_{s}^{n}f(r,X_{r}^{t,x},Y_{r}^{t,x,n})dr-\int_{s}^{n}g(r,X_{r}^{t,x},Y_{r}^{t,x,n})d^\dagger{\hat{B}}_r-\int_{s}^{n}\langle
Z_{r}^{t,x,n},dW_r\rangle.\ \ \ \ \ \ \ \ 
\end{eqnarray}
It is easy to verify that BDSDE (\ref{zz62}) satisfies conditions of
Theorem \ref{30}. Hence, for each $n$, there exists
$({Y}_\cdot^{t,\cdot,n},{Z}_\cdot^{t,\cdot,n})\in
S^{2p,-K}([t,n];L_{\rho}^{2p}({\mathbb{R}^{d}};{\mathbb{R}^{1}}))\times
M^{2,-K}([t,n];L_{\rho}^2({\mathbb{R}^{d}};{\mathbb{R}^{d}}))$ and
$({Y}_s^{t,x,n},{Z}_s^{t,x,n})$ is the unique solution of BDSDE
(\ref{zz62}). Therefore, for arbitrary $\varphi\in
C_c^{0}(\mathbb{R}^d;\mathbb{R}^1)$,
$({Y}_s^{t,x,n},{Z}_s^{t,x,n})$ satisfies 
\begin{eqnarray}\label{cf}
\int_{\mathbb{R}^{d}}{\rm e}^{-Ks}Y_s^{t,x,n}\varphi(x)dx&=&\int_{s}^{n}\int_{\mathbb{R}^{d}}{\rm e}^{-Kr}f(r,X_{r}^{t,x},Y_{r}^{t,x,n})\varphi(x)dxdr+\int_{s}^{n}\int_{\mathbb{R}^{d}}K{\rm e}^{-Kr}Y_{r}^{t,x,n}\varphi(x)dxdr\nonumber\\
&&-\int_{s}^{n}\int_{\mathbb{R}^{d}}{\rm e}^{-Kr}g(r,X_{r}^{t,x},Y_{r}^{t,x,n})\varphi(x)dxd^\dagger{\hat{B}}_r\nonumber\\
&&-\int_{s}^{n}\langle \int_{\mathbb{R}^{d}}{\rm
e}^{-Kr}Z_{r}^{t,x,n}\varphi(x)dx,dW_r\rangle\ \ \ P-{\rm a.s.}
\end{eqnarray}
Let $(Y_{s}^{n}, Z_{s}^{n})_{s>n}=(0, 0)$. Then
$({Y}_\cdot^{t,\cdot,n},{Z}_\cdot^{t,\cdot,n})\in S^{2p,-K}\bigcap
M^{2p,-K}([t,\infty);L_{\rho}^{2p}({\mathbb{R}^{d}};{\mathbb{R}^{1}}))\times
M^{2,-K}([t,\infty);\\
L_{\rho}^2({\mathbb{R}^{d}};{\mathbb{R}^{d}}))$.
We use a similar argument as in the proof of Theorem 5.1 in
\cite{zh-zh1} to prove that $({Y}_{\cdot}^{t,{\cdot},n},{Z}_{\cdot}^{t,{\cdot},n})$ is a Cauchy sequence in $S^{2p,-K}\bigcap
M^{2p,-K}([t,\infty);L_{\rho}^{2p}({\mathbb{R}^{d}};{\mathbb{R}^{1}}))\times
M^{2,-K}([t,\infty);L_{\rho}^2({\mathbb{R}^{d}};{\mathbb{R}^{d}}))$. Assume without loss of any generality $m\geq n$. On the interval $[n,m]$, by Conditions (H.5), (H.7)--(H.9) we can prove that
\begin{eqnarray*}
&&\big(2p\mu-{K}-{p(2p-1)}\sum_{j=1}^{\infty}L_j-\varepsilon\big)E[\int_{n}^{m}\int_{\mathbb{R}^{d}}{\rm e}^{{-{K}}r}{|{Y}_r^{t,x,m}-Y_r^{t,x,n}|}^{2p}\rho^{-1}(x)dxdr]\nonumber\\
&&+{p(2p-1)}E[\int_{n}^{m}\int_{\mathbb{R}^{d}}{\rm e}^{{-{K}}r}{{|{Y}_r^{t,x,m}-Y_r^{t,x,n}|}^{2p-2}}|{Z}_r^{t,x,m}-Z_r^{t,x,n}|^2\rho^{-1}(x)dxdr]\nonumber\\
&\leq&C_p\int_n^m\int_{\mathbb{R}^d}{\rm e}^{-Kr}(|f(r,x,0)|^{2p}+(\sum_{j=1}^\infty|g_j(r,x,0)|^{2})^p)\rho^{-1}(x)dxdr\longrightarrow0\ \ {\rm as}\ m,n\to\infty.
\end{eqnarray*}
By above estimates and B-D-G inequality, we further have
\begin{eqnarray}\label{zz63}
E[\sup_{n\leq s\leq m}\int_{\mathbb{R}^{d}}{\rm e}^{-{K}s}|{Y}_s^{t,x,m}-Y_s^{t,x,n}|^{2p}\rho^{-1}(x)dx]\longrightarrow0\ \ {\rm as}\ m,n\to\infty.
\end{eqnarray}
On the interval $[t,n]$, a similar calculation, together with (\ref{zz63}), leads to
\begin{eqnarray*}
&&\big(2p\mu-{K}-{p(2p-1)}\sum_{j=1}^{\infty}L_j\big)E[\int_{0}^{n}\int_{\mathbb{R}^{d}}{\rm e}^{{-{K}}r}{|{Y}_r^{t,x,m}-Y_r^{t,x,n}|}^{2p}\rho^{-1}(x)dxdr]\nonumber\\
&&+{p(2p-1)}E[\int_{0}^{n}\int_{\mathbb{R}^{d}}{\rm e}^{{-{K}}r}{{|{Y}_r^{t,x,m}-Y_r^{t,x,n}|}^{2p-2}}|{Z}_r^{t,x,m}-Z_r^{t,x,n}|^2\rho^{-1}(x)dxdr]\nonumber\\
&\leq&E[\int_{\mathbb{R}^d}{\rm e}^{-Kn}|Y_n^{t,x,m}|^{2p}\rho^{-1}(x)dx]\longrightarrow0\ \ {\rm as}\ m,n\to\infty,
\end{eqnarray*}
and
\begin{eqnarray*}
E[\sup_{t\leq s\leq n}\int_{\mathbb{R}^{d}}{\rm e}^{-{K}s}|{Y}_s^{t,x,m}-Y_s^{t,x,n}|^{2p}\rho^{-1}(x)dx]\longrightarrow0\ \ {\rm as}\ m,n\to\infty.
\end{eqnarray*}
Taking account of calculations on both $[t,n]$ and $[n,m]$ we know that $({Y}_s^{t,x,n},{Z}_s^{t,x,n})$ is a Cauchy sequence in $S^{2p,-K}\bigcap
M^{2p,-K}([t,\infty);L_{\rho}^{2p}({\mathbb{R}^{d}};{\mathbb{R}^{1}}))\times
M^{2,-K}([t,\infty);L_{\rho}^2({\mathbb{R}^{d}};{\mathbb{R}^{d}}))$. Let
$({Y}_s^{t,x},{Z}_s^{t,x})$ be the limit of
$({Y}_s^{t,x,n},{Z}_s^{t,x,n})$, then we show that $({Y}_s^{t,x},{Z}_s^{t,x})$ is a solution of BDSDE
(\ref{qi30}). We only need to verify that
$({Y}_s^{t,x},{Z}_s^{t,x})$ satisfies (\ref{cg}). For this, we prove
that
along a subsequence (\ref{cf}) converges to (\ref{cg}) in
$L^1(\Omega;\mathbb{R}^1)$ term by term as $n\to\infty$. Here we
only check the drift term which is of polynomial growth, i.e. we show that
along a subsequence, as $n\to\infty$,
\begin{eqnarray*}
E[\ |\int_{s}^{n}\int_{\mathbb{R}^{d}}{\rm
e}^{-Kr}f(r,X_{r}^{t,x},Y_r^{t,x,n})\varphi(x)dxdr-\int_{s}^{\infty}\int_{\mathbb{R}^{d}}{\rm
e}^{-Kr}f(r,X_{r}^{t,x},Y_r^{t,x})\varphi(x)dxdr|]\longrightarrow0.
\end{eqnarray*}
For this, note that for arbitrary $0<\delta<1$,
\begin{eqnarray}\label{ck}
&&E[\ |\int_{s}^{n}\int_{\mathbb{R}^{d}}{\rm e}^{-Kr}f(r,X_{r}^{t,x},Y_r^{t,x,n})\varphi(x)dxdr-\int_{s}^{\infty}\int_{\mathbb{R}^{d}}{\rm e}^{-Kr}f(r,X_{r}^{t,x},Y_r^{t,x})\varphi(x)dxdr|]\nonumber\\
&\leq&
C_pE[\int_{s}^{\infty}\int_{\mathbb{R}^{d}}{\rm
e}^{-Kr}|f(r,X_{r}^{t,x},Y_r^{t,x,n})-f(r,X_{r}^{t,x},Y_r^{t,x})|^{1+\delta}\rho^{-1}(x)dxdr]\nonumber\\
&&+C_pE[\int_{n}^{\infty}\int_{\mathbb{R}^{d}}{\rm
e}^{-Kr}(|f_0(r,X_r^{t,x})|^2+|Y_r^{t,x}|^{2p})\rho^{-1}(x)dxdr].
\end{eqnarray}
Both terms on the right hand side of the above inequality converge
to $0$ along a subsequence as $n\to\infty$. The 
convergence of the first term in (\ref{ck}) is not obvious, but
can be deduced similarly as the proof of (\ref{cv}).
After verifying other terms of (\ref{cf}) converges to the
corresponding terms of (\ref{cg}) in $L^1(\Omega;\mathbb{R}^1)$ as
$n\to\infty$, we can see that $(Y_s^{t,x}, Z_s^{t,x})$ satisfies
(\ref{cg}) and the proof of Theorem \ref{zz47} is completed.
$\hfill\diamond$

\begin{rmk}\label{zhang315}
The uniqueness of solution of BDSDE (\ref{qi30}) implies if
$(\hat{Y},\hat{Z})$ is another solution, then
$Y_s^{t,\cdot}=\hat{Y}_s^{t,\cdot}$ for all $s\geq t$ a.s. and
$Z_s^{t,\cdot}=\hat{Z}_s^{t,\cdot}$ for a.e. $s\geq t$ a.s. But we
can modify the $Z$ at the measure zero exceptional set of $s$ s.t.
$Z_s^{t,\cdot}=\hat{Z}_s^{t,\cdot}$ for all $s\geq t$ a.s.
\end{rmk}

Consider the case when $f$ and $g$ are time-independent coefficients.
\begin{cor}\label{15} Under Conditions (H.5)--(H.9), BDSDE (\ref{bp}) has
a unique solution $(Y_{\cdot}^{t,\cdot},Z_{\cdot}^{t,\cdot})$ in the
space $S^{2p,-K}\bigcap
M^{2p,-K}([t,\infty];L_{\rho}^{2p}({\mathbb{R}^{d}};{\mathbb{R}^{1}}))\times
M^{2,-K}([t,\infty];L_{\rho}^2({\mathbb{R}^{d}};{\mathbb{R}^{d}}))$.
\end{cor}



We construct a measurable metric dynamical system through defining a
measurable and probability preserving shift operator. Let
$\hat{\theta}_t=\theta'_t\circ\theta''_t$, $t\geq0$, where
$\theta'_t,\theta''_t:\Omega\longrightarrow\Omega$ are measurable
mappings on $(\Omega, {\mathscr{F}}, P)$ defined by
\begin{eqnarray*}
\theta'_t\left(\begin{array}{c}\hat{B}\\
W\\
\end{array}\right)(s)=\left(\begin{array}{c}\hat{B}_{s+t}-\hat{B}_{t}\\
W_s\\
\end{array}\right),\ \ \ \ \theta''_t\left(\begin{array}{c}\hat{B}\\
W\\
\end{array}\right)(s)=\left(\begin{array}{c}\hat{B}_{s}\\
W_{s+t}-W_t\\
\end{array}\right).
\end{eqnarray*}
Then for any $s$, $t\geq0$,
(i). $P=\hat{\theta}_{t}P$;
(ii). $\hat{\theta}_{0}=I$, where $I$ is the identity transformation on $\Omega$;
(iii). $\hat{\theta}_{s}\circ\hat{\theta}_{t}=\hat{\theta}_{s+t}$.
Also for an arbitrary $\mathscr{F}$ measurable $\phi$ and $t\geq0$, set
\begin{eqnarray*}
\hat{\theta}_t\circ\phi(\omega)=\phi\big(\hat{\theta}_t(\omega)\big).
\end{eqnarray*}
For any $r\geq0$, applying $\hat{\theta}_r$ to SDE (\ref{a}), by the
uniqueness of the solution and a perfection procedure (cf. Arnold
\cite{ar}) we have
\begin{eqnarray*}\label{qi18}
\hat{\theta}_r\circ X_{s}^{t,\cdot}=\theta''_r\circ
X_{s}^{t,\cdot}=X_{s+r}^{t+r,\cdot}\ \ {\rm for}\ {\rm all}\
r,s,t\geq0\ \ {\rm a.s.}
\end{eqnarray*}

Firstly, we consider the stationarity of BDSDE (\ref{bp}), which is
equivalent to
\begin{numcases}{}\label{bq}
Y_{s}^{t,x}=Y_{T}^{t,x}+\int_{s}^{T}f(X_{r}^{t,x},Y_{r}^{t,x})dr
-\int_{s}^{T}g(X_{r}^{t,x},Y_{r}^{t,x})d^\dagger{\hat{B}}_r-\int_{s}^{T}\langle Z_{r}^{t,x},dW_r\rangle\nonumber\\
\lim_{T\rightarrow\infty}{\rm e}^{-KT}Y_{T}^{t,x}=0\ \ \ {\rm a.s.}\nonumber
\end{numcases}

\begin{thm}\label{lz7}
Under Conditions (H.5)--(H.9), let
$(Y_{s}^{t,x},Z_{s}^{t,x})$ be the solution of BDSDE (\ref{bp}),
then $(Y_{s}^{t,x},Z_{s}^{t,x})$ satisfies for any $t\geq0$,
\begin{eqnarray*}
\hat{\theta}_r\circ Y_{s}^{t,\cdot}=Y_{s+r}^{t+r,\cdot}, \ \ \
\hat{\theta}_r\circ Z_{s}^{t,\cdot}=Z_{s+r}^{t+r,\cdot}\ \ \ {\rm
for}\ {\rm all}\ r\geq0,\ s\geq t\ {\rm a.s.}
\end{eqnarray*}
In particular, for any $t\geq0$,
\begin{eqnarray}\label{br}
\theta'_r\circ Y^{t,\cdot}_t=Y^{t+r,\cdot}_{t+r}\ \ \ {\rm for}\
{\rm all}\ r\geq0\ {\rm a.s.}
\end{eqnarray}
\end{thm}

The proof is similar to the proof of Theorem 1.7 in \cite{zh-zh2}. So it is omitted in this paper.

If we regard $Y_t^{t,\cdot}$ as a function of $t$, (\ref{br}) gives
a ``crude" stationary property of $Y_t^{t,\cdot}$. To make the
``crude" stationary property ``perfect", we need to prove the a.s.
continuity of $t\longrightarrow Y_t^{t,\cdot}$ to obtain one
indistinguishable version of $Y_t^{t,\cdot}$ with ``perfect"
stationary property w.r.t. $\hat{\theta}$. As the a.s. continuity
can be similarly proved if one follows the arguments of Theorem 2.11
in \cite{zh-zh1}, here we leave out the proof. Hence it comes
without a surprise that
\begin{thm}\label{lz9} Under Conditions (H.5)--(H.9), let
$(Y_{s}^{t,x},Z_{s}^{t,x})$ be the solution of BDSDE (\ref{bp}).
Then $Y_{t}^{t,\cdot}$ satisfies the ``perfect" stationary property
w.r.t. $\theta'$, i.e.
\begin{eqnarray}\label{zz18}
\theta'_r\circ Y^{t,\cdot}_t=Y^{t+r,\cdot}_{t+r}\ \ \ {\rm for}\
{\rm all}\ t,r\geq0\ {\rm a.s.}
\end{eqnarray}
\end{thm}

We can further prove an estimate following the
proof of Lemma \ref{1}.
\begin{lem}\label{8}
Under Conditions (H.5)--(H.9), if
$(Y_\cdot^{t,\cdot},Z_\cdot^{t,\cdot})$ is the solution of BDSDE
(\ref{bp}), then we have
\begin{eqnarray*}
&&E[\sup_{s\in[t,\infty)}\int_{\mathbb{R}^d}{\rm e}^{-Ks}|Y_{s}^{t,x}|^{{8p}}\rho^{-1}(x)dx]+E[\int_t^\infty\int_{\mathbb{R}^d}{\rm e}^{-Ks}|Y_s^{t,x}|^{{8p}}\rho^{-1}(x)dxds]\\
&&+\sup_nE[\int_{t}^{\infty}\int_{\mathbb{R}^{d}}{\rm e}^{-Kr}{{|{Y}_s^{t,x}|}^{{8p}-2}}|{Z}_{s}^{t,x}|^2\rho^{-1}(x)dxds]+\sup_nE[\big(\int_t^\infty\int_{\mathbb{R}^d}{\rm
e}^{-{Kr\over{4p}}}|Z_s^{t,x}|^{2}\rho^{-1}(x)dxds\big)^{4p}]<\infty.
\end{eqnarray*}
\end{lem}

Consider BDSDE (\ref{bp}) and its solution $Y_s^{t,\cdot}$ on
$[t,T]$. We choose $\hat{B}$ as the time reversal of $B$ from time
$T'$, i.e. $\hat{B}_s=B_{T'-s}-B_{T'}$ for $s\geq0$. Note that
the random variable $Y_{T'}^{T',\cdot}$ is
$\mathscr{F}^{\hat{B}}_{T',\infty}$ measurable which is independent
of $\mathscr{F}_t^W$. Changing variable in SPDE (\ref{ch}), we
can deduce from the Correspondence Theorem
\ref{20} and Remark \ref{6} that
$v(t,\cdot)=u(T'-t,\cdot)=Y_{T'-t}^{T'-t,\cdot}$ is a weak solution of
SPDE (\ref{ch}) on $[0,T']$ if $Y_{T'}^{T',x}$ satisfies Condition (H.4). Note $Y_T^{t,x}=Y_T^{T,X_T^{t,x}}$, so Condition (H.4) reads as
\begin{description}
\item[(H.4)$^*$.] $E[\int_{\mathbb{R}^d}|Y_{T'}^{T',x}|^{{8p}}\rho^{-1}(x)dx]<L(T')$ and $E[\int_{\mathbb{R}^d}|Y_{T'}^{t',x}-Y_{T'}^{t,x}|^q\rho^{-1}(x)dx]\leq
L(T')|t'-t|^{q\over2}$ for $2\leq q\leq{8p}$ and $X$ defined in (\ref{a}), where $L(T')$ is a constant which can depend on $T'$.
\end{description}
\begin{lem}\label{29}
Let $(Y_s^{t,x},Z_s^{t,x})$ be the solution of BDSDE (\ref{bp}). 
Then for arbitrary $T'$, $Y_{T'}^{T',x}$ satisfies Condition (H.4)$^*$.
\end{lem}
{\em Proof}. It follows immediately from Lemmas \ref{qi045} and
\ref{8} that
$E[\int_{\mathbb{R}^d}|Y_{T'}^{T',x}|^{{8p}}\rho^{-1}(x)dx]\leq L(T')$.
The proof of
$E[\int_{\mathbb{R}^d}|Y_{T'}^{t',x}-Y_{T'}^{t,x}|^q\rho^{-1}(x)dx]\leq
L(T')|t'-t|^{q\over2}$ is similar to Lemma 6.2 in \cite{zh-zh1}. 
$\hfill\diamond$\\

We can further prove a claim that $v(t,\cdot)=Y_{T'-t}^{T'-t,\cdot}$
does not depend on the choice of $T'$ as done in \cite{zh-zh1}, \cite{zh-zh2}. On the probability space $(\Omega,\mathscr{F},P)$, we
define ${\theta}_t:\Omega\longrightarrow\Omega$,
$t\in\mathbb{R}^{1}$, as the shift operator of Brownian motion $B$:
\begin{eqnarray*}
{\theta}_{t}\circ{B}_s={B}_{s+t}-{B}_{t},
\end{eqnarray*}
then $\theta$ satisfies the usual conditions:
(i). $P=P\circ{\theta}_{t}$;
(ii). ${\theta}_{0}=I$;
(iii). ${\theta}_{s}\circ{\theta}_{t}={\theta}_{s+t}$.
Noticing that $\hat{B}$ is chosen as the time reversal of $B$ and
$B,W$ are independent, we can define $\hat{\theta}$, served as the
shift operator of $\hat{B}$ and $W$, to be
$\hat{\theta}_{t}\triangleq(\theta_t)^{-1}\circ\theta''_t$,
$t\geq0$. Actually ${B}$ is a two-sided Brownian motion, so
$({\theta}_{t})^{-1}={\theta}_{-t}$ is well defined (see \cite{ar})
and it is easy to see that
$\theta'_{t}\triangleq({\theta}_{t})^{-1}$, $t\in\mathbb{R}^{1}$, is
a shift operator of $\hat{B}$.

Since $v(t,\cdot)=u(T'-t,\cdot)=Y_{T'-t}^{T'-t,\cdot}$ a.s., so by
(\ref{zz18}),
\begin{eqnarray*}
{\theta}_rv(t,\cdot,\omega)=\theta'_{-r}u(T'-t,\cdot,\hat{\omega})=\theta'_{-r}\theta'_{r}u(T'-t-r,\cdot,\hat{\omega})=u(T'-t-r,\cdot,\hat{\omega})=v(t+r,\cdot,\omega),
\end{eqnarray*}
for all $r\geq0$ and $T'\geq t+r$ a.s. In particular, let
$Y(\cdot,\omega)=v(0,\cdot,\omega)=Y_{T'}^{T',\cdot}(\hat{\omega})$, then
the above formula implies: 
\begin{eqnarray*}
{\theta}_tY(\cdot,\omega)=Y(\cdot,{\theta}_t\omega)=v(t,\cdot,\omega)=v(t,\cdot,\omega,v(0,\cdot,\omega))=v(t,\cdot,\omega,Y(\cdot,\omega))\
{\rm for}\ {\rm all}\ t\geq0\ {\rm a.s.}
\end{eqnarray*}
That is to say
$v(t,\cdot,\omega)=Y(\cdot,{\theta}_t\omega)=Y_{T'-t}^{T'-t,\cdot}(\hat{\omega})$
is a stationary solution of SPDE (\ref{ch}) w.r.t. $\theta$.
Therefore we obtain
\begin{thm}\label{qi047} Under Conditions (H.5)--(H.9), for arbitrary $T'$ and $t\in[0,T']$,
let $v(t,\cdot)\triangleq Y_{T'-t}^{T'-t,\cdot}$, where
$(Y_{\cdot}^{t,\cdot},Z_{\cdot}^{t,\cdot})$ is the solution of BDSDE
(\ref{bp}) with $\hat{B}_s={B}_{T'-s}-{B}_{T'}$ for all $s\geq0$. Then
$v(t,\cdot)$ is a "perfect" stationary solution of SPDE (\ref{ch}) independent of the choice of $T'$.
\end{thm}

It is not difficult to see that in the proof of Theorem \ref{zz47}, there is no need to take $h=0$. In fact we can consider BDSDE (\ref{spdespgrowth2}) with an arbitrary $h$ satisfying Condition (H.4). Its solution is denoted by $Y_s^{t,x,T}(h)$. Then under the same conditions as in Theorem \ref{qi047}, following the same procedure of this section, we can prove without real difficulty that $Y_\cdot^{t,\cdot,T}(h)\longrightarrow Y_\cdot^{t,\cdot}$ as $T\to\infty$ in $S^{2p,-K}\bigcap
M^{2p,-K}([t,\infty];L_{\rho}^{2p}({\mathbb{R}^{d}};{\mathbb{R}^{1}}))$ and $Y_\cdot^{t,\cdot}$ is the solution of infinite horizon BDSDE (\ref{bp}). This implies that $Y_t^{t,\cdot,T}(h)\longrightarrow Y_t^{t,\cdot}$ as $T\to\infty$ in $L^{2p}(\Omega;L_{\rho}^{2p}({\mathbb{R}^{d}};{\mathbb{R}^{1}}))$. In particular, from previous result of this section, we have
$$Y_{T'-t}^{T'-t,\cdot,T}(h)\longrightarrow Y_{T'-t}^{T'-t,\cdot}=v(t,\cdot)=Y(\theta_t\cdot)\ {\rm and}\ Y_{T'}^{T',\cdot,T}(h)\longrightarrow Y_{T'}^{T',\cdot}=Y(\cdot).$$
Noticing the correspondence of the solution of BDSDE (\ref{spdespgrowth2}) and the solution of the SPDE (\ref{bz}) in Theorem \ref{20}, we have $Y_{T'}^{T',\cdot,T}(h,\hat{\omega})=u(T',h,\hat{\omega})$. Note the correspondence of the forward SPDE (\ref{ch}) and the backward SPDE (\ref{bz}). Denote $\hat{\omega}^T$ the time reversal Brownian motion of $B$ at time $T$. Then
\begin{eqnarray*}
v(T-T',h,\theta_{-(T-T')}\omega)=u(T-(T-T'),h,\widehat{\theta_{-(T-T')}\omega}^T)=u(T',h,\hat{\omega})=Y^{T',\cdot,T}_{T'}(h)
\end{eqnarray*}
since
\begin{eqnarray*}
\widehat{(\theta_{-(T-T')}\omega^T)}(s)&=&(\theta_{-(T-T')}\omega)(T-s)-(\theta_{-(T-T')}\omega)(T)\nonumber\\&=&\omega(-(T-T')+T-s)-\omega(-(T-T')+T)
=\omega(T'-s)-\omega(T')=\hat{\omega}(s).
\end{eqnarray*}
Therefore as $T\to\infty$, $v(T-T',h,\theta_{-(T-T')}\cdot)=Y_{T'}^{T',\cdot,T}(h)\longrightarrow Y(\cdot)$ in $L^{2p}(\Omega;L_{\rho}^{2p}({\mathbb{R}^{d}};{\mathbb{R}^{1}}))$. The result does not depend on the choice of $T'$. So we have proved
\begin{thm}\label{32} Assume all conditions in Theorem \ref{qi047} and $h$ satisfies Condition (H.4). Then as $T\to\infty$, $v(T,h,\theta_{-T}\omega)\longrightarrow Y(\cdot)$ in $L^{2p}(\Omega;L_{\rho}^{2p}({\mathbb{R}^{d}};{\mathbb{R}^{1}}))$, and $Y(\theta_t\omega)$ is the stationary solution of the SPDE (\ref{ch}).
\end{thm}

\begin{rmk} The result in Theorem \ref{32} is also valid under the conditions in \cite{zh-zh1} and \cite{zh-zh2} respectively.
\end{rmk}

{\bf Acknowledgements}. QZ is partially supported by the National Basic Research Program of China (973 Program) with grant
no. 2007CB814904, the National Natural Science Foundation of China with grant nos. 10771122, 11101090 and the Specialized Research Fund for
the Doctoral Program of Higher Education of China with grant
no. 20090071120002. HZ would like to thank Prof. Shanjian Tang for
kind invitation to visit Laboratory of Mathematics for Nonlinear
Science, Fudan University, and thank Prof. Shige Peng for kind invitation to visit Shandong University.

\end{document}